\documentclass[a4paper,10pt,twoside]{article}

\usepackage[T1]{fontenc}
\usepackage[utf8]{inputenc}
\usepackage{latexsym}
\usepackage{amsfonts}
\usepackage{amstext}
\usepackage{amsmath}
\usepackage{graphicx}
\usepackage{hyperref}

\usepackage{amssymb,euscript,multicol,graphicx}
\usepackage[ec]{aeguill}       
\usepackage{delarray, fancyhdr}
\usepackage[active]{srcltx}

\newcommand{\be}{\begin{enumerate}}
\newcommand{\ee}{\end{enumerate}}
\newcommand{\beqn}{\begin{eqnarray*}}
\newcommand{\eeqn}{\end{eqnarray*}}
\newcommand{\disp}{\displaystyle}

\newcommand{\incl}[1][r]
      {\ar@<-0.2pc>@{^(-}[#1] \ar@<+0.2pc>@{-}[#1]}

\def\A{{\mathbb A}}
\def\C{{\mathbb C}}

\def\F{{\mathbb F}}

\def\N{{\mathbb N}}
\def\P{{\mathbb P}}
\def\Q{{\mathbb Q}}
\def\R{{\mathbb R}}
\def\Z{{\mathbb Z}}
\def\Ac{{\mathcal A}}
\def\Bc{{\mathcal B}}
\def\Cc{{\mathcal C}}

\def\Ec{{\mathcal E}}
\def\Fc{{\mathcal F}}
\def\Gc{{\mathcal G}}

\def\Lc{{\mathcal L}}
\def\Mc{{\mathcal M}}

\def\Oc{{\mathcal O}}
\def\Pc{{\mathcal P}}

\def\Sc{{\mathcal S}}

\def\Uc{{\mathcal U}}

\def\mgo{{\mathfrak m}}

\def\gm{{\mathbb G}_{\rm{m}}}

\def\eps{\varepsilon}
\def\fleche{\rightarrow}
\def\flechelongue{\longrightarrow}

\def\ov{\overline}
\newcommand{\cartesien}{\ar@{}[dr]|{\square}}

\newcommand{\double}{\ar@<2pt>[r] \ar@<-2pt>[r]}

\def\Ab{\hbox{\rm (Ab)}}
\def\Ann{\hbox{\rm (Ann)}}

\def\Card{\hbox{\rm Card}\,}
\def\Centr{\hbox{\rm Centr}}

\def\Coker{\hbox{\rm Coker}\,}

\def\Ens{\hbox{\rm (Ens)}\,}

\def\Ga{\mathbb{G}_\textrm{a}}

\def\GL{\hbox{\rm GL}\,}

\def\Gm{\mathbb{G}_\textrm{m}}
\def\Gr{\hbox{\rm (Gr)}}

\def\Grass{\hbox{\rm Grass}}
\def\Hom{\hbox{\rm Hom}}
\def\Homs{\hbox{\rm \underline{Hom}}}
\def\id{\hbox{\rm id}}
\def\Id{{\rm Id}}
\def\Im{\hbox{\rm Im}\,}

\def\Ker{\hbox{\rm Ker}\,}

\def\lind{\lim_{\longrightarrow}}

\def\Norm{\hbox{\rm Norm}}

\def\O{\hbox{\rm O}}

\def\Pic{\hbox{\rm Pic}\,}

\def\pr{{\rm pr}}

\def\Proj{\hbox{\rm Proj}\,}
\def\Qcoh{\hbox{${\mathfrak{Qcoh}}$}\,}

\def\rg{\hbox{\rm rg}\,}

\def\SL{\hbox{\rm SL}\,}
\def\Spec{{\rm Spec}\,}

	{\begin{pmatrix}}%
	{\end{pmatrix}}
	{\begin{vmatrix}}%
	{\end{vmatrix}}

	{\ensuremath{\left \{ \hskip -1.5 mm \begin{array}{#1@{\ =\ }l}}}%
	{\end{array}\right.}

	{\ensuremath{\left \{ \hskip -1.5 mm%
	 \begin{array}{#1@{\ }c@{\ }l}}}%
	{\end{array}\right.}

	{\ensuremath{\left \{ \hskip -1.5 mm \begin{array}{#1@{\quad}l}}}%
	{\end{array}\right.}
\newcommand{\fonction}[5]{%
        \ensuremath{#1\colon
        \left\{\hskip -1.5 mm                   
        \begin{array}{c@{\ }c@{\ }l}
        \medskip #2 & \flechelongue & #3 \\
        #4 & \longmapsto & #5 \\
        \end{array}
        \right .
        }}

\newtheorem{souscor}[subsubsection]{Corollaire}

\newtheorem{sousprop}[subsubsection]{Proposition}

\newtheorem{sousthm}[subsubsection]{Th\'eor\`eme}

\newtheorem{sousdefi}[subsubsection]{D\'efinition}

\newtheorem{souslem}[subsubsection]{Lemme}

\newtheorem{sousremarque}[subsubsection]{Remarque}

\newtheorem{sousexo}[subsubsection]{Exercice}

\newtheorem{sousexemple}[subsubsection]{Exemple}

\newtheorem{soustexte}[subsubsection]{\!\!}

\newenvironment{demo}{\noindent {\bf D\'emonstration.}}{$\square$ \vskip .3cm}

\newenvironment{etape}[1]{$\bullet$ \emph{#1}\\}{\vskip .2cm}
\newenvironment{etapefinale}[1]{$\bullet$ \emph{#1}\\}{}

\usepackage[all]{xy}
\entrymodifiers={!!<0pt,0.7ex>+}
\usepackage[francais]{babel}

\title{Topologies de Grothendieck, descente, quotients}
\date{\empty}
\author{Sylvain Brochard (I3M)\\
brochard@math.univ-montp2.fr}

\renewcommand{\to}{\fleche}
\newcommand{\lto}{\flechelongue}
\renewcommand{\fleche}{%
\xymatrix@C=1pc{\ar[r] &}}
\def\flechelongue{%
\xymatrix{\ar[r] &}}
\renewcommand{\mapsto}{%
\xymatrix@C=1pc{\ar@{|->}[r] &}}
\renewcommand{\longmapsto}{%
\xymatrix{\ar@{|->}[r] &}}

\def\Xpl{{X_{\textrm{pl}}}}
\def\Xzar{{X_{\textrm{Zar}}}}
\def\Xet{{X_{\textrm{ét}}}}

\def\Xcl{{X_{\textrm{cl}}}}

\def\FAb{{\mathcal{A}b}}
\def\EspAnn{(Esp. Ann.)}
\def\Div{\textrm{Div}}

\begin{document}

\maketitle
\begin{abstract}
Ces notes sont celles d'un cours donné à Luminy en 2011, dans le cadre d'une école d'été organisée à l'occasion de la réédition de SGA 3. L'objectif est de présenter quelques théorèmes d'existence du quotient d'un schéma par une action de groupe, ou plus généralement par une relation d'équivalence. Nous donnons dans un premier temps les bagages nécessaires de théorie des faisceaux et de théorie de la descente.
\end{abstract}

\tableofcontents

\section*{Introduction}
Si $R$ est une relation d'équivalence sur un espace topologique $X$, on peut munir l'ensemble quotient $X/R$ d'une topologie naturelle, la \og topologie quotient\fg. L'espace topologique quotient ainsi obtenu satisfait à la propriété universelle que l'on attend de lui, et l'on peut dire que la construction de quotients dans la catégorie des espaces topologiques est une trivialité. En géométrie algébrique, c'est une question beaucoup plus délicate. Un point de vue naïf consiste à construire d'abord un espace topologique quotient, puis à essayer d'en faire une variété algébrique. Mais ceci ne marche pas souvent. Les géomètres ont compris depuis longtemps qu'ils ne pouvaient pas séparer les orbites à leur guise, faute de disposer de suffisament de fonctions algébriques.

On peut bien sûr adopter \emph{a priori} une définition catégorique. Si par exemple un schéma en groupes $G$ agit sur un schéma $X$, on appelle \emph{quotient catégorique} de $X$ par $G$ tout schéma $Q$ muni d'un morphisme $G$-invariant $\pi : X \to Q$ qui soit initial parmi les morphismes de schémas ayant cette propriété. Cette définition, si elle garantit l'unicité du quotient, ne résout pas la question de son existence. On peut voir la catégorie des schémas de deux manières différentes, et chacune d'entre elles nous fournit un candidat pour jouer le rôle du schéma quotient. D'un côté, c'est une sous-catégorie de la catégorie des espaces annelés. De l'autre côté, on peut la voir comme sous-catégorie pleine de la catégorie des faisceaux d'ensembles.
$$\shorthandoff{;:!?} 
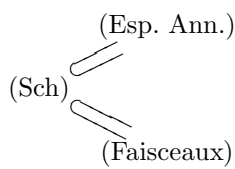
\xymatrix @!0 @R=2pc@C=4pc{& (\textrm{Esp. Ann.})\\
(\textrm{Sch})
\ar@<-0.2pc>@{^(-}[rd]
\ar@<+0.2pc>@{-}[rd]
\ar@<-0.2pc>@{^(-}[ru] 
\ar@<+0.2pc>@{-}[ru]
& \\
&(\textrm{Faisceaux})}$$
Dans chacune de ces deux sur-catégories, toutes les relations d'équivalence admettent des quotients. Nous disposons donc déjà d'un espace annelé quotient, disons $Q_{EA}$, et d'un faisceau quotient, disons $Q_F$. Si l'un de ces objets est un schéma, c'est alors automatiquement un quotient dans la catégorie des schémas. Dans les bons cas, les deux objets sont des schémas (nécessairement isomorphes), ce qui revient alors à dire qu'il existe un quotient catégorique, et que ce quotient a deux bonnes propriétés : il représente le faisceau quotient, et son espace annelé sous-jacent coïncide avec l'espace annelé quotient. Mais il peut arriver que ni $Q_{EA}$ ni $Q_F$ ne soient des schémas. On se retrouve alors avec deux objets quotients, auxquels peut éventuellement s'ajouter un quotient catégorique. Lequel privilégier ? En fait, lorsque $Q_{EA}$ n'est pas un schéma, il n'est pas d'une grande utilité en géométrie algébrique. Le quotient catégorique, quant à lui, peut être très éloigné de l'idée intuitive du quotient que l'on aurait envie de trouver. Par ailleurs, même si $Q_{EA}$ est un schéma, on ne sait en général pas décrire son foncteur des points, et nous aurons tendance à préférer malgré tout\footnote{c'est du moins le parti pris dans ces notes} le faisceau quotient $Q_F$ qui, à défaut d'être un schéma, garde de bonnes chances d'être un espace algébrique, ce qui rend presque les mêmes services en pratique. Dans ce cours, la question principale sera donc : le faisceau quotient est-il représentable par un schéma\footnote{ou par un espace algébrique} ? Subsidiairement, le schéma quotient et le morphisme de passage au quotient ont-ils de bonnes propriétés ?\footnote{Nous n'aborderons pas ici la théorie géométrique des invariants de Mumford, qui s'attelle justement à la construction de \og bons \fg\ quotients, satisfaisant à quelques exigences supplémentaires de nature géométrique. C'est bien sûr regrettable, mais cela aurait nécessité un cours à part entière. Le sujet était d'ailleurs abordé dans d'autres cours lors de l'école d'été.} 

Les faisceaux dont nous parlerons ici ne sont pas des faisceaux pour la topologie de Zariski, mais plutôt des faisceaux pour une \og topologie de Grothendieck\fg, en général \emph{fppf} ou étale. Comme c'est un des points-clés de la construction de quotients que nous allons présenter dans ces notes, nous commencerons par là. Le premier chapitre est donc consacré à quelques rudiments de théorie des faisceaux dans le cadre des topologies de Grothendieck. De même, la théorie de la descente sera fort utile, et nous lui consacrons le deuxième chapitre. Si les textes d'initiation aux topologies de Grothendieck et à la descente ont pu faire défaut par le passé, ce n'est plus le cas aujourd'hui : citons par exemple l'exposé de Vistoli dans~\cite{FGA_explained} et le livre en préparation~\cite{Behrend_Conrad_Edidin_and_Cie}. Nous avons donc essayé de ne pas nous étendre trop longuement sur ces sujets. 

Parmi les partis pris de ce cours, outre bien sûr le choix un peu arbitraire des thèmes abordés (notamment à la fin du chapitre 3), signalons le choix de se limiter, sauf en~\ref{groupoides}, à des actions \emph{libres}. Il se trouve que lorsqu'un groupe agit avec inertie sur un schéma, le faisceau quotient n'est en général pas représentable. Les énoncés traitant ce genre de situation foisonnent pourtant dans \cite{SGA3}~SGA 3, V. Ils assurent généralement, sous diverses hypothèses, que l'espace annelé quotient est un schéma, donnent un certain nombre de propriétés du schéma quotient ainsi obtenu, et montrent enfin que ce dernier représente le faisceau quotient lorsque l'action est libre. Nous avons consenti à la perte de généralité occasionnée par l'option retenue ici pour les deux raisons suivantes :
\begin{itemize}
 \item se limiter aux actions libres simplifie un peu les énoncés et les preuves des théorèmes ;
\item pour une action non libre, on considère aujourd'hui que \og le bon\fg\ objet quotient est plutôt un champ algébrique (le champ quotient), et éventuellement son espace de modules grossier lorsqu'il veut bien exister.
\end{itemize}

Les sujets présentés dans ce cours sont tout à fait classiques en géométrie algébrique et l'auteur de ces lignes ne prétend à aucune originalité, sauf peut-être dans la présentation. Quelques ouvrages ont particulièrement influencé l'écriture de ce texte : bien entendu, le séminaire SGA 3~\cite{SGA3} de Grothendieck, mais aussi~\cite{BRL} et  \cite{FGA_explained}. Ces notes sont celles d'un cours donné à Luminy en septembre 2011 dans le cadre d'une école d'été organisée à l'occasion de la réédition de SGA 3. Je remercie vivement les organisateurs P. Gille et P. Polo de m'avoir invité à donné ce cours. Je remercie aussi chaleureusement M. Raynaud pour son aide précieuse durant sa préparation. Merci enfin à B. Conrad pour la preuve de~\ref{quotient_par_norm}, à T.-Y. Lee pour la preuve détaillée dans l'exercice~\ref{exercice_cas_diag}, et à J. Oesterlé pour des commentaires éclairants. \vskip-2pc \null

\paragraph{Organisation du document.} Comme indiqué plus haut, les deux premiers chapitres sont consacrés aux topologies de Grothendieck et à la théorie de la descente. Nous n'abordons les quotients que dans le troisième. Les sections~\ref{schema_en_groupes} à~\ref{faisceau_quotient} définissent les objets élémentaires dont il sera question et étudient leurs premières propriétés. Nous donnons en~\ref{thm_repres} les principaux résultats de représentabilité du faisceau quotient que l'on peut trouver dans l'exposé~V de SGA~3~\cite{SGA3}, avec des preuves presque complètes. C'est là le cœur technique de ce cours. La section~\ref{espaces_algebriques} présente l'important théorème d'Artin de représentabilité par un espace algébrique, sans preuve. Nous nous spécialisons ensuite en~\ref{cas_ou_la_base_est_un_corps} au cas du quotient d'un groupe par un sous-groupe, lorsque la base est le spectre d'un corps. Des arguments de translation permettent alors d'obtenir un théorème d'existence du quotient (comme schéma) sous des hypothèses de finitude minimalistes. Après ces généralités sur les quotients, les sections~\ref{section_normalisateur} et~\ref{cas_groupe_diago} donnent deux exemples de situations spécifiques dans lesquelles on peut dire plus de choses. Nous effleurons enfin dans la dernière section~\ref{groupoides} le cas des actions non libres, en donnant d'une part un théorème d'algébricité du champ quotient, et d'autre part le théorème de Keel et Mori assurant l'existence d'un espace de modules grossier (ici aussi, sans preuves).

\section{Topologies de Grothendieck}

\subsection{Rappels sur les morphismes lisses, nets ou étales}

Cette section rassemble quelques résultats sur les morphismes lisses, nets ou étales, sans aucune démonstration. Nous renvoyons pour celles-ci à EGA~\cite{EGA} ou \og Néron Models\fg~\cite{BRL}, ou encore au livre~\cite{Milne_Etale_coh} de Milne.

\begin{sousdefi} Soit $f : X\to S$ un morphisme de schémas.
\begin{itemize}
 \item[(i)] On dit que $f$ est formellement lisse (resp. formellement étale, resp. formellement net) si pour tout $S$-schéma affine $T$ et tout sous-schéma fermé $T'$ de $T$ défini par un idéal de carré nul, l'application
$$\Hom_S(T,X) \lto \Hom_S(T',X)$$
est surjective (resp. bijective, resp. injective).
\item[(ii)] On dit que $f$ est lisse (resp. étale, resp. net) s'il est formellement lisse (resp. formellement étale, resp. formellement net) et localement de présentation finie. 
\end{itemize}
\end{sousdefi}

\begin{sousexo}\rm\

\noindent 1) Montrer que $\A_S^n$ est lisse sur $S$.\\
2) Montrer que $\Spec k[x,y]/(xy)$ n'est pas formellement lisse sur $k$ ($k$ un corps).\\
3) Montrer qu'une extension de corps $k'/k$ est finie séparable si et seulement si le morphisme correspondant $\Spec k' \to \Spec k$ est étale. 
\end{sousexo}

Nous listons ci-dessous quelques propriétés de ces morphismes.

\begin{sousprop}[\cite{EGA}~EGA IV 17.1 et 17.3]\ 

 \begin{itemize}
\item[(i)] Une immersion ouverte est étale. Un monomorphisme est formellement net.
\item[(ii)] La classe des morphismes lisses (resp. nets, resp. étales) est stable par composition, par changement de base et par produit.
\item[(iii)] Être lisse (resp. net, resp. étale) est une propriété locale au but comme à la source pour la topologie de Zariski. Autrement dit, si $f : X\to Y$ est un morphisme de schémas, $\{U_i\}$ un recouvrement ouvert de $Y$ et $\{V_{ij}\}$ un recouvrement ouvert de $f^{-1}(U_{i})$, alors $f$ a cette propriété si et seulement si pour tous $i, j$ le morphisme induit $V_{ij}\to U_i$ l'a.
\item[(iv)] Soit $X\to S$ un morphisme net (resp. net, resp. localement de présentation finie). Pour qu'un morphisme $Y\to X$ soit lisse (resp. étale, resp. net) il suffit que le composé $Y\to S$ le soit.  
 \end{itemize}
\end{sousprop}

\begin{sousprop}[\cite{EGA}~EGA IV 17.4.2]
 Soit $f : X\to Y$ un morphisme localement de présentation finie de schémas. Les propriétés suivantes sont équivalentes :
\begin{itemize}
 \item[(i)] $f$ est net ;
\item[(ii)] les fibres de $f$ sont nettes ;
\item[(iii)] le faisceau des différentielles relatives $\Omega_{X/Y}$ est nul ;
\item[(iv)] l'immersion diagonale $\Delta : X \to X\times_Y X$ est ouverte ;
\item[(v)] toute section de $f$ est une immersion ouverte, et ceci reste vrai après tout changement de base $Y'\to Y$. 
\end{itemize}
\end{sousprop}

En particulier, pour savoir si un morphisme est net, il suffit de regarder ses fibres. On est donc ramené à étudier la situation sur un corps. La proposition suivante caractérise les schémas nets sur un corps.

\begin{sousprop}
 Soit $X$ un schéma sur un corps $k$. Les propriétés suivantes sont équivalentes :
\begin{itemize}
 \item[(i)] $X$ est étale sur $k$ ;
\item[(ii)] $X$ est net sur $k$ ;
\item[(iii)] $X$ est une somme disjointe de spectres d'extensions finies séparables de $k$.
\end{itemize}
\end{sousprop}

Maintenant que nous comprenons mieux les morphismes nets, passons aux morphismes lisses. On dit qu'un morphisme $X\to S$ est lisse en $x\in X$ si $x$ admet un voisinage ouvert $U\subset X$ qui est lisse sur $S$ (\emph{idem} pour étale ou net).

\begin{sousprop}[critère jacobien, \cite{BRL}~\S~2.2, prop. 7]
 Soient $A$ un anneau et $B$ une $A$-algèbre de la forme $\frac{A[X_1, \dots, X_n]}{(f_1,\dots, f_r)}$. On note $S=\Spec A$ et $X=\Spec B$. Soit $x$ un point de $X$. Les conditions suivantes sont équivalentes.
\begin{itemize}
 \item[(i)] $X$ est lisse sur $S$ en $x$, de dimension relative $n-r$.
\item[(ii)] Évaluée en $x$, la matrice jacobienne $J(x)=\frac{\partial f}{\partial X}(x)$ est de rang $r$.
\end{itemize}
En particulier, lorsque $n=r$, le morphisme $X\to S$ est étale si et seulement si le déterminant jacobien est inversible dans $B$.
\end{sousprop}

\begin{sousexo}\rm\

\noindent 1) Soient $S=\Spec A$ un schéma affine non vide, $n\in \N^*$ un entier, $a$ un élément de $A$ et $X=\Spec A[T]/(T^n-a)$. Montrer que $X$ est lisse sur $S$ (resp. étale sur $S$) si et seulement si $n$ est inversible dans $A$. En particulier on voit que le groupe $\mu_n$ est étale sur $S$ (ou, ce qui revient au même, lisse) si et seulement si $n$ est inversible dans $S$.\\
 2) Soit $S$ un $\F_p$-schéma affine d'anneau $A$. On pose $B=A[T]/(T^p-T-a)$ pour un $a\in A$. Montrer que $\Spec B$ est étale sur $S$.
\end{sousexo}

\begin{sousprop}[\cite{BRL}~\S 2.2, Proposition~11]
 Soit $f : X\to S$ un morphisme lisse. Alors tout point $x\in X$ admet un voisinage ouvert $U \subset X$ tel que $f_{|_U}$ se factorise en :
$$\xymatrix{f_{|_U} : U \ar[r]^{g} & \A_S^n\ar[r]^p &S
}$$
où $g$ est étale et $p$ est la projection canonique.
\end{sousprop}

\begin{sousprop}[\cite{EGA}~EGA~IV 17.5.1]
 Soit $f : X\to Y$ un morphisme de schémas. Pour que $f$ soit lisse, il faut et il suffit qu'il soit localement de présentation finie, plat et à fibres lisses.
\end{sousprop}

En particulier, on voit qu'un morphisme est étale si et seulement s'il est net et plat. (Il est évident sur la définition que \og étale \fg\ équivaut à net et lisse.)

La propriété suivante dit que, en un point où l'extension résiduelle ne l'interdit pas, un morphisme lisse admet, localement au but pour la topologie étale, toujours une section. Comme il y a beaucoup de points comme ça (il y en a un ensemble dense dans toutes les fibres non vides par~\ref{ensemble_dense}) on voit qu'un morphisme lisse admet des sections passant à peu près partout (après changement de base étale).

\begin{sousprop}[\cite{BRL}~\S 2.2, Proposition 14]
 Soit $f : X\to S$ un morphisme lisse de schémas, $s$ un point de $S$, et $x\in X$ un point au-dessus de $s$ tel que $\kappa(x)$ soit une extension finie séparable de $\kappa(s)$. Alors il existe un morphisme étale $S'\to S$, un point~$s'$ de~$S'$ au-dessus de $s$ de corps résiduel $\kappa(s')=\kappa(x)$ et un $S$-morphisme $S'\to X$ qui envoie~$s'$ sur $x$. (Autrement dit, après le changement de base étale $S'\to S$, le morphisme~$f$ admet une section passant par $x$.)
\end{sousprop}

\begin{sousprop}[\cite{BRL}~\S 2.2, Corollaire 13]
\label{ensemble_dense}
 Si $X$ est un schéma lisse sur un corps~$k$, l'ensemble des points fermés $x$ de $X$ tels que $\kappa(x)$ soit une extension (finie) séparable de~$k$ est dense dans $X$.
\end{sousprop}

\begin{souscor}[\cite{EGA} EGA IV 17.16.3] 
\label{sections_morphismes_lisses}
 Soit $X\to S$ un morphisme lisse et surjectif de schémas. Alors il existe un morphisme $S'\to S$ étale et surjectif, et un $S$-morphisme $S'\to X$. Autrement dit, localement pour la topologie étale, un morphisme lisse et surjectif admet des sections.
\end{souscor}

\subsection{Topologies de Grothendieck}

Avant de définir de façon formelle les topologies de Grothendieck, rappelons ce qu'est un faisceau sur un espace topologique. Si $X$ est un espace topologique, il est possible de définir la notion de faisceau sur $X$ en oubliant complètement l'ensemble sous-jacent à~$X$ pour ne retenir que les ouverts de $X$ et les relations entre eux. On procède de la manière suivante. On note $\Xcl$ la catégorie dont les objets sont les ouverts de $X$ et dont les flèches sont les inclusions $\xymatrix@C=1pc{V \ar@{^(->}[r] &U}$ entre ouverts. Un préfaisceau sur $X$ est alors un foncteur contravariant
$$\Fc : \Xcl^o \flechelongue \Ens.$$
C'est un faisceau si, de plus, pour tout recouvrement ouvert $\{V_i \fleche U\}_{i\in I}$\,, le diagramme
$$\xymatrix{\Fc(U) \ar[r] & \disp \prod_{i\in I} \Fc(V_i) \ar@<2pt>[r]^-{p_1} \ar@<-2pt>[r]_-{p_2} &\disp \prod_{(i,j)\in I^2} \Fc(V_i\cap V_j)}$$
est exact, \emph{i.e.} la première flèche est injective, et son image est l'ensemble des familles $s=(s_i)_{i\in I}$ telles que $p_1(s)=p_2(s)$.
Si l'on souhaite réellement oublier l'ensemble sous-jacent à $X$ et ne retenir que les relations entre les ouverts, on préfèrera éviter le recours à l'intersection ensembliste $V_i\cap V_j$. Il suffit pour cela de remarquer que dans la catégorie $\Xcl$, l'intersection $V_i\cap V_j$ n'est autre que le produit fibré $V_i\times_U V_j$ de $V_i$ et $V_j$ au-dessus de $U$. Cet exemple va nous servir de guide pour les définitions qui suivent.

\begin{sousdefi}
 Soit $\Cc$ une catégorie. Une (pré)topologie sur $\Cc$ est la donnée, pour tout objet $U$ de $\Cc$, d'une collection d'ensembles de flèches $\{V_i \fleche U\}_{i\in I}$ ayant toutes le même but $U$ (un tel ensemble de flèches sera appelé un \og recouvrement \fg\ de $U$ ou encore une \og famille couvrante\fg), ces données vérifiant les conditions suivantes :
\begin{itemize}
 \item[(i)] si $V\fleche U$ est un isomorphisme, alors $\{V \fleche U\}$ est un recouvrement ;
\item[(ii)] si $\{U_i \fleche U\}_i$ est un recouvrement et $V\fleche U$ une flèche, alors les produits fibrés $U_i\times_U V$ existent dans $\Cc$, et la collection $\{U_i \times_U V \fleche V\}_i$ est un recouvrement de $V$ ;
\item[(iii)] si $\{U_i \fleche U\}_i$ est un recouvrement et si pour tout $i$, on se donne un recouvrement $\{V_{ij} \fleche U_i\}_j$  de $U_i$, alors la collection
formée des composés $V_{ij} \fleche U_i \fleche U$ est un recouvrement.
\end{itemize}
\end{sousdefi}

\begin{sousremarque}\rm
\label{rem_voc}
Nous n'avons pas cru bon de préciser les questions relatives aux choix d'univers. Nous renvoyons pour ceci le lecteur à SGA 4 \cite{SGA4}. Par ailleurs, l'objet défini ci-dessus est habituellement appelé \og prétopologie\fg\ (par exemple dans SGA 4), et il y a une notion de topologie induite par une prétopologie. Nous n'utiliserons pas la notion de topologie au sens de SGA 4 : en effet lorsqu'une topologie est induite par une prétopologie, ce qui sera toujours le cas pour nous, les faisceaux peuvent être définis en termes de la prétopologie uniquement. Nous utiliserons dorénavant le terme \og topologie\fg\ pour désigner les \og prétopologies\fg\ de SGA 4. 
\end{sousremarque}

\begin{sousdefi}
 Un \emph{site} est une catégorie munie d'une topologie.
\end{sousdefi}

\begin{sousexemple}\rm Si $X$ est un espace topologique, alors $\Xcl$ est un site (en prenant les recouvrements ouverts comme familles couvrantes).
\end{sousexemple}

\begin{sousexemple}\rm
 Soit $X$ un schéma et soient $C/X$ une sous-catégorie pleine de $(Sch/X)$ et $E$ une classe de morphismes de schémas vérifiant les conditions suivantes :
\begin{itemize}
 \item[(i)] $E$ est stable par changement de base, par composition et contient les isomorphismes ;
\item[(ii)] si $Y \fleche X$ est dans $C/X$ et si $U\fleche Y$ est dans $E$, alors $U\fleche X$ est dans $C/X$.
\end{itemize}
\footnote{Ces conditions impliquent que $C/X$ est partiellement stable par produit fibré, \emph{i.e.} si $U_1\to U$ et $U_2\to U$ sont des flèches dans $C/X$ dont l'une au moins est dans $E$, alors le produit fibré $U_1\times_U U_2$ est encore dans $C/X$.}On note $(C/X)_E$ le site dont la catégorie sous-jacente est $C/X$, et les recouvrements sont les familles $(\xymatrix@C=1pc{U_i \ar[r]^{g_i} &Y})_{i\in I}$
de flèches de $E$ telles que $Y=\cup_{i\in I} g_i(U_i)$.
Par la suite nous serons amené à considérer les cas suivants.
\begin{center}
\begin{tabular}{rcll}
 $E$&=& (Zar) & immersions ouvertes\\
$E$&=& (ét) & morphismes étales\\
$E$&=& (\emph{fppf}) & morphismes plats et localement de présentation finie.
\end{tabular}
\end{center}
En général on parle de \og petit site \fg\ lorsque l'on considère $(E/X)_E$ et de \og grand site\fg\ lorsque l'on considère $(Sch/X)_E$ (ou encore $(Ltf/X)_E$ si besoin, où $(Ltf/X)$ désigne la catégorie des schémas localement de type fini sur~$X$). On notera $\Xzar$, $\Xet$ les petits sites $(\textrm{Zar}/X)_{\textrm{Zar}}$, $(\textrm{Ét}/X)_{\textrm{ét}}$ et $\Xpl$ le grand site $(Sch/X)_{\textrm{fppf}}$.
\end{sousexemple}

\begin{sousexemple}[topologie \emph{fpqc}]\rm
 On dit qu'un morphisme de schémas $X\fleche Y$ est \emph{fpqc} s'il est fidèlement plat et si tout ouvert quasi-compact de $Y$ est l'image d'un ouvert quasi-compact de $X$. En particulier tout morphisme fidèlement plat et quasi-compact est \emph{fpqc}. De même tout morphisme fidèlement plat et universellement ouvert (donc tout morphisme fidèlement plat et localement de présentation finie) est \emph{fpqc}. La \og topologie \emph{fpqc}\fg\ sur $(Sch/S)$ est alors celle pour laquelle les recouvrements sont les familles $\{U_i\fleche U\}_i$ telles que le morphisme induit $\disp\coprod_{i\in I} U_i \fleche U$ soit \emph{fpqc}. Nous renvoyons à~\cite[2.3.2]{FGA_explained} pour les sorites sur les morphismes \emph{fpqc} et pour la preuve du fait que la topologie \emph{fpqc} que nous venons de définir est bien une topologie.
\end{sousexemple}

\begin{sousdefi}
 Soit $\Cc$ un site. Un préfaisceau sur $\Cc$ est un foncteur contravariant
$$\Fc : \Cc^o \flechelongue \Ens\, .$$
On dit que $\Fc$ est un \emph{faisceau} si pour toute famille couvrante $\{V_i \fleche U\}_{i\in I}$ le diagramme
$$\xymatrix{\disp \Fc(U) \ar[r] &\disp \prod_{i\in I} \Fc(V_i) \ar@<2pt>[r]^-{p_1} \ar@<-2pt>[r]_-{p_2} &\disp \prod_{(i,j)\in I^2} \Fc(V_i\times_U V_j)}$$
est exact. Nous noterons $\Pc(\Cc)$ la catégorie des préfaisceaux sur $\Cc$ et $\Sc(\Cc)$ la catégorie des faisceaux sur $\Cc$.
\end{sousdefi}

\begin{sousremarque}\rm On peut se donner des structures supplémentaires sur $\Fc$. Par exemple un \og faisceau de groupes\fg\ (resp. de groupes abéliens, d'anneaux) est un foncteur contravariant de $\Cc$ dans la catégorie $\Gr$ des groupes (resp. dans celle $\Ab$ des groupes abéliens, dans celle $\Ann$ des anneaux commutatifs et unitaires) dont la composition avec le foncteur d'oubli de $\Gr$ (resp. de $\Ab$, $\Ann$) vers $\Ens$ est un faisceau d'ensembles. De même si $\Ac$ est un faisceau d'anneaux sur $\Cc$, on dispose d'une notion naturelle de $\Ac$-module. Nous noterons $\FAb(\Cc)$ la catégorie des faisceaux de groupes abéliens sur $\Cc$.
\end{sousremarque}

\begin{sousexemple}[préfaisceaux constants]\rm
\label{prefaisceau_constant}
 Soient $A$ un ensemble et $S$ un schéma. On note $P_A$ le foncteur de $(Sch/S)^o$ vers $\Ens$ défini par $P_A(U)=A$ pour $U$ non vide et $P_A(U)=\{*\}$ sinon (avec les flèches de transition évidentes). Ceci définit un préfaisceau, mais $P_A$ n'est en général pas un faisceau. Si $A$ est un groupe (ou un anneau, \emph{etc.}) alors $P_A$ est clairement un préfaisceau de groupes (ou d'anneaux, \emph{etc.}). 
\end{sousexemple}

\begin{sousexo}\rm
 On suppose $\Card A\geq 2$ et $X$ non vide. Montrer que la restriction de $P_A$ au petit site Zariski $\Xzar$ est un faisceau si et seulement si $X$ est irréductible.
\end{sousexo}

\begin{sousexo}\rm
\label{foncteur_de_Picard_absolu}
 Soit $f : X\fleche S$ un morphisme de schémas. Montrer que le foncteur en groupes sur $(Sch/S)$ qui à $T$ associe $\Pic(X\times_S T)$ n'est jamais un faisceau pour la topologie de Zariski, sauf si $X$ est vide. [Prendre un $S$-schéma $T$ muni d'un faisceau inversible $\Lc$ tel que $f_T^* \Lc$ soit non trivial. Conclure en remarquant que, localement pour la topologie de Zariski sur $T$, le fibré $\Lc$ est trivial, donc $f_T^*\Lc$ aussi. Par exemple on pourra prendre $T=\P_X^1$ et $\Lc=\Oc(1)$. Alors le morphisme diagonal $X\fleche X\times_S X$ induit une section de $f_T$, ce qui prouve que $f_T^*\Lc$ est non trivial.]
\end{sousexo}

Le théorème suivant est un résultat de descente. Nous verrons plus précisément que c'est une conséquence immédiate du théorème~\ref{descente_fpqc_gen}. Il permet de produire de nombreux exemples de faisceaux. Il dit que, vu comme foncteur contravariant (donc comme préfaisceau) sur la catégorie $(Sch/S)$, tout $S$-schéma est un faisceau \emph{fpqc}. En particulier tout $S$-schéma en groupes est un faisceau de groupes.

\begin{sousthm}
\label{schema_implique_faisceau}
 Soit $X\fleche S$ un morphisme de schémas. Alors le foncteur
$$\fonction{h_X}{(Sch/S)^o}{\Ens}{U}{\Hom_S(U,X)}$$
est un faisceau pour la topologie \emph{fpqc} (donc aussi \emph{fppf}, étale et Zariski).
\end{sousthm}

On dit qu'un préfaisceau $P$ sur $(Sch/S)$ est \emph{représentable} s'il existe un $S$-schéma $X$ tel que $P$ soit isomorphe à $h_X$.

\begin{sousexemple}\rm\ 

\begin{itemize}
 \item[$\bullet$] Le préfaisceau $\Ga$ défini par $\Ga(U)=\Gamma(U,\Oc_U)$ est un faisceau de groupes abéliens sur $\Spec \Z$. En effet il est représentable par $\Spec \Z[X]$.
 \item[$\bullet$] Le préfaisceau $\Gm$ défini par $\Gm(U)=\Gamma(U,\Oc_U)^{\times}$ est un faisceau de groupes abéliens sur $\Spec \Z$. En effet il est représentable par $\Spec \Z[X,X^{-1}]$.
 \item[$\bullet$] Le préfaisceau $\mu_n$ défini par $\mu_n(U)=\{x\in \Gamma(U,\Oc_U)^{\times}\ |\ x^n=1\}$ est un faisceau de groupes abéliens sur $\Spec \Z$. En effet il est représentable par $\Spec\left( \frac{\Z[X]}{X^n-1}\right)$.
\end{itemize}

\noindent
Pour un schéma quelconque $S$, on notera ${\Gm}_{,S}$ (resp. ${\Ga}_{,S}$, $\mu_{n,S}$) le $S$-schéma en groupes $\Gm\times S$ (resp. $\Ga\times S$, $\mu_n\times S$).
\end{sousexemple}

\subsection{Comparaison entre les topologies \emph{fpqc}, \emph{fppf}, étale et Zariski}

Ce paragraphe est une digression sur les topologies introduites plus haut sur la catégorie des schémas et peut être omis en première lecture. Une famille couvrante au sens de Zariski (resp. étale, \emph{fppf}) est toujours une famille couvrante au sens étale (resp. \emph{fppf}, \emph{fpqc}). En particulier on voit que les topologies Zariski, étale, \emph{fppf} et \emph{fpqc} sont \og de plus en plus fines\fg. En fait chacune est \emph{strictement} plus fine que la précédente. Ainsi, un revêtement fini étale est un recouvrement étale mais n'est en général pas un recouvrement Zariski (par exemple $\Spec \C$ sur $\Spec \R$, ou bien le revêtement étale non trivial de degré~2 de la cubique nodale dessiné dans l'exercice~III~10.6 de~\cite{Hartshorne_AG}). Nous précisons un peu ces questions de comparaison dans les exercices suivants.

\begin{sousexo}[comparaison des morphismes couvrants]\rm
\ 

\label{def_couvrants}
 \noindent On dit qu'un morphisme $X\to Y$ est couvrant pour une topologie s'il existe des morphismes $U_i\to X$ tels que les composés $U_i\to Y$ forment une famille couvrante (\emph{i.e.} un recouvrement) de $Y$ pour la topologie considérée. Il est évident que si $X\to Y$ est un recouvrement alors c'est un morphisme couvrant, mais la réciproque est fausse en général. Par exemple les morphismes lisses et surjectifs sont couvrants pour la topologie étale d'après~\ref{sections_morphismes_lisses} mais en général ils ne sont pas eux-mêmes étales. Le fait que les topologies lisse et étale aient les mêmes morphismes couvrants entraîne que les catégories de faisceaux coïncident (c'est d'ailleurs pour cette raison que l'on ne parle pas beaucoup de la topologie lisse). Pour montrer qu'une topologie est réellement plus fine qu'une autre, il faut montrer qu'elles n'ont pas les mêmes morphismes couvrants. (On voit ici l'un des avantages des \og topologies \fg\ de SGA~4\dots)\\
a) Montrer que les recouvrements étales donnés en exemple à la fin du paragraphe précédent ne sont pas des morphismes couvrants pour la topologie de Zariski.\\
b) Soit $k$ un corps. Montrer que le morphisme $\A^1_k\to \A_k^1$ donné par $x\mapsto x^2$ est une famille couvrante pour la topologie \emph{fppf} mais n'est pas un morphisme couvrant pour la topologie étale. Même question pour une extension de corps finie et purement inséparable.\\
c) Une extension de corps est une famille couvrante pour la topologie \emph{fpqc} mais n'est pas un morphisme couvrant pour la topologie \emph{fppf} sauf si l'extension est finie. Voici un autre exemple intéressant de famille couvrante \emph{fpqc}. Soit $X$ un schéma noethérien contenant au moins deux points et soit $p$ un point fermé de $X$. On note $U$ le complémentaire de $p$ et $X_p$ le spectre de l'anneau local en $p$. Alors le morphisme naturel $U\coprod X_p \to X$ est un recouvrement \emph{fpqc} mais en général pas ce n'est pas un morphisme couvrant pour la topologie \emph{fppf} (le démontrer).
\end{sousexo}

\begin{sousexo}[comparaison des faisceaux]\rm Pour comparer les différentes topologies, on peut aussi comparer les catégories de faisceaux qu'elles définissent.\\
a) Montrer que tout faisceau \emph{fpqc} est un faisceau \emph{fppf}. De même, tout faisceau \emph{fppf} est un faisceau étale et tout faisceau étale est un faisceau Zariski.\\
b) Soient $S=\Spec \Q$ et $f : {\Gm}_{,S} \to {\Gm}_{,S}$ l'élévation à la puissance $n$. Montrer que le faisceau image de $f$ au sens de Zariski n'est pas un faisceau étale. (Voir la remarque~\ref{Ker_Im_Coker} pour la définition du faisceau image.)\\
c) Soient $p$ un nombre premier, $S=\Spec \F_p$ et $f : {\Gm}_{,S} \to {\Gm}_{,S}$ l'élévation à la puissance $p$. Montrer que le faisceau image de $f$ au sens étale n'est pas un faisceau \emph{fppf}.\\
d) Soit $k$ un corps. Montrer que le foncteur de Picard $\Pic_{\A_k^1/k}$ de la droite affine sur~$k$, défini en~\ref{foncteur_de_Picard}, est un faisceau \emph{fppf} mais n'est pas un faisceau \emph{fpqc}. (Ceci prouve en particulier qu'il n'est pas représentable.)
\end{sousexo}

\subsection{Faisceau associé à un préfaisceau}

\begin{sousthm}
\label{faisceau_associe}
 Soient $\Uc$ un univers, $\Cc$ un $\Uc$-site et $F : \Cc^o \fleche \Ens$ un préfaisceau. Il existe un faisceau $aF$ et un morphisme de préfaisceaux $\phi : F \fleche aF$ avec la propriété universelle suivante :  pour tout faisceau $G$ et tout morphisme de préfaisceaux $\phi'$ de $F$ vers~$G$, il existe un unique morphisme $\psi : aF \fleche G$ tel que $\psi \circ \phi = \phi'$.
$$\xymatrix@C=1.1pc{F \ar[rr]^{\phi} \ar[rd]_{\phi'} && aF \ar@{.>}[ld]^{\psi} \\&G& }$$
\end{sousthm}
\begin{demo}
 \cite{SGA4} SGA 4, II, théorème~3.4.
\end{demo}
Autrement dit, le foncteur d'oubli $\textrm{oub} : \Sc(\Cc) \fleche \Pc(\Cc)$ admet un adjoint à gauche $a : \Pc(\Cc) \fleche \Sc(\Cc)$. Le faisceau associé $aF$ est caractérisé (à unique isomorphisme près) par les deux propriétés suivantes :
\begin{itemize}
 \item[(i)] Deux sections de $F$ sont égales dans $aF$ si et seulement si elles sont \og localement égales\fg\ dans $F$. Plus précisément : soient $U$ un objet de $\Cc$ et $\xi, \eta\in F(U)$. Alors, pour que l'égalité $a\xi=a\eta$ ait lieu dans $(aF)(U)$, il faut et il suffit qu'il existe un recouvrement $\{U_i\fleche U\}_i$ tel que pour tout $i\in I$ on ait $\xi_{|_{U_i}}=\eta_{|_{U_i}}$.
\item[(ii)] Toute section de $aF$ provient localement d'une section de $F$. Plus précisément : soient $U$ un objet de $\Cc$ et $\ov{\xi}\in (aF)(U)$. Alors il existe un recouvrement $\{U_i \fleche U\}_i$ et des $\xi_i\in F(U_i)$ tels que pour tout $i$ on ait $a\xi_i=\ov{\xi}_{|_{U_i}}$.
\end{itemize}

\begin{sousremarque}\rm Le foncteur $a$ dépend de la topologie considérée sur $\Cc$. Si $F$ est déjà un faisceau, alors le morphisme $F\fleche aF$ est un isomorphisme.
\end{sousremarque}

{\footnotesize
\begin{sousremarque}[questions d'univers]\rm
Contrairement à ce qui était annoncé dans la remarque~\ref{rem_voc}, nous avons cru bon ici de préciser que le site $\Cc$ en question devait être un \og $\Uc$-site\fg, où $\Uc$ est un univers fixé. (Nous renvoyons le lecteur désireux de connaître les définitions d'un $\Uc$-site ou d'un univers à~SGA~4~\cite{SGA4}.) En fait, la preuve du théorème repose sur l'existence de certaines limites inductives indexées par la collection des familles couvrantes d'un objet donné $U$ de $\Cc$. Ces limites inductives existent pourvu que cette collection de familles couvrantes ne soit pas trop grosse. C'est essentiellement ce que garantit l'hypothèse selon laquelle $\Cc$ est un $\Uc$-site. Lorsque l'on travaille avec les topologies Zariski, étale ou \emph{fppf}, on peut ignorer ces problèmes logiques : le faisceau associé existe et ne dépend pas du choix de l'univers. En revanche avec la topologie \emph{fpqc} il faut prendre quelques précautions : le faisceau associé \emph{dépend a priori} de l'univers fixé. En un sens, on peut dire qu'en général, \emph{le faisceau fpqc associé n'existe pas !} Waterhouse donne d'ailleurs dans l'article~\cite{Waterhouse_Basically_bounded_functors}, théorème~5.5, un exemple explicite de préfaisceau qui n'admet pas de faisceautisé. C'est pourquoi on doit généralement se contenter de la topologie \emph{fppf} dans les processus de faisceautisation. Toutefois, pour un préfaisceau $F$ donné et un morphisme $\phi : F \to G$ vers un faisceau $G$, on peut parfaitement donner un sens à la phrase \og $G$ est le faisceau \emph{fpqc} associé à $F$\fg. Il suffit de demander à~$G$ de satisfaire à la propriété universelle du théorème~\ref{faisceau_associe}, ou ce qui revient au même, de vérifier les propriétés~(i) et~(ii) mentionnées ci-dessus. Un tel $G$ est évidemment unique à unique isomorphisme près. Signalons enfin que, si l'on souhaite réellement disposer, pour un préfaisceau $F$ donné, d'un faisceau \emph{fpqc} associé, c'est parfois possible en pratique. En effet, Waterhouse définit dans~\cite{Waterhouse_Basically_bounded_functors} une classe de préfaisceaux, qu'il qualifie de \og basically bounded\fg, et montre que dans le cas de ces préfaisceaux, le faisceau associé existe bel et bien.
\end{sousremarque}
}

\begin{sousexemple}[faisceau constant]\rm
\label{faisceau_constant}
 Sur $(Sch/S)_{\textrm{Zar}}$ (pour un schéma $S$ donné), on note $\underline{A}$ le faisceau associé au préfaisceau $P_A$ de l'exemple~\ref{prefaisceau_constant}. On montre facilement que~$\underline{A}$ s'identifie au faisceau qui à $U$ associe l'ensemble des fonctions continues de $U$ dans $A$ (où $A$ est muni de la topologie discrète). Si $U$ est localement connexe alors  $\underline{A}(U)=A^{\pi_0(U)}$. Il est inutile de faisceautiser davantage. En effet, le faisceau $\underline{A}$ ainsi construit est représentable, d'après~\ref{schema_implique_faisceau} c'est donc un faisceau pour toutes les topologies raisonnables sur la catégorie des schémas (en particulier \emph{fpqc}). Comme schéma, $\underline{A}$ est tout simplement isomorphe à une union disjointe de copies de $S$, indexée par les éléments de A. Maintenant si $A$ est un groupe alors $\underline{A}$ est naturellement un schéma en groupes. Le morphisme structural $\mu : \underline{A}\times_S \underline{A}\fleche \underline{A}$ peut être décrit de la manière suivante. On note $\underline{A}=\disp\coprod_{a\in A} S_a$ avec $S_a=S$. Le produit fibré $\underline{A}\times_S \underline{A}$ s'identifie canoniquement à $\disp \coprod_{a,b\in A} S_a\times_S S_b$. Le morphisme $\mu$ envoie simplement $S_a\times_S S_b$ (canoniquement isomorphe à $S$) sur la copie numéro $ab$ de $S$ dans $\underline{A}$.
\end{sousexemple}

Le processus de \og faisceautisation\fg\ est crucial dans un certain nombre de constructions géométriques. Ce sera le cas pour la construction de quotients par des actions de groupes, où nous commencerons par définir un faisceau quotient, dont nous étudierons ensuite la représentabilité. Voici un autre exemple fondamental de faisceautisation.

\begin{sousexemple}\rm
\label{foncteur_de_Picard}
 Soit $f : X\fleche S$ un morphisme de schémas. On note $\Pic_{X/S}$ le faisceau \emph{fppf} associé au préfaisceau $T\mapsto \Pic(X\times_S T)$ de l'exercice~\ref{foncteur_de_Picard_absolu}. Ce faisceau est appelé le foncteur de Picard relatif de $X/S$. S'il est représentable par un schéma, le schéma qui le représente est un schéma en groupes appelé le \emph{schéma de Picard} de $X/S$. C'est l'objet géométrique naturel qui paramètre les fibrés en droites sur $X$.
\end{sousexemple}

\subsection{Fonctorialité}

Pour un schéma $X$, on note $\Cc_X$ l'un des sites (gros ou petit) Zariski, étale, ou \emph{fppf}, et $\Pc(\Cc_X)$ (resp. $\Sc(\Cc_X)$) la catégorie des préfaisceaux (resp. des faisceaux) sur ce site. Un morphisme de schémas $f : X\to Y$ induit alors un foncteur \og image directe\fg\ $f_*$ de~$\Pc(\Cc_X)$ vers $\Pc(\Cc_Y)$ défini comme suit. Soit $F$ un préfaisceau sur $X$. Pour tout $U\in \Cc_Y$, on pose $$(f_*F)(U)=F(X\times_Y U).$$ On vérifie facilement que, si $F$ est un faisceau, alors $f_*F$ est encore un faisceau, donc $f_*$ induit un foncteur, encore noté $f_*$, de $\Sc(\Cc_X)$ dans $\Sc(\Cc_Y)$. Ce foncteur commute clairement aux limites projectives, et il a un adjoint à gauche $f^*$ que l'on peut décrire de la manière suivante. Soient $G$ un faisceau sur $Y$. Alors $f^*G$ est le faisceau associé au préfaisceau qui à $U\in \Cc_X$ associe $\disp \lind F(V)$, où la limite inductive est prise sur les carrés commutatifs
$$\xymatrix@R=1pc@C=1pc{U\ar[r]\ar[d] &V\ar[d] \\ X\ar[r]& Y}$$
avec $V\in \Cc_Y$.

Nous laissons au lecteur le soin de vérifier que le foncteur $f^*$ ainsi défini est bien l'adjoint à gauche de $f_*$. En particulier il commute aux limite inductives. De plus, comme les sites en question ont des objets finals et des produits fibrés, on vérifie que $f^*$ est exact (voir par exemple~\cite{Tamme_Etale_cohomology}~3.6.7). Enfin, pour des morphismes composables $f : X\to Y$ et $g : Y\to Z$, on voit immédiatement sur les définitions que l'on a des isomorphismes canoniques $g_*f_* \simeq (gf)_*$ et $f^*g^*\simeq (gf)^*$. Nous utiliserons peu ces foncteurs par la suite, le lecteur est donc invité à consulter les références classiques pour plus d'informations. 

\subsection{La catégorie des faisceaux abéliens}

\begin{sousprop} Soit $\Cc$ un site. La catégorie $\FAb(\Cc)$ des faisceaux de groupes abéliens sur $\Cc$ est une catégorie abélienne avec suffisament d'injectifs.
\end{sousprop}
\begin{demo}
Voir par exemple \cite{Tamme_Etale_cohomology}~3.2.2.
\end{demo}

\begin{sousremarque}\rm
\label{Ker_Im_Coker}
 Soit $f : F\to G$ un morphisme dans $\FAb(\Cc)$. On peut décrire le noyau, le conoyau et l'image de $f$ de la manière suivante. Le noyau est simplement donné par $(\Ker f)(U)=\Ker(f(U) : F(U)\to G(U))$. Le conoyau est $\Coker f= aC$, le faisceau associé au préfaisceau $C$ donné par $C(U)=\Coker(f(U))$. De même le faisceau image $\Im f$ est le faisceau associé au préfaisceau image $I$ défini par $I(U)=\Im(f(U))$.

L'objet nul est ce que l'on pense, de même que les produits quelconques de faisceaux. Plus généralement les limites projectives de faisceaux se calculent terme à terme. En revanche, pour les sommes directes (et plus généralement pour les limites inductives), après avoir fait le calcul terme à terme il faut prendre le faisceau associé.
\end{sousremarque}

\begin{sousexo}\rm
 Soit $A\!\to\! B\!\to\! C$ une suite de faisceaux abéliens \emph{fppf} sur un schéma~$X$. Montrer que, si cette suite est exacte dans $\FAb(\Xet)$, alors elle est aussi exacte dans~$\FAb(\Xpl)$.
\end{sousexo}

\begin{sousexemple}[suite exacte de Kummer]\rm
 Soit $X$ un schéma. Nous allons montrer que la suite
$$\xymatrix{0\ar[r] &\mu_n \ar[r]& \gm \ar[r]^n &\gm \ar[r] & 0}$$
est exacte dans la catégorie $\FAb(\Xpl)$, et que cette suite est même exacte dans la catégorie $\FAb(\Xet)$ si aucune caractéristique résiduelle de $X$ ne divise $n$.
D'abord, quel que soit le site avec lequel on travaille, étale ou \emph{fppf}, il est clair que $\mu_n$ est le noyau de l'élévation à la puissance $n$ de $\gm$ dans lui-même. Il faut donc seulement montrer l'exactitude à droite. Pour cela, il suffit de montrer que, si $\xi\in\gm(U)$ pour un certain $X$-schéma $U$, alors $\xi$ est localement (pour la topologie considérée) une puissance $n$-ième. La question étant locale sur $U$, on peut supposer $U$ affine, disons $U=\Spec A$. On pose $B=\frac{A[T]}{T^n-\xi}$ et $U'=\Spec B$. Alors $U'\to U$ est une famille couvrante \emph{fppf} (la platitude est évidente car $B$ est libre sur~$A$) et, par construction, la restriction $\xi_{|_{U'}}$ de $\xi$ à $U'$ admet une racine $n$-ième. Ceci prouve la première assertion. Enfin, si $n$ est inversible dans $\Gamma(X,\Oc_X)$, le critère jacobien montre immédiatement que $B$ est étale sur $A$. Le morphisme $U'\to U$ est donc un recouvrement étale et ceci prouve la seconde assertion.
\end{sousexemple}

\begin{sousexo}\rm
 Montrer que pour un schéma $X$ et un entier $n$, on a équivalence entre :
\begin{itemize}
 \item[a)] aucune caractéristique résiduelle de $X$ ne divise $n$ ;
\item[b)] $n$ est inversible dans $\Gamma(X,\Oc_X)$. 
\end{itemize}
\end{sousexo}

\begin{sousexo}[suite exacte d'Artin-Schreier]\rm
 Soit $X$ un $\F_p$-schéma. Le morphisme
$F-\Id : {\Ga}_{,X} \to {\Ga}_{,X}$
donné par $x\mapsto x^p-x$ est alors un morphisme de groupes.

\noindent 1) Soit $U$ un $X$-schéma et soit $x\in \Ker((F-\Id)(U))$. Montrer qu'il existe des sous-schémas ouverts et fermés $F_0, \dots, F_{p-1}$ de $U$ tels que $x=i$
 dans $\Oc_U(F_i)$.

\noindent 2) En déduire que le noyau de $F-\Id$ est le faisceau constant $\Z/p\Z$. 

\noindent 3) Montrer que la suite
$$\xymatrix{0 \ar[r]& \Z/p\Z \ar[r] &{\Ga}_{,X}\ar[r]^{F-\Id} &{\Ga}_{,X}\ar[r] &0 }$$
est exacte pour la topologie étale (donc aussi pour la topologie \emph{fppf}).
\end{sousexo}

\begin{sousexo}\rm
 Soit $X$ un $\F_p$-schéma et soit $\alpha_p$ le sous-groupe de ${\Ga}_{,X}$ défini fonctoriellement par $\alpha_p(U)=\{x\in \Gamma(U,\Oc_U)\ |\ x^p=0\}$. Il est représentable, c'est donc un sous-schéma en groupes de ${\Ga}_{,X}$. Si $F$ désigne le Frobénius $x\mapsto x^p$, montrer que la suite de faisceaux \emph{fppf}
$$\xymatrix{0 \ar[r]& \alpha_p \ar[r] &{\Ga}_{,X}\ar[r]^{F} &{\Ga}_{,X}\ar[r] &0 }$$
est exacte dans la catégorie $\FAb(\Xpl)$. 
\end{sousexo}

Soit $\Cc$ un site. Comme la catégorie $\FAb(\Cc)$ des faisceaux de groupes abéliens sur $\Cc$ a suffisament d'injectifs, on peut définir des foncteurs dérivés à droite $R^qf$, $q\geq 0$, pour tout foncteur $f : \FAb(\Cc) \to \Bc$ additif et exact à gauche à valeurs dans une catégorie abélienne~$\Bc$. Sur un schéma $X$, on peut donc définir les groupes de cohomologie étale (resp. \emph{fppf}) en dérivant le foncteur \og sections globales\fg\ sur la catégorie des faisceaux abéliens étales (resp. \emph{fppf}). De même si $f : X\to Y$ est un morphisme de schémas, et si~$\Fc$ est un faisceau (étale, \emph{fppf},  \emph{etc.}) on peut définir ses \og images directes supérieures\fg~$R^qf_*\Fc$ pour la topologie considérée.

\begin{sousexemple}\rm
 La suite exacte de Kummer ci-dessus induit une suite exacte longue de cohomologie \emph{fppf} :
$$\xymatrix @R=1mm{0\ar[r]& H^0(\Xpl,\mu_n) \ar[r]& H^0(\Xpl,\Gm)\ar[r]& H^0(\Xpl,\Gm)   \\
\ar[r]& H^1(\Xpl,\mu_n)\ar[r]& H^1(\Xpl,\Gm)\ar[r]& H^1(\Xpl,\Gm) \\
\ar[r]& H^2(\Xpl,\mu_n)\ar[r]&\parbox{1.9cm}{\quad \dots}
}$$
\end{sousexemple}

\begin{sousremarque}\rm
 Soit $X$ un schéma. Nous avons vu ci-dessus que la catégorie des faisceaux de groupes abéliens (disons au sens \emph{fppf} pour fixer les idées) est abélienne. En général la sous-catégorie des faisceaux représentables, c'est-à-dire la catégorie des $X$-schémas en groupes commutatifs, n'est pas abélienne. Cependant, nous verrons plus loin que la catégorie des schémas en groupes commutatifs et de type fini sur le spectre d'un corps est abélienne (\ref{categories_abeliennes}).
\end{sousremarque}

\subsection{Anneaux locaux (pour la topologie étale)}

Nous n'aurons pas besoin dans ce cours d'anneaux locaux pour la topologie étale, mais il nous a semblé difficile de passer totalement sous silence leur existence. Le lecteur pressé d'arriver aux quotients de schémas peut donc allègrement sauter cette section.
Pour la topologie étale, on peut définir le germe d'un faisceau en un point de manière tout à fait analogue à ce que l'on fait pour la topologie de Zariski, et l'objet ainsi construit rend les mêmes services.

Soient $X$ un schéma et $F$ un faisceau étale abélien sur $X$. Soit $\ov{x} : \Spec \Omega \to X$ un point géométrique de $X$ (\emph{i.e.} $\Omega$ est un corps séparablement clos). Un \emph{voisinage étale} de $\ov{x}$ est un couple $(U,u)$ où $U$ est un $X$-schéma étale et $u$ est un $\Omega$-point de $U$ qui relève $\ov{x}$.
$$\xymatrix{&U\ar[d]^{\textrm{étale}} \\
\Spec \Omega \ar[r]_-{\ov{x}} \ar[ur]^u &X
}$$
Un morphisme de $(U,u)$ vers un autre voisinage étale $(V,v)$ est un $X$-morphisme $\varphi : U\to V$ tel que $\varphi\circ u=v$. Les voisinages étales de $\ov{x}$ forment ainsi une catégorie cofiltrante au sens SGA~4~\cite[I~2.7]{SGA4}.

\begin{sousdefi}
 Le germe de $F$ en $\ov{x}$ est
$$F_{\ov{x}}=\lind F(U)$$
où la limite inductive est prise sur la catégorie des voisinages étales de $\ov{x}$. 
\end{sousdefi}

Nous listons ci-dessous quelques propriétés utiles, sans preuve.

\begin{sousprop}\ 

 \begin{itemize}
  \item[(i)] Cette définition ne dépend pas du choix du corps séparablement clos $\Omega$, mais seulement du point $x\in X$ image de $\ov{x}$.\footnote{Plus précisément, si $L$ est une extension séparablement close de $\Omega$ et si $\ov{x}_L$ est le composé $\Spec L \to \Spec \Omega \to X$, le morphisme naturel $F_{\ov{x}}\to F_{\ov{x}_L}$ est un isomorphisme.} On note parfois simplement $F_x$ au lieu de $F_{\ov{x}}$ le germe de $F$ en $x$ si aucune confusion n'en résulte.
\item[(ii)] Un morphisme $v : F' \to F$ de faisceaux abéliens sur $\Xet$ est un isomorphisme (resp. un monomorphisme, un épimorphisme) si et seulement si pour tout $x\in X$, le morphisme induit $F'_{x}\to F_{x}$ en est un.
\item[(iii)] Un morphisme $v$ comme ci-dessus est nul si et seulement si pour tout $x\in X$, le morphisme induit $F'_{x}\to F_{x}$ est nul.
\item[(iv)] Pour une section $s\in F(X)$, on a $s=0$ si et seulement si son image dans $F_{x}$ est nulle pour tout $x\in X$.
\item[(v)] Une suite $F'\to F\to F''$ est exacte si et seulement si pour tout $x\in X$ la suite $F'_{x}\to F_{x}\to F''_{x}$ est exacte.
\item[(vi)] Pour $x\in X$ fixé, le foncteur qui à $F$ associe $F_{x}$ commute aux limites inductives.
\item[(vii)] Si $f : Y\to X$ est un morphisme de schémas et si $y$ est un point de $Y$, alors
$$(f^*F)_{y}\simeq F_{f(y)}.$$
 \end{itemize}
\end{sousprop}

\begin{sousremarque}\rm
 Soit $\Omega$ un corps séparablement clos et $p=\Spec \Omega$. On rappelle que $\FAb(p_{\textrm{ét}})$ désigne la catégorie des faisceaux de groupes abéliens sur le petit site étale de $p$. Le foncteur naturel
$$
\left\{
\begin{array}{ccc}
           {\FAb(p_{\textrm{ét}})} & \lto& {\Ab} \\{F} & \longmapsto& {F(p)}
          \end{array}
\right.
$$
est une équivalence de catégories. Maintenant si $F$ est un faisceau étale abélien sur un schéma $X$ quelconque et si $\ov{x} : p \to X$ est un point géométrique, le groupe abélien $F_{\ov{x}}$ correspond \emph{via} l'équivalence ci-dessus au faisceau étale abélien $\ov{x}^*F$ sur $\Spec \Omega$.
\end{sousremarque}

\begin{sousdefi}
 Soient $X$ un schéma et $x$ un point de $X$. On appelle \emph{hensélisé strict de $X$ en $x$}, et on note $\Oc_{X,\ov{x}}$ ou parfois $\Oc_{X,x}^{sh}$, le germe du faisceau structural $\Oc_X$ au point géométrique $\ov{x} : \Spec \kappa(x)^s \to X$ où $\kappa(x)^s$ est une clôture séparable du corps résiduel de $x$.\footnote{On évite ici la notation $\Oc_{X,x}$ pour ne pas confondre avec l'anneau local en $x$ au sens de la topologie de Zariski.} 
\end{sousdefi}

On a alors un morphisme naturel
$$\Spec \Oc_{X,\ov{x}} \lto X.$$

\begin{sousdefi}
 Un anneau local $A$ est dit strictement hensélien si le morphisme naturel de $A$ dans $\Oc_{X,\ov{x}}$ est un isomorphisme (avec $X=\Spec A$ et $x$ son point fermé). Un schéma est dit \emph{strictement local} si c'est le spectre d'un anneau local strictement hensélien.
\end{sousdefi}

\begin{sousprop}\ 

\begin{itemize}
 \item[(i)] $\Oc_{X,\ov{x}}$ est strictement hensélien, de corps résiduel la clôture séparable de $\kappa(x)$ dans~$\Omega$.
\item[(ii)] Le morphisme naturel $\Oc_{X,x}\to \Oc_{X,\ov{x}}$ est fidèlement plat. Par ailleurs si $\mgo_x$ (resp.~$\mgo_{\ov{x}}$) désigne l'idéal maximal de $\Oc_{X,x}$ (resp. de $\Oc_{X,\ov{x}}$), alors
$$\mgo_{\ov{x}}=\mgo_x.\Oc_{X,\ov{x}}.$$
\item[(iii)] Pour que $\Oc_{X,\ov{x}}$ soit noethérien, il faut et il suffit que $\Oc_{X,x}$ le soit.
\end{itemize}
\end{sousprop}

\begin{sousremarque}\rm
 Dans la définition de l'hensélisé strict, on pourrait avoir envie de localiser un peu moins en imposant que l'extension résiduelle soit triviale dans les voisinages étales, c'est-à-dire en ne prenant la limite inductive que sur la famille des voisinages étales dont l'extension résiduelle est triviale. On obtient alors \emph{l'hensélisé} (non strict) et les anneaux locaux $\Oc_{X,x}^h$ ainsi obtenus sont dits \emph{henséliens}. Si $A$ est un anneau local, il est strictement hensélien si et seulement s'il est hensélien et si son corps résiduel est séparablement clos. Il y a un certain nombre de caractérisations utiles des anneaux henséliens ou strictement henséliens. Par exemple un anneau local $A$ est hensélien si et seulement s'il vérifie le lemme de Hensel, ou encore si et seulement si toute $A$-algèbre finie est un produit fini d'anneaux locaux.
\end{sousremarque}

Il y aurait encore beaucoup à dire sur les anneaux henséliens ou strictement henséliens, mais il ne nous a pas semblé opportun de nous étendre davantage sur ce sujet ici. Nous renvoyons donc le lecteur intéressé aux ouvrages classiques, par exemple~\cite{Milne_Etale_coh}, EGA~IV~\cite{EGA} ou encore~\cite{Raynaud_ALH}.

\section{Descente}
\label{descente}

Avant de traiter le c\oe ur de cette section, à savoir la descente fidèlement plate et ses applications, nous allons essayer d'expliquer la problématique sur un exemple simple. Essentiellement, la descente est un exercice de \og recollement\fg.

\paragraph*{Recollement Zariski de morphismes.}
Soient $U$ un schéma, $\{U_i\}_i$ un recouvrement Zariski de $U$ et $X$, $Y$ deux $U$-schémas. Pour tout $i$ on note $X_i$ et $Y_i$ les images réciproques de l'ouvert $U_i$ par les morphismes structuraux $X\to U$ et $Y\to U$ (c'est-à-dire les produits fibrés $X\times_U U_i$ et $Y\times_U U_i$). Donnons-nous pour tout $i$ un morphisme de schémas $\varphi_i : X_i \to Y_i$ et supposons que ces morphismes coïncident sur les intersections, \emph{i.e.} que pour tout couple $(i,j)$, les morphismes $\varphi_i$ et $\varphi_j$ soient égaux sur $X_i\cap X_j$. Alors il existe un unique morphisme de schémas $\varphi : X \to Y$ tel que pour tout $i$, $\varphi_{|_{X_i}}=\varphi_i$. La preuve est immédiate et laissée en exercice au lecteur.

On peut reformuler ce résultat de la manière suivante. Soient $S$ un schéma et $X, Y$ deux $S$-schémas. Alors le foncteur
$$\Homs_S(X,Y) : (Sch/S)^o \lto \Ens$$
qui à tout $S$-schéma $U$ associe l'ensemble des morphismes de $U$-schémas de $X\times_S U$ vers $Y\times_S U$, est un faisceau pour la topologie de Zariski. En particulier, pour $X=S$, on voit que le foncteur des points $h_Y$ de $Y$ est un faisceau pour la topologie de Zariski.

\paragraph*{Recollement Zariski de schémas.}
Soient $U$ un schéma et $\{U_i\}_i$ un recouvrement Zariski de $U$. Pour tout triplet d'indices $i, j, k$ on note $U_{ij}=U_i\cap U_j$ et $U_{ijk}=U_i\cap U_j \cap U_k$. Pour tout $i$, soit $f_i : X_i \to U_i$ un schéma au-dessus de $U_i$. On suppose que pour tout couple d'indices $i, j$, on a un isomorphisme $\varphi_{ij} : f_j^{-1}(U_{ij}) \to f_i^{-1}(U_{ij})$. On suppose de plus que ces isomorphismes vérifient la \og condition de cocycle\fg\ :
$$\varphi_{ik}=\varphi_{ij}\circ\varphi_{jk} : f_k^{-1}(U_{ijk}) \lto f_j^{-1}(U_{ijk}) \lto f_i^{-1}(U_{ijk})\, .$$
Alors il existe un unique $U$-schéma $f : X \to U$  et des isomorphismes $\varphi_i : f^{-1}(U_i) \to\! X_i$ tels que pour tous $i, j$ on ait :
$$\varphi_{ij} = \varphi_i \circ \varphi_j^{-1} : f_j^{-1}(U_{ij}) \lto f^{-1}(U_{ij}) \lto f_i^{-1}(U_{ij})\, .$$
Ici aussi la preuve est élémentaire et nous en dispenserons le lecteur. On peut de la même manière recoller bien d'autres objets : des faisceaux, des modules-quasi-cohérents, des courbes elliptiques, \emph{etc}.

\paragraph*{Et pour d'autres topologies ?} Une question naturelle se pose : les résultats ci-dessus sont-ils encore valables (\emph{mutatis mutandis}) si l'on remplace la famille couvrante Zariski~$\{U_i\}_i$ de~$U$ par une famille couvrante pour une topologie plus fine, par exemple la topologie \emph{fpqc} ? Nous verrons plus bas que la réponse est affirmative pour les morphismes de schémas~(\ref{descente_fpqc_gen}). En revanche le recollement de schémas est un peu plus délicat et ne marche pas toujours (voir~\ref{ddd_non_effectives}). Pour obtenir des résultats positifs, nous devrons faire des hypothèses supplémentaires sur les $X_i$, par exemple supposer qu'ils sont quasi-affines sur les $U_i$ (\ref{descente_quaff}), ou bien munis de faisceaux inversibles relativement amples (\ref{descente_ample}). La descente fidèlement plate des modules quasi-cohérents (\ref{fpqc_qcoh}) est la pièce maîtresse sur laquelle repose tout l'édifice. Nous donnons en~\ref{section_ddd} le vocabulaire nécessaire à l'énoncé des résultats de~\ref{descente_fpqc}.

\paragraph*{Descente de propriétés de morphismes.} On utilise encore le même terme de \og descente\fg\ pour un exercice légèrement différent. Considérons comme précédemment une famille couvrante~$\{U_i \to U\}$ pour une certaine topologie. Soit $X\to Y$ un $U$-morphisme. Supposons que pour tout $i$, le morphisme $X_i\to Y_i$ induit au-dessus de $U_i$ vérifie une certaine propriété $\Pc$. Alors le morphisme $X\to Y$ vérifie-t-il $\Pc$ ? Autrement dit, la propriété~$\Pc$ \og descend-elle\fg\ à $X\to Y$ ? Lorsque c'est le cas, on dit que $\Pc$ est \og de nature locale \fg\ au but pour la topologie considérée. Nous verrons en~\ref{section_permanence} que de nombreuses propriétés intéressantes de morphismes sont de nature locale au but pour la topologie~\emph{fpqc}.   

\subsection{Données de descente}
\label{section_ddd}

\begin{sousdefi}
 Une \emph{catégorie fibrée} $\Fc$ sur une catégorie $\Cc$ consiste en les données suivantes.
\begin{itemize}
 \item[(i)] Pour tout objet $U$ de $\Cc$, une catégorie $\Fc(U)$.
\item[(ii)] Pour tout morphisme $f : U\to V$ de $\Cc$, un foncteur $f^* : \Fc(V) \to \Fc(U)$.
\item[(iii)] Pour tout objet $U$ de $\Cc$, un isomorphisme de foncteurs $\eps_U : (\id_U)^*\simeq \id_{\Fc(U)}$.
\item[(iv)] Pour tout couple de flèches $\xymatrix@C=1pc{U\ar[r]^f &V \ar[r]^g&W}$ un isomorphisme de foncteurs $\alpha_{f,g} : f^*g^*\simeq (gf)^* : \Fc(W)\to \Fc(U)$.
\end{itemize}
Ces données sont assujetties à des conditions de compatibilité évidentes dont nous laissons les détails au lecteur (si $f : U\to V$ est une flèche dans $\Cc$, compatibilité de $\alpha_{\id_U,f}$ avec $\eps_U$ et de $\alpha_{f,\id_V}$ avec $\eps_V$ ; si $f, g, h$ sont trois flèches composables, diagramme d'\og associativité\fg\ pour les foncteurs $f^*, g^*, h^*$).
\end{sousdefi}

\begin{sousremarque}\rm
 Le lecteur attentif aura sans doute remarqué que la définition donnée ici n'est pas équivalente à la définition usuelle d'une catégorie fibrée. En fait, ce que nous avons défini ci-dessus est plutôt appelé \og pseudo-foncteur\fg\ dans la littérature. La donnée d'un pseudo-foncteur est équivalente à la donnée d'une catégorie fibrée munie d'un \og clivage\fg\ (avec la terminologie de Vistoli dans~\cite{FGA_explained}). Avec l'axiome du choix, toute catégorie fibrée admet un clivage. L'abus de langage ci-dessus est donc sans dommages, et -- je l'espère -- compensé par le fait que la notion de pseudo-foncteur semble plus intuitive que celle de catégorie fibrée (ce point de vue subjectif n'engage bien sûr que moi). Pour plus de détails sur ces questions, on pourra consulter l'exposé de Vistoli dans~\cite{FGA_explained}. 
\end{sousremarque}

\begin{sousexemple}\rm
\label{def_qcoh}
 Soient $S$ un schéma et $\Cc$ la catégorie des $S$-schémas. La catégorie fibrée $\Qcoh$ des modules quasi-cohérents sur $\Cc$ est la catégorie fibrée qui à tout $S$-schéma $U$ associe la catégorie $\Qcoh(U)$ des modules quasi-cohérents sur $U$. Si $f : U\to V$ est un morphisme de $S$-schémas, le foncteur de changement de base $f^*$ est le foncteur qui à~$\Mc$ associe $\Mc\otimes_{\Oc_V}\Oc_U$. Les isomorphismes $\eps_U$ et $\alpha_{f,g}$ sont les isomorphismes canoniques $\Mc\otimes_{\Oc_U}\Oc_U\simeq \Mc$ et $(\Mc\otimes_{\Oc_W}\Oc_V)\otimes_{\Oc_V}\Oc_U \simeq \Mc\otimes_{\Oc_W}\Oc_U$.
\end{sousexemple}

\begin{sousexemple}\rm
\label{cat_fibree_des_schemas}
 Toujours sur la catégorie $\Cc$ des $S$-schémas, la catégorie fibrée des schémas sur $\Cc$ est la catégorie fibrée qui à tout $S$-schéma $U$ associe la catégorie des $U$-schémas. Si~$f : U\to V$ est un morphisme de $S$-schémas, le foncteur de changement de base $f^*$ est le foncteur qui à $X\to V$ associe $X\times_V U \to U$. 
\end{sousexemple}

\begin{sousdefi}[donnée de descente]
 Soient $\Cc$ un site (autrement dit, une catégorie munie d'une topologie de Grothendieck) et $\Fc$ une catégorie fibrée sur $\Cc$. Soit $U$ un objet de $\Cc$ et soit $\Uc=\{U_i\to U\}_i$ un recouvrement de $U$. On note $U_{ij}$ (resp. $U_{ijk}$) le produit fibré $U_i\times_U~U_j$ (resp. $U_i\times_UU_j\times_U U_k$) et $\pr_1, \pr_2, \pr_{12}, \pr_{13}, \pr_{23}, q_1, q_2, q_3$ les projections canoniques (avec $\pr_1 : U_{ij} \to U_i$, $\pr_{12} : U_{ijk} \to U_{ij}$, $q_1 : U_{ijk} \to U_i$, \emph{etc.}). Pour tout~$i$, soit $\xi_i$ un objet de~$\Fc(U_i)$. Une \emph{donnée de descente} sur la famille $(\xi_i)$ est un ensemble d'isomorphismes $$\varphi_{ij} : \pr_2^*\xi_j \lto \pr_1^*\xi_i$$ dans $\Fc(U_{ij})$ vérifiant la condition de cocycle suivante : pour tout triplet d'indices $(i,j,k)$, on a dans $\Fc(U_{ijk})$ l'égalité\footnote{aux isomorphismes structuraux de $\Fc$ près, $q_1^*\simeq \pr_{12}^*p_1^*\simeq \pr_{13}^*p_1^*$, \emph{etc.} On remarque que cette condition de cocycle est l'analogue naturel de celle qui apparaissait dans le \og recollement Zariski de schémas\fg.}
$$\pr_{13}^*(\varphi_{ik}) = \pr_{12}^*(\varphi_{ij}) \circ \pr_{23}^*(\varphi_{jk}) : 
q_3^*(\xi_k) \lto q_1^*(\xi_i)\, .$$
On dit que $((\xi_i), (\varphi_{ij}))$ est un \emph{objet muni d'une donnée de descente} relativement au recouvrement $\Uc$.

Si $((\xi_i), (\varphi_{ij}))$ et $((\eta_i), (\psi_{ij}))$ sont deux objets munis de données de descente, un morphisme du premier dans le second est un ensemble de flèches $\alpha_i : \xi_i \to \eta_i$ dans $\Fc(U_i)$ tel que pour tout couple d'indices $i,j$ le diagramme suivant commute.
$$\xymatrix{\pr_2^*\xi_j \ar[r]^{\pr_2^*\alpha_j} \ar[d]^{\varphi_{ij}} &
\pr_2^*\eta_j \ar[d]^{\psi_{ij}}\\
\pr_1^*\xi_i \ar[r]^{\pr_1^*\alpha_i} & \pr_1^* \eta_i
}$$

On note $\Fc(\Uc/U)$ la catégorie des objets munis d'une donnée de descente relativement au recouvrement $\Uc$ de $U$ ainsi constituée.
\end{sousdefi}

\begin{soustexte}\rm
\label{foncteur_descente}
Si $\xi$ est un objet de $\Fc(U)$, on peut lui associer naturellement un objet muni d'une donnée de descente de la manière suivante. Pour tout $i$, on pose $\xi_i=u_i^*\xi$ où $u_i$ est le morphisme donné $U_i \to U$, et pour tout couple $i,j$ d'indices, on prend pour $\varphi_{ij}$ l'isomorphisme canonique de $\pr_2^*u_j^*\xi$ dans $\pr_1^*u_i^*\xi$.

Si $\xi \to \eta$ est un morphisme dans $\Fc(U)$, on a un morphisme naturel de l'objet muni d'une donnée de descente associé à $\xi$, vers celui associé à $\eta$. Ceci définit donc un foncteur $$\Fc(U)\lto \Fc(\Uc/U)\, .$$

Les exemples vus au début du paragraphe~\ref{descente} montrent que pour la catégorie fibrée des schémas sur $(Sch/S)$ de~\ref{cat_fibree_des_schemas}, le foncteur $\Fc(U) \to \Fc(\Uc/U)$ est une équivalence de catégories pour toute famille couvrante $\Uc$ de $U$ au sens de Zariski. Plus précisément, le \og recollement de morphismes\fg\ montre que ce foncteur est pleinement fidèle, et le \og recollement de schémas\fg\ montre qu'il est essentiellement surjectif.
\end{soustexte}

\begin{sousdefi}
 Une donnée de descente sur une famille d'objets $(\xi_i)$ est dite \emph{effective} s'il existe un objet $\xi$ de $\Fc(U)$ qui induit (à isomorphisme près) la famille $(\xi_i)$ avec sa donnée de descente. 
\end{sousdefi}

Il revient donc au même de dire qu'une donnée de descente est effective, ou que le couple constitué de la famille d'objets et de cette donnée de descente est dans l'image essentielle du foncteur $\Fc(U)\to \Fc(\Uc/U)$. Dans la pratique, on se ramène souvent au cas où la famille couvrante $\Uc$ est constituée d'une seule flèche $V\to U$. On note alors parfois $\Fc(V/U)$ au lieu de $\Fc(\Uc/U)$. 

Donnons encore un peu de vocabulaire dans ce contexte. Soient $\Fc$ une catégorie fibrée sur $(Sch/S)$ et $\pi : V\to U$ un morphisme de $S$-schémas. On dit que $\pi$ est un \emph{morphisme de descente} pour $\Fc$ si le foncteur $\Fc(U) \to \Fc(V/U)$ est pleinement fidèle. On dit que $\pi$ est un \emph{morphisme de descente effective} pour $\Fc$ si de plus toute donnée de descente est effective, \emph{i.e.} si $\Fc(U) \to \Fc(V/U)$ est une équivalence.

\begin{sousdefi}
 \label{def_champ}
Soit $\Fc$ une catégorie fibrée sur un site $\Cc$. On dit que $\Fc$ est un \emph{champ} si pour tout objet $U$ de $\Cc$ et toute famille couvrante $\Uc$ de $U$, le foncteur $\Fc(U)\to \Fc(\Uc/U)$ est une équivalence de catégories.
\end{sousdefi}

\subsection{Descente \emph{fpqc} des modules quasi-cohérents et applications}
\label{descente_fpqc}

\begin{sousthm}
\label{fpqc_qcoh}
 Soient $U$ un schéma et $\Uc=\{U_i\to U\}_i$ une famille couvrante \emph{fpqc}. Alors le foncteur
$$\Qcoh(U) \lto \Qcoh(\Uc/U)$$
défini en \ref{foncteur_descente} est une équivalence de catégories. (Autrement dit, la catégorie fibrée des modules quasi-cohérents définie en~\ref{def_qcoh} est un champ pour la topologie \emph{fpqc}.)
\end{sousthm}

La pleine fidélité signifie que, pour tout schéma $S$ et tout couple de modules quasi-cohérents $\Fc, \Gc$ sur $S$, le foncteur $\Homs_{\Oc_S}(\Fc, \Gc)$ est un faisceau \emph{fpqc}. L'essentielle surjectivité signifie que toute donnée de descente sur une collection de $\Oc_{U_i}$-modules quasi-cohérents est effective.

\vskip 2mm

\noindent \begin{demo}
 Nous donnons seulement une esquisse. Des arguments certes un peu longs et fastidieux, mais élémentaires et de nature essentiellement formelle\footnote{D'ailleurs, les mêmes arguments sont valables pour n'importe quelle catégorie fibrée sur $(Sch/S)$. Autrement dit une telle catégorie fibrée $\Fc$ est un champ \emph{fpqc} si et seulement si c'est un champ Zariski et si pour tout $V\to U$ fidèlement plat avec $U$ et $V$ affines, le foncteur $\Fc(U)\to \Fc(V/U)$ est une équivalence.}, montrent qu'il suffit de prouver le théorème dans les deux cas particuliers suivants.
\begin{itemize}
 \item[a)] La famille couvrante $\{U_i\to U\}_i$ est une famille couvrante pour la topologie de Zariski.
\item[b)] La famille couvrante $\{U_i\to U\}_i$ est réduite à un morphisme $V \to U$ avec $U$ et~$V$ affines.
\end{itemize}
Le cas a) est trivial. Reste le b). Notons $U=\Spec A$ et $V=\Spec B$. \emph{Via} l'équivalence bien connue entre la catégorie des modules sur un anneau et la catégorie des faisceaux de modules quasi-cohérents sur son spectre, on se ramène à une question portant uniquement sur des modules. Pour la pleine fidélité, il faut montrer que si $N$ et $M$ sont deux $A$-modules, alors le diagramme naturel (dont nous laissons au lecteur le soin d'expliciter les flèches)
$$\xymatrix{\Hom_A(M,N) \ar[r] & \Hom_B(M\otimes_A B, N\otimes_A B) \ar@<2pt>[r]^-{p_1^*} \ar@<-2pt>[r]_-{p_2^*} &
\Hom_{B^{\otimes 2}}(M\otimes_A B^{\otimes 2}, N\otimes_A B^{\otimes 2})
}$$
 est exact. C'est une conséquence facile du lemme~\ref{lemme_descente_fpqc_qcoh} ci-dessous.

Pour l'essentielle surjectivité, on part d'un $B$-module $M'$ muni d'une donnée de descente, c'est-à-dire d'un isomorphisme de $B^{\otimes 2}$-modules $\varphi : M'\otimes_A B \simeq B\otimes_A M'$ assujetti à une condition de cocycle. Il faut montrer que le couple $(M',\varphi)$ provient d'un $A$-module $M$ (à isomorphisme près). En regardant d'une part le morphisme canonique $M' \to M'\otimes_A B$, et d'autre part le composé du morphisme canonique $M'\to B\otimes_A M'$ et de $\varphi^{-1}$, on a un couple de flèches
$\xymatrix@C=1pc{M' \ar@<2pt>[r] \ar@<-2pt>[r] & M'\otimes_A B}$. On note $M$ le noyau de ce couple de flèches. C'est naturellement un $A$-module. Il reste à vérifier qu'il convient, ce qui demande encore du travail et est laissé en exercice (on pourra, comme dans la preuve du lemme ci-dessous, se ramener \emph{via} le changement de base fidèlement plat $A\to B$ au cas où $A\to B$ a une rétraction). Pour plus de détails, consulter par exemple~\cite{BRL}~6.1.
\end{demo}

\begin{souslem}
 \label{lemme_descente_fpqc_qcoh}
Soit $A\to B$ un morphisme fidèlement plat d'anneaux. Alors, pour tout $A$-module $M$, le diagramme naturel
$$\xymatrix{M \ar[r] & M\otimes_A B \ar@<2pt>[r] \ar@<-2pt>[r] &
M\otimes_A B^{\otimes 2}
}$$
est exact.
\end{souslem}
\begin{demo}
 Si $A\to B$ a une rétraction, toutes les flèches du diagramme en ont et l'exactitude est immédiate. On ramène le cas général à ce cas particulier par le changement de base fidèlement plat $A\to B$.
\end{demo}

\begin{sousremarque}\rm
 Le théorème~\ref{fpqc_qcoh} serait faux avec une topologie fidèlement plate trop fine, sans hypothèse de finitude sur les recouvrements. Un contre-exemple est donné dans le cours de Vistoli~\cite{FGA_explained}~4.24.
\end{sousremarque}

\begin{sousexo}\rm
 Déduire de~\ref{fpqc_qcoh} le fait que les catégories fibrées suivantes sont des champs \emph{fpqc} :
\begin{itemize}
 \item la catégorie fibrée des modules cohérents ;
\item la catégorie fibrée des modules localement libres de rang $n$.
\end{itemize}
\end{sousexo}

\begin{sousthm}
\label{descente_fpqc_gen}
 Soit $\pi : S'\to S$ un morphisme \emph{fpqc} de schémas et soient $X$ et $Y$ deux $S$-schémas. On note $X'$, $Y'$, $X''$, $Y''$ les produits fibrés auxquels on pense. Alors le diagramme
$$\xymatrix{\Hom_S(X,Y) \ar[r]^{\pi^*} &\Hom_{S'}(X',Y') \ar@<2pt>[r]^{p_1^*} \ar@<-2pt>[r]_{p_2^*}&
\Hom_{S''}(X'',Y'') }$$
est exact. 
\end{sousthm}
\begin{demo}
 La question est locale sur $S$ et $Y$. On peut donc supposer qu'ils sont affines. De plus, quitte à remplacer $S'$ par une somme disjointe finie d'affines, on peut supposer que $S'$ est affine. On peut alors reformuler le problème en termes d'algèbres quasi-cohérentes et appliquer le théorème~\ref{fpqc_qcoh}.
\end{demo}

Ce théorème montre que $\Homs_S(X,Y)$ est un faisceau \emph{fpqc}. Autrement dit, les morphismes \emph{fpqc} sont des morphismes de descente pour la catégorie fibrée des schémas sur $(Sch/S)$, ou encore : pour $\pi$ comme dans l'énoncé et $\Fc$ la catégorie fibrée des schémas, le foncteur $\Fc(S) \to \Fc(S'/S)$ de~\ref{foncteur_descente} est pleinement fidèle. En prenant $X=S$, on obtient le cas particulier important mentionné en~\ref{schema_implique_faisceau} : $h_Y$ est un faisceau \emph{fpqc}.
En général, la descente \emph{fpqc} n'est pas effective. Elle l'est cependant dans plusieurs cas utiles, le premier d'entre eux étant le cas des morphismes affines (ou même quasi-affines).

\begin{sousthm}
\label{descente_quaff}
 Soit $\pi : S'\to S$ un morphisme fidèlement plat et quasi-compact de schémas et soit $(X', \varphi)$ un $S'$-schéma muni d'une donnée de descente, \emph{i.e.} $\varphi$ est un isomorphisme $p_1^*X' \to p_2^*X'$ qui vérifie la condition de cocycle usuelle (avec les notations usuelles). On suppose que le morphisme $X'\to S'$ est affine (resp. une immersion ouverte, resp. quasi-affine). Alors la donnée de descente est effective. De plus le morphisme $X\to S$ qui induit $X'\to S'$ est affine (resp. une immersion ouverte, resp. quasi-affine).
\end{sousthm}
\begin{demo}
 Ici aussi nous donnons seulement une esquisse, pour les détails nous renvoyons par exemple à~\cite{BRL}.

Puisque la catégorie des schémas affines sur une base $S$ est anti-équivalente à la catégorie des $\Oc_S$-algèbres quasi-cohérentes, le cas affine est une conséquence du théorème~\ref{fpqc_qcoh}.

Supposons maintenant que $X'\to S'$ est une immersion ouverte, ou ce qui revient au même que $X'$ est un ouvert de $S'$. On note $p_1, p_2$ les projections canoniques de $S''=S'\times_S S'$ sur $S'$. La présence d'une donnée de descente sur $X'$ montre alors que $p_1^{-1}(X')=p_2^{-1}(X')$. On en déduit que $X'=\pi^{-1}(\pi(X'))$. Comme le morphisme $\pi$ est submersif (la topologie de $S$ est quotient de celle de $S'$), ceci entraîne que $\pi(X')$ est ouvert dans $S$. On pose alors $X=\pi(X')$ et on vérifie que ce choix convient.

Comme un morphisme quasi-affine est le composé d'une immersion ouverte quasi-compacte et d'un morphisme affine, le cas général va se déduire des deux précédents. Plus précisément, si $f' : X'\to S'$ est quasi-affine, alors dans la factorisation de Stein
$$X'\lto Z'=\Spec f'_*(\Oc_{X'}) \lto S'\, ,$$
la première flèche est une immersion ouverte quasi-compacte. Comme la formation de $f'_*(\Oc_{X'})$ commute au changement de base plat, on voit que la donnée de descente sur $X'$ induit une donnée de descente sur  $f'_*(\Oc_{X'})$, donc sur $Z'$. Cette donnée de descente est effective d'après le cas affine et on trouve $Z$ affine sur $S$ qui induit $Z'$. Le cas d'une immersion ouverte permet alors de trouver $U$ ouvert de $Z$ qui convient. Il reste juste à montrer que l'immersion ouverte $U\to Z$ est quasi-compacte. C'est une conséquence de~\ref{permanence}.
\end{demo}

Ainsi, une donnée de descente sur un schéma quasi-affine est effective. La proposition suivante montre que l'on peut même se contenter de recouvrir $X'$ par des schémas quasi-affines stables par la donnée de descente. Plus précisément, soit $\pi : S'\to S$ un morphisme fidèlement plat et quasi-compact de schémas et soit $(X', \varphi)$ un $S'$-schéma muni d'une donnée de descente. Un ouvert $U'$ de $X'$ est dit \emph{stable} par $\varphi$ si $\varphi$ induit une donnée de descente sur $U'$, c'est-à-dire si $\varphi$ induit un isomorphisme $\pr_2^*U' \simeq \pr_1^*U'$.

\begin{sousprop}
\label{descente_locale}
 On suppose que $X'$ est recouvert par des ouverts $X'_i$ stables par la donnée de descente. Pour que la donnée de descente sur $X'$ soit effective, il faut et il suffit qu'il en soit de même des données de descente induites sur les $X'_i$.
\end{sousprop}
\begin{demo}
 SGA 1~\cite{SGA1}~VIII~7.2
\end{demo}

\begin{sousexo}[SGA 1~\cite{SGA1}~VIII~7.5]\rm
 Soit $S'\to S$ un morphisme fidèlement plat, quasi-compact et radiciel (\emph{i.e.} universellement injectif). Soit $X'$ un $S'$-schéma muni d'une donnée de descente $\varphi$.

1) Montrer que $S'$ est séparé sur $S$. 

2) Montrer que tout ouvert $U$ de $X'$ est stable par $\varphi$.

3) En déduire que la donnée de descente est effective. En d'autres termes : un morphisme fidèlement plat, quasi-compact et radiciel est de descente effective pour la catégorie fibrée des schémas.
\end{sousexo}

\begin{sousdefi}\ 

\begin{itemize}
 \item[a)] Soit $X$ un schéma quasi-compact et quasi-séparé. Un faisceau inversible $\Lc$ sur $X$ est dit ample s'il existe un entier $n>0$ tel que $\Lc^{\otimes n}$ soit engendré par des sections globales $s_1, \dots, s_r$ telles que pour tout $i$, le lieu $X_{s_i}$ où $s_i$ engendre $\Lc^{\otimes n}$ soit quasi-affine. (Alors, pour tout $n>0$ et toute section globale $s$ de $\Lc^{\otimes n}$, l'ouvert $X_s$ de $X$ est quasi-affine.)
\item[b)] Soient $S$ un schéma de base et $X$ un $S$-schéma quasi-compact et quasi-séparé. Un faisceau inversible $\Lc$ sur $X$ est dit $S$-ample s'il existe un recouvrement ouvert affine~$\{S_j\}$ de $S$ tel que pour tout $j$, la restriction de $\Lc$ à $X\times_S S_j$ soit ample. (On montre qu'alors, pour tout ouvert affine $U$ de $S$, la restriction de $\Lc$ à $X\times_S U$ est ample.)
\end{itemize}
\end{sousdefi}

\begin{sousthm}
\label{descente_ample}
  Soit $\pi : S'\to S$ un morphisme fidèlement plat et quasi-compact de schémas et soit $(X', \varphi)$ un $S'$-schéma muni d'une donnée de descente. Soit $\Lc'$ un faisceau inversible $S'$-ample sur $X'$ et soit $\lambda$ une donnée de descente sur $\Lc'$ compatible avec $\varphi$. Alors la donnée de descente sur $X'$ est effective, et le couple $(X',\Lc')$ descend en un couple~$(X,\Lc)$ où $\Lc$ est un faisceau inversible $S$-ample sur $X$.
\end{sousthm}

Avant d'esquisser la démonstration, précisons un peu l'énoncé. La donnée de descente sur $X'$ est un $(S'\times_S S')$-isomorphisme $\varphi : S'\times_S X' \to X'\times_S S'$ qui vérifie une certaine condition de cocycle. En notant $q_1 : X'\times_S S'\to X'$ la projection sur le premier facteur et~$q_2~:~X'\times_S~S'~\to~X'$ la composée de $\varphi^{-1}$ et de la projection $S'\times_S X'\to X'$ sur le second facteur, on obtient un couple de flèches $\xymatrix@C=1pc{X'\times_S S' \ar@<2pt>[r] \ar@<-2pt>[r]& X'}$. La donnée de descente $\lambda$ \og compatible avec $\varphi$\fg\ est alors un isomorphisme $\lambda : q_2^*\Lc' \to q_1^*\Lc'$ qui vérifie une certaine condition de cocycle.

\begin{sousexo}\rm
 Expliciter cette condition de cocycle.
\end{sousexo}
\begin{demo}
 Encore une fois, nous donnons seulement une esquisse (tirée d'ailleurs de~\cite{BRL}). On peut supposer $S$ et $S'$ affines. On note $f' : X'\to S'$ le morphisme structural de $X'$. On note $\Mc'$ la $S'$-algèbre graduée $\bigoplus_{n\geq 0} f_*'(\Lc'^{\otimes n})$. C'est une $S'$-algèbre quasi-cohérente. D'après~\ref{fpqc_qcoh} elle descend donc à une algèbre quasi-cohérente $\Mc$ sur $S$. De plus la graduation naturelle sur $\Mc'$ induit une graduation $\Mc=\bigoplus_{n\geq 0} \Mc_n$ sur $\Mc$. Soit~$s'$ une section globale  de $\Lc'^{\otimes n}$ pour un certain $n$. On peut écrire
$$s' = \sum_{i=1}^m a_i\otimes s_i$$
où les $a_i$ sont des sections globales de $\Oc_{S'}$ et les $s_i$ sont des sections globales de $\Mc_n$. Si, en un point $x'\in X'$, la section $s'$ engendre $\Lc'^{\otimes n}$, alors au moins une des sections $1\otimes s_i$ doit engendrer $\Lc'^{\otimes n}$ en $x'$. On peut donc recouvrir $X'$ par des ouverts quasi-affines $X_{s}'$ où $s$ est une section globale d'un $\Lc'^{\otimes n}$ qui descend en une section globale de $\Mc$. Ceci montre que les $X'_s$ sont stables par $\varphi$, donc cette donnée de descente sur $X'$ est effective d'après~\ref{descente_locale}. Enfin, en appliquant de nouveau~\ref{fpqc_qcoh}, on voit que le faisceau inversible $\Lc'$ descend à $X$.
\end{demo}

Lorsque l'on souhaite utiliser ce théorème en pratique, il n'est pas toujours aisé de vérifier que le faisceau inversible~$\Lc'$ est ample relativement à~$S'$. À toutes fins utiles, rappelons ici le critère d'amplitude relative EGA IV~\cite{EGA}~9.6.5. Soit~$X$ un schéma propre et de présentation finie sur une base~$S$, et soit~$\Lc$ un faisceau inversible sur~$X$. Alors $\Lc$ est $S$-ample si et seulement si sa restriction à chaque fibre de $X\to S$ est ample.

\begin{sousexo}\rm
 Soit $\Mc_g$ la catégorie fibrée sur $(Sch/S)$ définie ainsi. Pour tout $U$, les objets de la catégorie $\Mc_g(U)$ sont les morphismes propres et lisses $\pi : C\to U$ dont les fibres géométriques sont des courbes connexes de genre $g$, et les flèches de $\Mc_g(U)$ sont les $U$-isomorphismes. Montrer que, si $g\neq 1$, alors $\Mc_g$ est un champ pour la topologie \emph{fpqc}. [Indication : On utilisera le théorème précédent. Pour $g\geq 2$, on pourra remarquer que, si $\pi : C\to U$ est un objet de $\Mc_g(U)$, alors le fibré canonique $\Omega_{C/U}$ sur $X$ est ample relativement à $\pi$. Pour $g=0$ on peut prendre son dual.]
\end{sousexo}

\begin{sousremarque}\rm
 Pour $g=1$ il n'y a pas de fibré relativement ample \emph{canonique} sur une famille de courbes $\pi : C\to U$. On ne pourra donc pas appliquer le théorème précédent. Et de fait, un exemple de Raynaud (voir~\cite{Raynaud_Faisceaux_amples_sur}~XIII~3.2) montre que $\Mc_1$ \emph{n'est pas} un champ pour la topologie \emph{fpqc}.
\end{sousremarque}

\begin{sousexemple}[données de descente non effectives]\rm
\label{ddd_non_effectives}
 Nous donnons en~\ref{exemple_espace_alg} un exemple d'espace algébrique qui n'est pas un schéma. On en déduit l'existence d'une donnée de descente non effective sur la droite affine complexe avec origine dédoublée, relative au recouvrement $\Spec \C \to \Spec \R$. Un autre exemple est donné en~4.4.2 dans l'exposé de Vistoli de~\cite{FGA_explained}.
\end{sousexemple}

Nous terminons cette partie par une application des théorèmes de descente ci-dessus à la question de la représentabilité d'un faisceau.

\begin{souslem}[descente de la représentabilité]\ 
\label{descente_representabilite}

\noindent Soient $S$ un schéma et $F : (Sch/S)^o \to \Ens$ un foncteur.
\begin{itemize}
 \item[(i)] On suppose que $F$ est un faisceau pour la topologie de Zariski. Soit $\{S_i\}_i$ un recouvrement ouvert (Zariski) de $S$ tel que la restriction $F_i=F\times_S S_i$ de $F$ à chaque~$S_i$ soit représentable par un $S_i$-schéma $X_i$. Alors $F$ est représentable par un~$S$-schéma~$X$. (Autrement dit : pour un faisceau Zariski, être représentable par un schéma est une condition de nature locale sur $S$ pour la topologie de Zariski.)
\item[(ii)] On suppose que $F$ est un faisceau \emph{fpqc} (resp. \emph{fppf}, resp. étale). Soit $\{S_i\to S\}_i$ une famille couvrante \emph{fpqc} (resp. \emph{fppf}, resp. étale) telle que la restriction $F_i=F\times_S S_i$ de $F$ à chaque $S_i$ soit représentable par un $S_i$-schéma $X_i$. Alors la famille des $X_i$ est munie d'une donnée de descente relativement à la famille couvrante $\{S_i\to S\}_i$.

Si de plus cette donnée de descente est effective (ce qui est le cas par exemple si chaque $X_i$ est quasi-affine sur $S_i$), alors $F$ est représentable par un $S$-schéma $X$ (quasi-affine sur $S$ dans le cas particulier de la parenthèse précédente).
\end{itemize}
\end{souslem}
\begin{demo}
 Voir la réédition de SGA~3, \cite{SGA3}~VIII~1.7.2
\end{demo}

\subsection{Propriétés de permanence}
\label{section_permanence}

Si $\pi : S'\to S$ est un morphisme de changement de base, et si $f$ est un $S$-morphisme de schémas, on peut se demander dans quelle mesure une propriété du morphisme $f'$ obtenu par changement de base \og descend\fg\ à $f$. Nous avons rassemblé ci-dessous un certain nombre de cas traités dans les EGA. En fait la plupart des propriétés intéressantes descendent par morphismes \emph{fpqc}. Lorsqu'une hypothèse plus faible sur $\pi$ est suffisante, nous l'avons signalé entre parenthèses. 

\begin{sousprop}
 \label{permanence}
Soit $\pi : S'\to S$ un morphisme fidèlement plat et quasi-compact de schémas et soit $f : X \to Y$ un morphisme de $S$-schémas. On note $f' : X' \to Y'$ le morphisme obtenu par le changement de base $\pi$. Considérons, pour un morphisme, la propriété d'être :
\begin{itemize}
 \item injectif ($\pi$ surjectif suffit) ;
\item surjectif ($\pi$ surjectif suffit) ;
\item à fibres ensemblistement finies ($\pi$ surjectif suffit) ;
\item bijectif ($\pi$ surjectif suffit) ;
\item radiciel ($\pi$ surjectif suffit) ;
\item ouvert ;
\item universellement ouvert ;
\item fermé ;
\item universellement fermé ;
\item un homéomorphisme ;
\item universellement un homéomorphisme ;
\item quasi-compact ($\pi$ quasi-compact et surjectif suffit) ;
\item quasi-compact et dominant ;
\item séparé ;
\item quasi-séparé ;
\item de type fini ;
\item localement de type fini ;
\item de présentation finie ;
\item localement de présentation finie ;
\item propre ;
\item un isomorphisme ;
\item un monomorphisme ;
\item une immersion ouverte ;
\item une immersion fermée ;
\item affine ;
\item quasi-affine ;
\item fini ;
\item quasi-fini ;
\item entier ;
\item plat ;
\item fidèlement plat ;
\item lisse ;
\item net (\emph{i.e.} non-ramifié) ;
\item étale.
\end{itemize}
Alors, si $\Pc$ désigne l'une des propriétés précédentes, pour que $f$ ait la propriété $\Pc$ il suffit que $f'$ la possède.
\end{sousprop}
\begin{demo}
 EGA IV \cite{EGA}, 2.2.13, 2.6.1, 2.6.2, 2.6.4, 2.7.1, 17.7.3.
\end{demo}

\begin{sousremarque}\rm
\label{exemple_Hironaka}
 Les grands absents de la liste ci-dessus sont les morphismes projectifs ou quasi-projectifs. De fait, Hironaka a donné dans~\cite{Hironaka_nonprojective_example} un exemple de morphisme propre non projectif $f : \widetilde{X} \to X$, où $X$ est réunion de deux ouverts $U_1$ et $U_2$ tels que les deux morphismes $f^{-1}(U_i)\to U_i$ soient projectifs. Dans cet exemple, $\widetilde{X}$ est une variété de dimension~3 sur~$\C$, obtenue à partir d'une variété projective lisse de dimension~3 en faisant des éclatements astucieux le long de certaines courbes, puis un recollement. On pourra consulter~\cite{Hartshorne_AG}~App.~B, 3.4.1 pour plus de détails et un joli dessin.
\end{sousremarque}

\section{Quotients}

\subsection{Schémas en groupes}
\label{schema_en_groupes}

\begin{sousdefi}
Soit $S$ un schéma. Un schéma en groupes sur $S$ est un objet en groupes dans la catégorie des $S$-schémas. 
\end{sousdefi}

De manière équivalente, c'est un $S$-schéma $G$ muni d'un morphisme de $S$-schémas
$$\mu : G\times_S G \lto G$$
qui vérifie un certain nombre d'axiomes : existence d'un morphisme inverse $i : G\to G$ et commutativité de quelques diagrammes. D'après le lemme de Yoneda, la donnée d'une structure de schéma en groupes sur $G$ est encore équivalente à la donnée d'une factorisation
$$h_G : (Sch/S)^o \lto \Gr \lto \Ens$$
de son foncteur des points à travers la catégorie des groupes.

\begin{sousexemple}\rm
 Nous avons déjà parlé plus haut des schémas en groupes $\Ga, \Gm$ et $\mu_n$. De nombreux autres exemples sont donnés par les groupes linéaires. Ainsi $GL_n$ est défini fonctoriellement de la manière suivante. Pour tout schéma $U$, $GL_n(U)$ est l'ensemble des matrices carrées de taille $n\times n$ à coefficients dans $\Gamma(U,\Oc_U)$ et dont le déterminant est un inversible de $\Gamma(U,\Oc_U)$. Il est représentable par le spectre de $\Z[X_{11}, X_{12}, \dots, X_{nn}]/(D)$ où $D$ est le polynôme en les variables $X_{11}, X_{12}, \dots, X_{nn}$ qui représente le déterminant. Les sous-groupes $O_n$, $\SL_n$, \emph{etc.} sont définis de manière analogue.
\end{sousexemple}

\begin{sousexemple}[groupes constants]\rm
  Soit $M$ un groupe commutatif ordinaire et $S$ un schéma. Nous noterons ici $M_S$ le faisceau constant associé à l'ensemble $M$ (pour éviter toute ambiguïté sur la base). Nous avons vu en~\ref{faisceau_constant} que $M_S$ est naturellement muni d'une structure de schéma en groupes sur $S$.
\end{sousexemple}

\begin{sousexemple}[groupes diagonalisables] \rm Soit $M$ un groupe commutatif ordinaire. On note alors $D_S(M)$ le foncteur sur $(Sch/S)^o$ défini par
$$D_S(M)(U)=\Hom_{U-Gr}(M_U, {\Gm}_{,U})=\Hom_{Gr}(M,\Gamma(U,\Oc_U)^{\times}).$$
On vérifie alors facilement que $D_S(M)$ est représentable par le spectre de l'algèbre $\Oc_S[M]$ du groupe $M$ sur $\Oc_S$. Un \emph{groupe diagonalisable} est un groupe de la forme $D_S(M)$ pour un certain groupe commutatif $M$. Par exemple, si $M=\Z/n\Z$, on obtient $\mu_n$, et si $M=\Z$, on a $D_S(M)={\Gm}_{,S}$. Plus généralement, si $M$ est un groupe abélien de type fini, disons le produit de $\Z^r$ par un produit fini de facteurs $\Z/n_i\Z$, $n_i>0$, alors $D_S(M)$ s'identifie à ${\Gm}^r_{,S}\times_S \prod_{i} {\mu_{n_i}}_{,S}$.

Revenons au cas général. L'égalité $D_S(M)=\Spec \Oc_S[M]$ montre que $D_S(M)$ est affine sur $S$. Son algèbre étant libre sur $\Oc_S$, il est de plus fidèlement plat. Nous retiendrons en particulier que le morphisme structural $D_S(M)\to S$ est fidèlement plat et quasi-compact, ce qui sera utile dans les questions de descente. Nous donnons ci-dessous quelques propriétés élémentaires. Nous renvoyons à SGA~3~\cite{SGA3}, ou au cours de J. Oesterlé sur les groupes de type multiplicatif dans la même école d'été, pour de plus amples développements. Nous donnerons en~\ref{cas_groupe_diago} un théorème d'existence du quotient d'un schéma affine par un groupe diagonolisable agissant librement.
\end{sousexemple}

\begin{sousprop}[SGA~3~VIII~2.1]
\label{prop_gpes_diagonalisables}
 Soient $S$ un schéma, $M$ un groupe commutatif ordinaire, et $G=D_S(M)$ le $S$-groupe diagonalisable défini par $M$.
\begin{itemize}
 \item[a)] Si $M$ est de type fini, alors $G$ est de présentation finie sur $S$.
\item[b)] Si $M$ est fini, alors $G$ est fini sur $S$.
\item[c)] Si $M$ est de torsion, alors $G$ est entier sur $S$. 
\end{itemize}
De plus, si $S$ est non vide, les réciproques sont vraies. Enfin dans le cas b), le degré de $G$ sur $S$ est égal au cardinal de $M$.
\end{sousprop}

\subsection{Conoyaux et espaces annelés quotients}
\label{espaces_anneles}

Soit $\Cc$ une catégorie (par exemple la catégorie des schémas). Si un objet en groupes $G$ de $\Cc$ agit sur un objet $X$, la première idée qui vient à l'esprit pour définir un quotient de $X$ par $G$ est de le définir par propriété universelle. On dit que $\pi : X \to Z$ est un \emph{quotient catégorique}, ou plus précisément un \emph{quotient dans la catégorie $\Cc$}, si c'est un morphisme invariant sous $G$ (ce qui se traduit par la commutativité d'un certain diagramme) et s'il est universel parmi les morphismes invariants sous $G$. Il revient au même de dire que $\pi$ est un conoyau du couple de flèches $\xymatrix@C=1pc{\mu, p_2 : G\times X \ar@<2pt>[r] \ar@<-2pt>[r] & X}$ où $\mu$ est l'action et $p_2$ la projection sur le second facteur.

 \begin{sousdefi}
Soit $\xymatrix@C=1pc{p_1, p_2 : X \ar@<2pt>[r] \ar@<-2pt>[r] & Y}$ un couple de flèches dans une catégorie $\Cc$. Un morphisme $u : Y\to Z$ est dit \emph{compatible} avec $(p_1, p_2)$ si $up_1=up_2$. C'est un \emph{conoyau} s'il est universel parmi les morphismes compatibles avec $(p_1, p_2)$, \emph{i.e.} si pour tout morphisme $v : Y \to T$ compatible avec $(p_1, p_2)$, il existe un unique $w : Z\to T$ tel que $v=wu$.
$$\xymatrix{X \ar@<2pt>[r]^{p_1} \ar@<-2pt>[r]_{p_2} & Y \ar[r]^u \ar[rd]_v&Z \ar@{.>}[d]^w \\ &&T}$$
\end{sousdefi}

\begin{sousremarque} \rm Le conoyau $Z$ ci-dessus représente naturellement un foncteur \emph{covariant}.
\end{sousremarque}

\begin{sousexemple}\rm
\label{fppf_epim_effectifs}
 Soit $f : X\to Y$ un morphisme de schémas qui forme une famille couvrante pour la topologie \emph{fpqc}, et soient $p_1, p_2$ les projections canoniques de $X\times_Y X$ vers~$X$. Alors $f$ est un conoyau de $(p_1,p_2)$ dans la catégorie des schémas. En effet, le fait que tout schéma $T$ soit un faisceau \emph{fpqc} donne immédiatement la propriété universelle du conoyau. En fait, cette propriété est vraie dès que $f$ est un morphisme \emph{couvrant} (voir exercice~\ref{epimorphismes_effectifs}).
\end{sousexemple}

Si $\xymatrix@C=1pc{p_1, p_2 : X \ar@<2pt>[r] \ar@<-2pt>[r] & Y}$ est une double flèche dans la catégorie des espaces annelés, elle a toujours un conoyau $Z$ dans cette catégorie. En effet, $Z$ a pour espace topologique sous-jacent le quotient de l'espace topologique sous-jacent à $Y$ obtenu en identifiant $p_1(x)$ et $p_2(x)$, pour tout~$x\in X$. On note $\pi : Y\to Z$ la projection canonique. Le faisceau d'anneaux sur $Z$ est construit de la manière suivante. Si $U$ est un ouvert de $Z$, $\Oc_Z(U)$ est le sous-anneau de $\Oc_Y(\pi^{-1}(U))$ formé des éléments $s$ tels que $p_1^{\sharp}(s)=p_2^{\sharp}(s)$ où $p_i^{\sharp}$ est le morphisme $\Oc_Y \to p_{i*}\Oc_X$ associé à $p_i$. 

\begin{sousexemple}\rm Au-dessus d'un corps $k$ de caractéristique différente de~2, on fait agir $\Z/2\Z$ sur la droite affine $\A^1_k$ par $x \mapsto -x$. On vérifie facilement que $\pi : \A_k^1 \to \A_k^1$ défini par $\pi(x)=x^2$ est un quotient au sens des espaces annelés. C'est aussi un quotient dans la catégorie des schémas.
\end{sousexemple}

\begin{sousprop}
\label{comparaison_conoyaux}
 Soit $\xymatrix@C=1pc{p_1, p_2 : X \ar@<2pt>[r] \ar@<-2pt>[r] & Y}$ une double flèche dans la catégorie des schémas, et soit $\pi : Y \to Z$ son conoyau dans la catégorie des espaces annelés. Si $Z$ est un schéma et $\pi$ un morphisme de schémas, alors $\pi$ est un conoyau dans la catégorie des schémas.
\end{sousprop}
\begin{demo}
 La courte preuve est sans surprise, voir~SGA~3~\cite{SGA3}~V~1.2.
\end{demo}

Pour construire le quotient d'un schéma sous l'action d'un groupe, l'espace annelé quotient semble donc \emph{a priori} être un bon candidat : il existe toujours, et si c'est un schéma alors c'est bien un quotient dans la catégorie des schémas. Cependant, ce point de vue présente au moins deux inconvénients. D'une part, l'espace annelé quotient n'est pas toujours un schéma (ceci peut arriver même s'il existe un conoyau dans la catégorie des schémas, voir exemple ci-dessous). Il ne sert alors pas à grand-chose en géométrie algébrique. D'autre part, même si ce quotient est un schéma, on ne sait pas décrire son foncteur des points. Le point de vue dominant dans la suite sera donc celui du \emph{faisceau quotient}, plus facile à décrire comme foncteur et dont on peut espérer qu'il soit représentable par un espace algébrique, sinon un schéma.

\begin{sousexemple}\rm  (\emph{cf.} \cite{Demazure_Gabriel}~III~\S 2~no~3.1)  Soit $k$ un corps algébriquement clos. On considère le groupe additif sous-jacent à $k$ (comme groupe ordinaire) et on note $G$ le groupe constant associé au-dessus de $\Spec k$. Pour tout $\gamma \in k$, on note $\rho_{\gamma}$ l'automorphisme de~$\A_k^1$ correspondant à l'automorphisme $x\mapsto x+\gamma$ de $k[x]$. Ceci définit une action de $G$ sur~$\A_k^1$, donc un couple de flèches $\xymatrix@C=1pc{p_1, p_2 : G\times_k \A_k^1 \ar@<2pt>[r] \ar@<-2pt>[r] &\A_k^1}$
où $p_1$ est l'action de $G$ et $p_2$ est la projection sur le second facteur. On vérifie facilement que l'espace topologique quotient est un ensemble à deux points muni de la topologie grossière. En particulier l'espace annelé quotient ne peut pas être un schéma. Par ailleurs, on démontre sans difficultés que le morphisme structural $\A_k^1 \to \Spec k$ est un conoyau de $(p_1,p_2)$ dans la catégorie des schémas.
\end{sousexemple}

\begin{sousexemple}\rm Soit $A$ un anneau de valuation discrète contenant $\C$, de corps des fractions $K$ et soit $\pi$ une uniformisante. On note $L$ le corps $K[T]/(T^n-\pi)$ et $A'$ la fermeture intégrale de $A$ dans $L$. L'anneau $A'=A[T]/(T^n-\pi)$ est encore de valuation discrète, d'uniformisante $T$. L'extension $L/K$ est galoisienne de groupe $G=\Z/n\Z$ et totalement ramifiée. Le groupe $G$ agit naturellement sur $A'$, donc sur son spectre $X'$, et on vérifie facilement que le quotient catégorique $X'/G$ s'identifie au spectre $X$ de l'anneau des invariants $A'^G=A$. Notons $\mgo'$ l'idéal maximal de $A'$. Il est naturellement muni d'une action de~$G$, compatible avec l'action sur $A'$, donc il définit un $\Oc_{X'}$-module $G$-équivariant sur~$X'$. Supposons que ce module $G$-équivariant provienne d'un $\O_X$-module. Il existerait alors un~$A$-module $M$ tel que $\mgo'$ soit isomorphe à $M\otimes_A A'$ en tant que $A'$-module avec action de~$G$. Le $A$-module $M$ serait nécessairement libre de rang 1 par descente si bien que $\mgo'$ serait engendré, comme $A'$-module, par $\mgo'^G$, donc par $\pi$, ce qui fournit une contradiction. On voit en particulier que la catégorie des modules quasi-cohérents sur le quotient $X'/G$ \emph{n'est pas équivalente} à la catégorie des $\Oc_{X'}$-modules quasi-cohérents $G$-équivariants.
\end{sousexemple}

\begin{sousexo}[épimorphismes effectifs]
\label{epimorphismes_effectifs} 
 Soit $f : X\to Y$ un morphisme de schémas. On dit que $f$ est un épimorphisme effectif s'il est égal au conoyau du couple de projections $\xymatrix@C=1pc{X\times_Y X\double &X}$. On a vu en~\ref{fppf_epim_effectifs} que si $f$ est une famille couvrante \emph{fpqc} c'est un épimorphisme effectif. \\
{\bf 1.} Montrer que si $f$ a une section $s$ alors $f$ est un épimorphisme effectif. [Indication : On vérifiera simplement la propriété universelle du conoyau. On pourra utiliser le morphisme de $X$ vers $X\times_Y X$ qui vaut l'identité sur le premier facteur et $sf$ sur le second.]\\
{\bf 2.} Montrer qu'un morphisme couvrant\footnote{Voir définition en~\ref{def_couvrants}.} pour la topologie \emph{fpqc} est un épimorphisme effectif.
\end{sousexo}

\subsection{Groupoïdes et relations d'équivalence}
\label{gpds}

\begin{sousdefi}
Soient $S$ un schéma et $X$ un foncteur contravariant de $(Sch/S)$ vers la catégorie des ensembles. Une relation d'équivalence sur $X$ est un sous-foncteur $R$ de $X\times_S X$ tel que pour tout $S$-schéma $U$, $R(U)$ soit le graphe d'une relation d'équivalence sur $X(U)$.
\end{sousdefi}

Remarquons qu'il revient au même de se donner l'inclusion de $R$ dans $X\times_S X$, ou de se donner le couple de flèches $\xymatrix@C=1pc{p_1, p_2 : R \ar@<2pt>[r] \ar@<-2pt>[r] &X}$ obtenues en composant l'inclusion avec les projections canoniques. Nous dirons donc parfois (par un léger abus de langage clairement innofensif) qu'une relation d'équivalence est la donnée d'un foncteur $R$ et d'un couple de $S$-morphismes $\xymatrix@C=1pc{p_1, p_2 : R \ar@<2pt>[r] \ar@<-2pt>[r] &X}$ tels que le morphisme $(p_1, p_2) : R \to X\times_S X$ soit un monomorphisme et que pour tout $S$-schéma $U$, l'image de $(p_1, p_2)(U)$ soit le graphe d'une relation d'équivalence sur $X(U)$.

\begin{sousremarque}\rm
\label{rqes_relation}
 Si $\xymatrix@C=1pc{p_1, p_2 : R \ar@<2pt>[r] \ar@<-2pt>[r] &X}$ est une relation d'équivalence, on voit facilement qu'il existe un automorphisme~$\sigma$ de $R$ tel que $\sigma^2=\id$ et $p_1\sigma=p_2$ (on envoie un couple $(x,y)\in X(T)\times X(T)$ sur $(y,x)$). De même, la diagonale de $X \to S$ induit un morphisme $\eps : X \to R$ qui est une section de $p_1$ et de $p_2$.
\end{sousremarque}

\begin{sousdefi}[changement de base]
\label{chgt_base}
Soit $\xymatrix@C=1pc{p_1, p_2 : R \ar@<2pt>[r] \ar@<-2pt>[r] &X}$ une relation d'équivalence sur un foncteur $X$ et soit $u : X \to Q$ un morphisme compatible avec $(p_1, p_2)$, \it{i.e.} $up_1=up_2$. Soit $\varphi : Q'\to Q$ un morphisme. Notons $R'=R\times_Q Q'$\, , $X'=X\times_Q Q'$ et $p_1', p_2'$ les morphismes obtenus à partir de $p_1, p_2$ par le changement de base $\varphi$. Alors $\xymatrix@C=1pc{p_1', p_2' : R' \ar@<2pt>[r] \ar@<-2pt>[r] &X'}$ est une relation d'équivalence sur $X'$. On dit que la relation d'équivalence $R'$ est \emph{obtenue à partir de $R$ par le changement de base} $\varphi$.  
 
\end{sousdefi}

En pratique, dans les cas qui nous intéressent, $X$ et $R$ sont au moins des espaces algébriques, voire des schémas. Dans ce cas, on appelle souvent \emph{quotient catégorique} de $X$ par $R$ le conoyau de $(p_1, p_2)$ dans la catégorie des schémas (s'il existe). On fera attention de ne pas confondre avec le \emph{faisceau quotient} de $X$ par $R$ dont nous parlerons au prochain paragraphe (chacun de ces deux objets est parfois appelé simplement  \og quotient\fg...). Il se trouve que, dans plusieurs cas utiles, ces deux objets coïncident (voir les théorèmes de représentabilité à partir de~\ref{thm_repres}). Dans la suite nous essaierons de préciser à chaque fois à quel objet nous faisons référence.

\begin{sousexemple}\rm
 Soit $G$ un $S$-schéma en groupes agissant sur $X$ (disons à gauche). On dit que l'action est \emph{libre} si pour tout $U$, le groupe $G(U)$ agit librement sur $X(U)$. De manière équivalente, une action est libre si et seulement si le morphisme induit
$$\begin{array}{rcl}
   G\times_S X& \lto &X\times_S X \\
(g,x) &\longmapsto& (x,gx)
  \end{array}
$$
est un monomorphisme. Il est immédiat que ce monomorphisme est alors une relation d'équivalence sur $X$.
\end{sousexemple}

Si l'action n'est pas libre, le morphisme ci-dessus n'est pas un monomorphisme. Dans ce cas la notion \emph{ad hoc} pour tenir compte des groupes d'inertie est celle de $S$-groupoïde. Comme le montre la remarque~\ref{gpd_vs_relation}, cette notion n'apporte rien pour étudier une action libre. Nous n'en aurons donc pas besoin car nous nous limiterons dans ce qui suit aux actions libres (à l'exception de quelques remarques dans le paragraphe~\ref{groupoides}). Nous donnons tout de même la définition ci-dessous à titre culturel. Rappelons qu'une catégorie $\Cc$ est un \emph{groupoïde} si tous ses morphismes sont des isomorphismes. 

\begin{sousdefi}
\label{def_groupoides}
Soit $S$ un schéma. Un $S$-groupoïde $X_*$ est la donnée de deux foncteurs contravariants $X_0$ et $X_1$ munis de $S$-morphismes \emph{source} $s : X_1 \to X_0$, \emph{but} $b : X_1 \to X_0$, \emph{neutre} $\eps : X_0 \to X_1$, \emph{composition} $m : X_1\times_{s,X_0,b} X_1 \to X_1$ et \emph{inverse} $i : X_1 \to X_1$ de telle sorte que pour tout $S$-schéma $U$, $X_*(T)$ soit un groupoïde dont l'ensemble des flèches est $X_1(T)$, l'ensemble des objets est $X_0(T)$, et les applications source, but, composition, inverse, identité sont données par $s(T), b(T), m(T), i(T), \eps(T)$.
\end{sousdefi}

Notons $X_2$ le produit fibré $X_1\times_{s,X_0,b} X_1$ et $p_0, p_2$ les projections sur le premier et le second facteur. On a alors un diagramme :
$$\xymatrix{
X_2 \ar@<4pt>[r]^{p_2} \ar@<-4pt>[r]_{p_0} \ar[r]|m &X_1  \ar@<2pt>[r]^s \ar@<-2pt>[r]_b &X_0\, .
}$$
En fait il revient au même (par un léger abus de langage) de se donner un groupoïde $X_*$, ou de se donner un diagramme comme ci-dessus, en exigeant que les trois carrés
$$\xymatrix{X_2 \ar[r]^{p_0} \ar[d]_{p_2} & X_1 \ar[d]_s \\
X_1 \ar[r]^b &X_0}\qquad \quad
\xymatrix{X_2 \ar[r]^m \ar[d]_{p_0} & X_1 \ar[d]_b \\
X_1 \ar[r]^b &X_0}\qquad \quad
\xymatrix{X_2 \ar[r]^m \ar[d]_{p_2} & X_1 \ar[d]_s \\
X_1 \ar[r]^s &X_0}$$
soient cartésiens, qu'un diagramme évident traduisant l'associativité de la composition $m$ soit commutatif, et qu'il existe une flèche $\eps : X_0\to X_1$ qui est à la fois une section de $s$ et une section de $b$ (voir SGA~3~\cite{SGA3}~V~1 pour plus de détails).

\begin{sousexemple}\rm
 Soient $X$ un $S$-schéma et $G$ un $S$-groupe agissant sur $X$ (à gauche pour fixer les idées). On note $\mu : G\times_S X \to X$ l'action de $G$. On appelle \emph{groupoïde associé à l'action de $G$ sur $X$} le groupoïde
$$\xymatrix{
G\times_S G\times_S X \ar@<4pt>[r]^-{\pr_{23}} \ar@<-4pt>[r]_-{id_G\times\mu} \ar[r]|-m &G\times_S X  \ar@<2pt>[r]^-{\pr_2} \ar@<-2pt>[r]_-{\mu} &X\, .
}$$
où la composition $m$ est le morphisme $\mu\times id_X$ ($\pr_2$ et $\pr_{23}$ désignant les projections évidentes). 
\end{sousexemple}

\begin{sousremarque}\rm
\label{gpd_vs_relation}
 Soit $X_*$ un groupoïde. Si $b\times s  : X_1\to X_0\times X_0$ est un monomorphisme, alors le couple de flèches $(b, s)$ est une relation d'équivalence. Réciproquement, si on part d'une relation d'équivalence $\xymatrix@C=1pc{b, s : X_1 \ar@<2pt>[r] \ar@<-2pt>[r]& X_0}$, on pose $X_2:=X_1\times_{b,X_0,s}X_1$. Il existe alors une unique flèche $m : X_2 \to X_1$ faisant de $X_*$ un groupoïde. (Exercice : le démontrer, et décrire $m$. La réponse se trouve dans SGA~3 V \S 2 b)). Il revient donc au même de se donner une relation d'équivalence ou un groupoïde dans lequel $b\times s$ est un monomorphisme.
\end{sousremarque}

Nous définissons ci-dessous l'image réciproque d'une relation d'équivalence $\xymatrix@C=1pc{R\double &X}$ par un morphisme $X'\to X$, puis nous en étudions quelques propriétés dans les exercices qui suivent. Cette notion servira essentiellement dans la preuve du lemme~\ref{cas_quasi_section} et le lecteur peut éventuellement omettre la fin de ce paragraphe en première lecture. 

\begin{sousdefi}
 \label{image_reciproque}
Si $(p_1,p_2) : R\to X\times_S X$ est une relation d'équivalence et si $f : X'\to X$ est un $S$-morphisme, on définit une relation d'équivalence $R'$ sur $X'$ en choisissant $R'$ égal au produit fibré suivant.
$$\xymatrix{R'\ar[r]\ar[d]_{(p'_1, p'_2)} &R \ar[d]^{(p_1, p_2)} \\ X'\times_S X' \ar[r]^{f\times f}& X\times_S X}$$
On dit que $R'$ est \emph{l'image réciproque} de $R$ par le morphisme $f$.
\end{sousdefi}

\begin{sousremarque}\rm
 On prendra garde de ne pas confondre cette notion d'image réciproque avec la notion de \emph{changement de base} (d'ailleurs plus utile) de~\ref{chgt_base}. En particulier, avec les notations ci-dessus, si les diagrammes 
$$\xymatrix{R'\ar[r] \ar[d]_{p'_i} & R \ar[d]^{p_i} \\
X' \ar[r] & X}$$
sont toujours commutatifs, il ne sont en général pas cartésiens.
\end{sousremarque}

\begin{sousexo}\rm
\label{images_reciproques_identiques}
 Soient $X$ un préfaisceau sur $(Sch/S)^o$ et $R$ un sous-foncteur de $X\times_S X$ qui définit une relation d'équivalence $\xymatrix@C=1pc{p_1, p_2 : R \double & X}$. Montrer que la relation d'équivalence $R' \to R\times_S R$, image réciproque de $R$ par le morphisme $p_1$, coïncide, comme sous-foncteur de $R\times_S R$, avec l'image réciproque de $R$ par $p_2$. [Indication : Les deux sont égales au sous-foncteur de $R\times_S R$ des couples $((x_1, x_2), (x_3,x_4))$ tels que les $x_i$ soient tous équivalents modulo $R$.]
\end{sousexo}

\begin{sousexo}\rm
 \label{descente_fpqc_existence_conoyau} Soit $\xymatrix@C=1pc{p_1, p_2 : R \double & X}$ une relation d'équivalence où $X$ et $R$ sont des schémas. Soit $f : X' \to X$ un morphisme couvrant (\emph{cf.}~\ref{def_couvrants}) pour la topologie \emph{fpqc}. On note $\xymatrix@C=1pc{p'_1, p'_2 : R' \double & X'}$ la relation d'équivalence sur $X'$ image réciproque de $R$ par $f$.\\
{\bf 1.} Soit $\pi : X \to Q$ un conoyau de $(p_1,p_2)$ dans $(Sch/S)$. Montrer que $\pi f$ est un conoyau de $(p'_1,p'_2)$. [Indication : En notant $X''=X'\times_X X'$ et $q_1, q_2$ ses projections sur $X'$, construire $\Delta$ de $X''$ vers $R'$ tel que $p'_i\Delta=q_i$, puis jouer avec les propriétés universelles.]\\
{\bf 2.} Réciproquement, soit $\pi' : X'\to Q$ un conoyau de $(p'_1,p'_2)$. Montrer qu'il existe un unique morphisme $\pi$ de $X$ vers $Q$ tel que $\pi'=\pi\circ f$, puis que $\pi$ est un conoyau de~$(p_1,p_2)$.
\end{sousexo}

\begin{sousexo}\rm
 Reprenons les notations de~\ref{chgt_base} et notons $f : X'\to X$ le morphisme induit par $\varphi$.\\
{\bf 1.} Si $\varphi$ est un monomorphisme, montrer que la relation d'équivalence $R'$ sur $X'$ obtenue par le changement de base $\varphi$ coïncide avec l'image réciproque (\ref{image_reciproque}) de $R$ par $f$.\\
{\bf 2.} Donner un contre-exemple dans le cas général.
\end{sousexo}

\subsection{Faisceau quotient}
\label{faisceau_quotient}

\begin{sousdefi}
Soient $S$ un schéma et soit $\xymatrix@C=1pc{p_1, p_2 : R \ar@<2pt>[r] \ar@<-2pt>[r]& X}$ une relation d'équivalence où $R$ et $X$ sont des schémas (voire des espaces algébriques). Soit (Top) une topologie sur $(Sch/S)^o$. On appelle \emph{faisceau quotient (Top)} de $X$ par $R$ le faisceau (Top) associé au préfaisceau qui à $U$ associe l'ensemble quotient $X(U)/R(U)$. En général on choisira la topologie \emph{fppf}, et lorsque l'on parlera \emph{du} faisceau quotient, il s'agira du faisceau quotient \emph{fppf}.
\end{sousdefi}

\begin{sousremarque}\rm
 On vérifie facilement avec la propriété universelle du faisceau associé à un préfaisceau que le faisceau quotient (Top) est en fait le conoyau du couple de flèches $\xymatrix@C=1pc{R \double & X}$ dans la catégorie des faisceaux (Top). En particulier, si $R$, $X$ et $X/R$ sont des schémas, alors $X/R$ est le conoyau de $\xymatrix@C=1pc{R \double & X}$ dans la catégorie des schémas (donc un \emph{quotient catégorique} au sens du paragraphe~\ref{espaces_anneles}).
\end{sousremarque}

\begin{sousremarque}\rm
 Plus généralement, si $X_*$ est un groupoïde, on peut lui associer un faisceau quotient. En effet, le sous-foncteur $R$ de $X_0\times_S X_0$ défini en prenant pour $R(T)$ l'image du morphisme $X_1(T) \to (X_0\times_S X_0)(T)$ est une relation d'équivalence sur $X_0$ (mais en général $R$ n'est pas représentable). On définit alors un faisceau quotient comme ci-dessus. Ici aussi, on vérifie que le faisceau quotient~(Top) est le conoyau du couple de flèches $\xymatrix@C=1pc{X_1\double & X_0}$ dans la catégorie des faisceaux~(Top). Pour un $S$-schéma en groupes~$G$ agissant sur un $S$-schéma $X$, le faisceau quotient est donc le faisceau \emph{fppf} associé au préfaisceau qui à $U$ associe l'ensemble quotient $X(U)/G(U)$.
\end{sousremarque}

\begin{souslem}
 \label{lem_faisceaux}
Soit $X\to Y$ un morphisme de faisceaux sur un site $\Cc$. Les propositions suivantes sont équivalentes :
\begin{itemize}
 \item[(i)] $X\to Y$ est \og localement surjectif\fg\ (\emph{i.e.} toute section de $Y$ provient localement de $X$) ;
\item[(ii)] $X\to Y$ est un épimorphisme dans la catégorie $\Sc(\Cc)$ des faisceaux sur $\Cc$ ;
\item[(iii)] $X\to Y$ est un épimorphisme effectif dans la catégorie $\Sc(\Cc)$, autrement dit c'est un conoyau du couple de flèches $\xymatrix@C=1pc{X\times_Y X \double & X}$.
\end{itemize}
\end{souslem}
\begin{demo}
L'équivalence entre (i) et (ii) est facile et laissée en exercice au lecteur. Par ailleurs il est évident que (iii) implique (ii). Réciproquement supposons que $X\to Y$ soit couvrant et montrons (iii). On note $G$ le sous-préfaisceau de $Y$ image de $X$. Comme $X\to Y$ est localement surjectif, on voit que $Y$ s'identifie au faisceau $aG$ associé à $G$. On a donc une suite de morphismes de préfaisceaux :
$$\xymatrix{X\ar[r] & G \ar[r]& aG \ar[r]^{\sim}& Y.}$$
Par construction de $G$, pour tout $S$-schéma $U$ le morphisme $X(U)\to G(U)$ est surjectif, autrement dit c'est un conoyau (dans la catégorie des ensembles) du couple de flèches $\xymatrix@C=1pc{X(U)\times_{G(U)} X(U) \double & X(U)}$. Ceci montre que $X\to G$ est un conoyau, dans la catégorie des préfaisceaux, du couple de flèches $\xymatrix@C=1pc{X\times_G X \double & X}$. Par ailleurs, comme $G\to Y$ est un monomorphisme, le produit fibré $X\times_G X$ s'identifie à $X\times_Y X$, si bien que $X\to G$ est le conoyau, dans la catégorie des préfaisceaux, du couple de flèches $\xymatrix@C=1pc{X\times_Y X \double & X}$. Le faisceau associé à $G$ est donc le conoyau de ce même couple de flèches dans la catégorie des faisceaux (conséquence immédiate des propriétés universelles), cqfd.
\end{demo}

\begin{sousprop}
\label{prop_faisceau_quotient}
 Soient $S$ un schéma et $\xymatrix@C=1pc{R \ar@<2pt>[r]^{p_1} \ar@<-2pt>[r]_{p_2}& X}$ une relation d'équivalence avec $R, X$ des faisceaux (Top) sur $(Sch/S)$, où (Top) est l'une des topologies \emph{fpqc}, \emph{fppf}, étale ou Zariski. On note $Q$ le faisceau quotient (Top) de $X$ par $R$ (on suppose son existence dans le cas \emph{fpqc}) et $\pi : X\to Q$ la projection canonique.
\begin{itemize}
 \item[(i)] Le morphisme $\pi$ est un épimorphisme dans la catégorie des faisceaux~(Top).
\item[(ii)] Le carré
$$\xymatrix{R\ar[r]^{p_1} \ar[d]_{p_2}&X \ar[d]^{\pi} \\ X\ar[r]_{\pi} &Q}$$
est cartésien, autrement dit le morphisme naturel $(p_1, p_2) : R \to X\times_Q X$ est un isomorphisme.
\item[(iii)] Si $Q'\to Q$ est un morphisme de faisceaux (Top), alors $\pi' : X\times_Q Q' \to Q'$ est le faisceau quotient (Top) de la relation d'équivalence $\xymatrix@C=1pc{R\times_Q Q' \ar@<2pt>[r]^{p'_1} \ar@<-2pt>[r]_{p'_2}& X\times_Q Q'}$. Autrement dit, la formation du faisceau quotient commute au changement de base.
\end{itemize}
\end{sousprop}
\begin{demo}
Le point (i) est clair, car un conoyau est toujours un épimorphisme. Le carré de (ii) induit un morphisme $R\to X\times_Q X$. C'est un monomorphisme puisque le composé $R\to X\times_Q X\to X\times_S X$ en est un. Maintenant, soit $U$ un $S$-schéma et soit $(x_1,x_2)\in (X\times_QX)(U)$. Alors les éléments $x_1, x_2$ de $X(U)$ ont la même image dans $Q(U)$. Il existe donc une famille couvrante $U'\to U$ telle que ${x_1}_{|_{U'}}$ et ${x_2}_{|_{U'}}$ aient la même image dans le préfaisceau quotient $X/R$ (car $Q$ est le faisceau associé à ce préfaisceau). On en déduit que $({x_1}_{|_{U'}},{x_2}_{|_{U'}})\in R(U')$ puis, comme $R$ est un faisceau, que $(x_1, x_2)\in R(U)$, ce qui achève de prouver (ii). Le point (iii) est facile et laissé en exercice au lecteur.
\end{demo}

\begin{sousprop}
\label{comparaison_conoyau_faisceau}
 Soit $\xymatrix@C=1pc{R \ar@<2pt>[r]^{p_1} \ar@<-2pt>[r]_{p_2}& X}$ une relation d'équivalence avec $R, X$ des schémas sur une base $S$ fixée. Les propriétés suivantes sont équivalentes :
\begin{itemize}
 \item[(i)] il existe un schéma qui représente le faisceau quotient $X/R$ au sens \emph{fpqc} (resp. \emph{fppf}, resp. étale) ;
\item[(ii)] la relation d'équivalence admet un conoyau $\pi : X\to Q$ dans la catégorie des schémas, le morphisme naturel $(p_1, p_2) : R \to X\times_Q X$ est un isomorphisme, et le morphisme $\pi$ est couvrant\footnote{On rappelle qu'un morphisme de schémas $X\to Y$ est dit couvrant pour une topologie (sous-canonique) s'il existe des morphismes $U_i\to X$ tels que les composés $U_i\to Y$ forment une famille couvrante de $Y$ pour la topologie considérée. Cela revient à dire que c'est un épimorphisme dans la catégorie des faisceaux pour cette topologie.} pour la topologie \emph{fpqc} (resp. \emph{fppf}, resp. étale).
\end{itemize}
S'il en est ainsi, $Q$ représente le faisceau $X/R$.
\end{sousprop}

\begin{sousremarque}\rm
 L'énoncé est encore valable, avec la même preuve, en remplaçant partout le mot \og schéma\fg\ par \og espace algébrique\fg\ (\emph{cf.} plus bas pour la définition).
\end{sousremarque}
\begin{demo}
Le fait que (i) implique (ii) résulte de ce qui précède. Réciproquement, supposons (ii). Comme $\pi : X\to Q$ est couvrant, le lemme~\ref{lem_faisceaux} montre que c'est un conoyau de $\xymatrix@C=1pc{R=X\times_Q X \double & X}$ dans la catégorie des faisceaux. Le schéma $Q$ représente donc le faisceau quotient $X/R$. 
\end{demo}

\noindent{\sc Comparaison étale versus \emph{fppf}}

Dans le cas d'une relation d'équivalence donnée par l'action d'un groupe \emph{lisse et de présentation finie}, en admettant les résultats (difficiles) d'Artin de représentabilité par des espaces algébriques (voir plus bas), on déduit de la proposition~\ref{comparaison_conoyau_faisceau} un théorème de comparaison du faisceau quotient étale et du faisceau quotient \emph{fppf}.

\begin{sousthm}
 Soient $S$ un schéma, et $G$ un $S$-schéma en groupes qui agit (disons à droite) librement sur un $S$-schéma~$X$ quasi-séparé. On suppose que le morphisme structural $G\to S$ est lisse et de présentation finie. Alors le faisceau quotient étale $X/G$ est représentable par un espace algébrique. En particulier il coïncide avec le faisceau quotient au sens \emph{fppf}.
\end{sousthm}
\begin{demo}
 D'après~\ref{repres_par_esp_alg}, on sait que le faisceau quotient (\emph{fppf}) est représentable par un espace algébrique $Q$. Notons $\pi : X \to Q$ le morphisme de passage au quotient, $\mu : X\times_S G \to X$ l'action de $G$ sur $X$ et $\pr_1 :  X\times_S G \to X$ la projection sur le premier facteur. D'après~\ref{comparaison_conoyau_faisceau}, le carré
$$\xymatrix{X\times_S G \ar[r]^{\mu} \ar[d]_{\pr_1}&X \ar[d]^{\pi} \\ X\ar[r]_{\pi} &Q}$$
est cartésien et $\pi$ est couvrant pour la topologie \emph{fppf}. Or $\pr_1$ est lisse et surjectif. Par descente fidèlement plate, on en déduit que $\pi$ est lui-même lisse et surjectif, donc couvrant pour la topologie étale par~\ref{sections_morphismes_lisses}. La proposition~\ref{comparaison_conoyau_faisceau} montre alors que le quotient au sens étale est déjà un espace algébrique.
\end{demo}

\begin{sousexo}\rm
\label{exercice_pptes_1}
 Soit $G$ un schéma en groupes agissant librement sur un schéma $X$ au-dessus d'une base $S$ fixée. On suppose que le faisceau quotient $X/G$ (au sens \emph{fppf}) est représentable par un $S$-schéma $Q$ et on note $\pi : X\to Q$ le morphisme quotient. Montrer les assertions suivantes.
\begin{itemize}
 \item[(1)] Soit $S' \to S$ un morphisme de schémas. Posons $G'=G\times_S S'$ et $X'=X\times_S S'$. Alors $X'/G'$ est représentable par $Q\times_S S'$.
\item[(2)] Si $X$ est réduit, alors $Q$ aussi.
\item[(3)] Le monomorphisme $(\mu, \pr_2) : G\times_S X \to X\times_S X$ est une immersion.
\item[(4)] Pour que $Q$ soit séparé sur $S$, il faut et il suffit que $(\mu, \pr_2)$ soit une immersion fermée.
\item[(5)] Le morphisme $\pr_2 : G\times_S X\to X$ est plat et localement de présentation finie si et seulement si $\pi$ l'est. Sous ces conditions (c'est notamment le cas si $G$ lui-même est plat et localement de présentation finie sur $S$), si $X$ est localement de type fini (resp. de type fini, resp. plat, resp. lisse, resp. étale, resp. net, resp. localement quasi-fini, resp. quasi-fini) sur $S$, il en est de même de $Q$.
\end{itemize}
\end{sousexo}

\begin{sousexo}[\emph{cf.} SGA 3 \cite{SGA3} VI B \S 9]\rm
 Soit $u : G\to H$ un monomorphisme de $S$-schémas en groupes. On fait agir $G$ sur $H$ (librement) par translations à droite. On suppose que le faisceau quotient $H/G$ (au sens \emph{fppf}) est représentable par un $S$-schéma $Q$ et on note $\pi : H\to Q$ le morphisme quotient. Montrer les assertions suivantes.
\begin{itemize}
 \item[(1)] On note $\eps : S \to H$ la section unité de $H$ et $\eps_Q : S\to Q$ son composé avec $\pi$. On appelle $\eps_Q$ la \og section unité\fg\ de $Q$. Montrer que le diagramme
$$\xymatrix{G\ar[r]^{u} \ar[d]&H \ar[d]^{\pi} \\ S\ar[r]_{\eps_Q} &Q}$$
est cartésien. En particulier, $u$ est une immersion (car $\eps_Q$ en est une, puisque c'est une section de $Q$).
\item[(2)] $H$ agit à gauche sur $Q$, et le morphisme $\pi$ est compatible avec cette action et avec l'action de $H$ sur lui-même par translations à gauche.
\item[(3)] Si $G$ est invariant dans $H$, il existe sur $Q$ une unique structure de $S$-groupe qui fait de $\pi$ un morphisme de $S$-groupes.
\item[(4)] Pour que $Q$ soit séparé sur $S$, il faut et il suffit que $u$ soit une immersion fermée, ou encore que $\eps_Q$ soit une immersion fermée.
\end{itemize}
\end{sousexo}

Nous terminons cette section par l'énoncé d'un théorème qui permet de se débarrasser des nilpotents de la base dans les questions de représentabilité du faisceau quotient.

\begin{sousthm}[\cite{Demazure_Gabriel}~III~\S 2, no 7, thm 7.1]\ 
\label{representabilite_et_nilpotents}

\noindent
 Soit $S$ un schéma de base et soit $\xymatrix@C=1pc{p_1, p_2 : R \ar@<2pt>[r] \ar@<-2pt>[r]& X}$ une relation d'équivalence où $R$ et $X$ sont des schémas et où les $p_i$ sont fidèlement plats et de présentation finie. Soit $S_0$ un sous-schéma fermé de $S$ défini par un idéal nilpotent. On note $X_0=X\times_S S_0$ et $R_0$ la relation induite par $R$ sur $X_0$. Si le faisceau quotient $X_0/R_0$ est représentable par un schéma, alors $X/R$ l'est aussi.
\end{sousthm}

\subsection{Passage au quotient par une relation d'équivalence}
\label{thm_repres}

Nous donnons dans ce paragraphe les principaux résultats de SGA~3~\cite{SGA3}~V, avec des esquisses de preuves. Nous renvoyons à \emph{loc. cit.} pour des preuves complètes. Nous avons préféré nous limiter aux relations d'équivalence, cas dans lequel les résultats sont plus forts (représentabilité du faisceau quotient, contre la seule existence d'un quotient catégorique autrement) et les preuves un peu plus simples. Les résultats présentés ici ont tous des analogues dans le cas des groupoïdes, que nous évoquerons dans la section~\ref{groupoides}.

\begin{sousthm}[cas fini et localement libre]
\label{cas_fini}
 Soit $(p_1,p_2) : \xymatrix@C=1pc{R \double &X}$ une relation d'équivalence, avec $X, R$ des schémas. On suppose que
\begin{itemize}
 \item[a)] $p_1$ est fini et localement libre (alors $p_2$ l'est aussi) ;
\item[b)] pour tout $x\in X$, $p_1p_2^{-1}(x)$ est contenu dans un ouvert affine de $X$.
\end{itemize}
Alors :
\begin{itemize}
 \item[(i)] Il existe un morphisme $\pi : X \to Q$ qui est un conoyau de $(p_1,p_2)$ dans $(Sch/S)$. De plus $\pi$ est un conoyau dans la catégorie \EspAnn\ de tous les espaces annelés.
\item[(ii)] $\pi$ est fini et localement libre.
\item[(iii)] Le morphisme $R \to X\times_Q X$ de composantes $p_1$ et $p_2$ est un isomorphisme.
\item[(iv)] $Q$ représente le faisceau quotient \emph{fppf} de $X$ par la relation d'équivalence $(p_1, p_2)$.
\item[(v)] Pour tout morphisme $Q'\to Q$, $Q'$ est le conoyau du couple de flèches $(p_1', p'_2)$ déduit de $(p_1,p_2)$ par le changement de base $Q' \to Q$. Autrement dit, \og la formation du quotient commute au changement de base\fg.
\end{itemize}
\end{sousthm}
\begin{demo}
 On rappelle (\ref{rqes_relation}) qu'il existe une involution $\sigma$ de $R$ telle que $p_1\sigma=p_2$ (donc $p_2\sigma=p_1$). En particulier, on voit que $p_1$ est fini et localement libre si et seulement si $p_2$ l'est. Par ailleurs, l'ensemble $p_1p_2^{-1}(x)$ est aussi $p_1\sigma\sigma^{-1}p_2^{-1}(x)=p_2p_1^{-1}(x)$, donc les hypothèses sont en réalité symétriques en $p_1$ et $p_2$. Il suffit de montrer (i), (ii) et (iii). Les assertions (iv) et (v) en sont des conséquences par~\ref{comparaison_conoyau_faisceau} et~\ref{prop_faisceau_quotient}.

\begin{etape}{Cas où $X$ est affine et $p_1$ localement libre de rang constant $n$.}
 $R$ est alors affine aussi puisque $p_1$ est fini. On note $X=\Spec A$ et $\delta_i : A \to \Oc(R)$ le morphisme d'anneaux correspondant à $p_i$. On note $B$ le sous-anneau de $A$ formé des éléments $a\in A$ tels que $\delta_1(a)=\delta_2(a)$. On note enfin $Q=\Spec B$ et $\pi : X \to Q$ le morphisme induit par l'inclusion de $B$ dans $A$. Nous allons montrer que $\pi$ convient. En utilisant les propriétés d'une relation d'équivalence, on montre sans trop de difficultés (voir SGA~3 pour les détails) que pour $a\in A$, le polynôme caractéristique de $\delta_1(a)$ lorsqu'on considère $\Oc(R)$ comme algèbre sur $A$ \emph{via} $\delta_2$ est à coefficients dans $B$ et qu'il annule $a$. Ceci prouve d'une part que $A$ est entier sur $B$, et d'autre part 
que pour $a\in A$, la norme $N_{p_2}(\delta_1(a))$ est dans $B$.

Assertion : si deux points $x, y\in X$ ont même image dans $Q$, alors il existe $z\in R$ tel que $p_1(z)=x$ et $p_2(z)=y$. On raisonne par l'absurde. Dans le cas contraire, pour $z\in p_2^{-1}(y)$ on aurait $p_1(z)\neq x$. Par ailleurs, on a aussi $\pi p_1(z)=\pi p_2(z)=\pi(y)=\pi(x)$. Comme $\pi$ est entier, ceci implique que $p_1(z)\notin \ov{\{x\}}$ (par Cohen-Seidenberg, si $x_0\in \ov{\{x_1\}}$ et $\pi(x_0)=\pi(x_1)$ alors $x_0=x_1$). On a donc montré que l'ensemble $p_1p_2^{-1}(y)$ ne rencontre pas $\ov{\{x\}}$. Comme $X$ est affine et $p_1p_2^{-1}(y)$ fini, il existe alors une fonction $a\in A$ qui s'annule en $x$ mais pas aux points de $p_1p_2^{-1}(y)$ (voir exercice~\ref{lemme_affine} ci-dessous). Alors la fonction $\delta_1(a)$ s'annule sur la fibre $p_1^{-1}(x)$, mais en aucun point de $p_2^{-1}(y)$. Notant $Z(\delta_1(a))$ le lieu des zéros de cette fonction, on a $p_2(Z(\delta_1(a)))=Z(N_{p_2}(\delta_1(a)))$ (voir exercice~\ref{norme}). Alors $N_{p_2}(\delta_1(a))$ s'annule sur $p_2(p_1^{-1}(x))$, donc en $x$, mais pas en $y$. Ceci est absurde car $N_{p_2}(\delta_1(a))\in B$ et $x$ et $y$ ont la même image dans $Q$.

Considérons le morphisme $R \to X\times_Q X$ de composantes $p_1$ et $p_2$. Il est immédiat que c'est un morphisme fini (en composant avec $\pr_1$ par exemple, on obtient $p_1$ qui est fini). De plus c'est un monomorphisme puisque $(p_1, p_2)$ est une relation d'équivalence. C'est donc une immersion fermée. Par des arguments d'algèbre un peu longs pour être reproduits ici (voir SGA~3), on montre que c'est même un isomorphisme, et au passage que $\pi$ est fini et localement libre de rang $n$ (tout ceci est admis ici), donc surjectif.

Nous pouvons maintenant conclure. Vu l'assertion ci-dessus, on voit que l'ensemble sous-jacent à $Q$ est le quotient de l'ensemble sous-jacent à $X$ par la relation d'équivalence définie par $(p_1,p_2)$. De plus, comme $\pi$ est fini et localement libre, il est surjectif et ouvert, donc submersif, si bien que la topologie de $Q$ est la topologie quotient de celle de $X$. Enfin, il est clair vu la définition de $B$ que le faisceau structural $\Oc_Q$ de $Q$ est formé des sections~$s$ de $\pi_*\Oc_{X}$ telles que $\delta_1(s)=\delta_2(s)$. Il en résulte que $\pi : X \to Q$ est le conoyau dans la catégorie des espaces annelés (voir la construction de ce conoyau en~\ref{espaces_anneles}). Mais comme de plus $\pi$ est un morphisme de schémas, c'est un conoyau dans la catégorie des schémas par~\ref{comparaison_conoyaux}.

\end{etape}
 
\begin{etape}{Cas où $p_1$ est localement libre de rang constant $n$.}
 On montre d'abord l'assertion suivante.

\emph{Assertion : Tout point $x\in X$ possède un voisinage ouvert affine et saturé.}

\noindent Par hypothèse, il existe $V$ ouvert affine contenant la classe d'équivalence $p_2p_1^{-1}(x)$. On note $V'$ la réunion des classes d'équivalences incluses dans $V$. C'est un ouvert de $X$ car c'est le complémentaire de $p_2p_1^{-1}(X\setminus V)$ qui est fermé puisque $p_2$ est fini. De plus il est saturé par construction et c'est le plus grand ouvert saturé inclus dans $V$. Comme $V$ est affine et comme l'ensemble $p_2p_1^{-1}(x)$ est fini et contenu dans $V'$, il existe une fonction $f\in \Gamma(V,\Oc_V)$ qui s'annule sur $V\setminus V'$ mais pas aux points de $p_2p_1^{-1}(x)$. Autrement dit l'ouvert principal $D(f)$ est inclus dans $V'$ et contient $p_2p_1^{-1}(x)$. On note $V''$ la réunion des classes d'équivalences incluses dans $D(f)$. Comme précédemment c'est un ouvert saturé de $X$. Il est contenu dans $D(f)$ et contient $p_2p_1^{-1}(x)$. Il suffit de montrer qu'il est affine. On note $Z(f)$ le lieu d'annulation de $f$ dans $V'$. Alors $p_1^{-1}(Z(f))$ est le lieu d'annulation de $p_1^*(f)$ dans $f^{-1}(V')$. Donc $p_2p_1^{-1}(Z(f))$ est le lieu d'annulation de $N_{p_2}(p_1^*(f))$ dans $V'$ (exercice~\ref{norme}). Or son complémentaire dans $V'$ est précisément $V''$ (par construction). $V''$ est donc l'ensemble des points de $D(f)$ où $N_{p_2}(p_1^*(f))$ ne s'annule pas, ce qui prouve qu'il est affine et achève la preuve de notre assertion.

Montrons maintenant (i), (ii) et (iii). Soit $\pi : X \to Q$ le conoyau de $(p_1,p_2)$ dans la catégorie des espaces annelés. Soit $x\in X$ et soit $U_0$ un voisinage ouvert affine et saturé de $x$. On note $U_1=p_1^{-1}(U_0)=p_2^{-1}(U_0)$. Alors $(p_1, p_2) : \xymatrix@C=1pc{U_1\double &U_0}$ est une relation d'équivalence dont $V=\pi(U_0)$ (qui est un ouvert de $Q$) est le conoyau dans \EspAnn. D'après le cas affine, $V$ est un schéma. Vu que de tels $V$ recouvrent $Q$, on voit que $\pi$ est un morphisme de schémas, puis qu'il vérifie les conclusions (i), (ii) et (iii), locales sur $Q$.
\end{etape}

\begin{etapefinale}{Cas général.}
 Si $x\in X$, on note $\rg_x(p_1)$ le rang de la fibre $p_1^{-1}(x)$. Pour $n\in \N$, on note alors $U^n$ l'ensemble des points de $X$ tels que $\rg_x(p_1)=n$. Comme $p_1$ est fini et localement libre, $X$ est la somme disjointe des $U^n$. Par ailleurs, en utilisant l'associativité de la relation d'équivalence, on montre que $p_1^{-1}(U^n)=p_2^{-1}(U^n)$, autrement dit, $U^n$ est \emph{saturé}. On note $V^n=p_1^{-1}(U^n)$. Le couple d'équivalence $(p_1,p_2) : \xymatrix@C=1pc{R \double &X}$ est alors la somme disjointe des couples d'équivalence $(p_1,p_2) : \xymatrix@C=1pc{V^n \double &U^n}$, ce qui nous ramène au cas précédent et termine la preuve.
\end{etapefinale}
\end{demo}

\begin{sousexo}[lemme d'évitement]\rm
\label{lemme_affine}
 Soit $X$ un schéma affine. Montrer que si $F$ est un fermé de $X$ et si $x_1, \dots, x_n$ sont des points de $X\setminus F$, il existe une fonction $f\in \Gamma(X, \Oc_X)$ qui s'annule sur~$F$ mais pas aux points $x_i$.
\end{sousexo}

\begin{sousexo}[norme]\rm
 \label{norme}
Soit $f : X\to S$ un morphisme de schémas fini et localement libre de rang constant $n$. Soit $b\in \Gamma(X,\Oc_X)$. Soit $S_i$ un ouvert affine de $S$ au-dessus duquel $f$ est libre de rang $n$. La multiplication par $b_{|_{f^{-1}(S_i)}}$ est un endomorphisme $\Oc_S(S_i)$-linéaire de $\Oc_X(f^{-1}(S_i))$. On note $s_i\in \Oc_S(S_i)$ son déterminant. Il est clair que les $s_i$ se recollent en une section de $\Gamma(S, \Oc_S)$ que l'on appelle la \emph{norme} de $b$ relativement à $f$ et que l'on note $N_f(b)$.

1) Montrer que la formation de $N_f(b)$ commute au changement de base. 

2) Montrer que $b\in \Gamma(X,\Oc_X)^{\times}$ si et seulement si $N_f(b)\in \Gamma(S,\Oc_S)^{\times}$.

3) Montrer que si $Z(b)$ désigne le lieu des zéros de $b$ dans $X$, alors $f(Z(b))$ est le lieu des zéros de $N_f(b)$ dans $S$. 
\end{sousexo}

Les autres résultats majeurs de cet exposé~V de SGA~3 sont l'existence du quotient dans le cas propre et plat (théorème~\ref{cas_propre_et_plat}) et l'existence du quotient \emph{génériquement} dans le cas d'une relation d'équivalence seulement plate (théorème~\ref{cas_plat}). Ils se ramènent tous deux au cas fini et localement libre vu ci-dessus par des techniques de quasi-section.

\begin{sousdefi}
 Soit $(p_1,p_2) : \xymatrix@C=1pc{R \double &X}$ une relation d'équivalence. Une \emph{quasi-section} de $(p_1, p_2)$ est un sous-schéma $U\to X$ qui vérifie les deux conditions suivantes.
\begin{itemize}
 \item[(1)] La restriction $v$ de $p_2$ à $p_1^{-1}(U)$ est un morphisme fini, localement libre et surjectif de $p_1^{-1}(U)$ sur $X$.
\item[(2)] Pour tout $x\in U$, l'ensemble $p_1p_2^{-1}(x)\cap U$ (fini d'après (1)) est contenu dans un ouvert affine de $U$.
\end{itemize}
\end{sousdefi}

\begin{souslem}[passage au quotient en présence d'une quasi-section]\ 
\label{cas_quasi_section}

\noindent Soit $(p_1,p_2) : \xymatrix@C=1pc{R \double &X}$ une relation d'équivalence qui possède une quasi-section. Alors :
\begin{itemize}
 \item[(i)] Il existe un morphisme $\pi : X \to Q$ qui est un conoyau de $(p_1,p_2)$ dans $(Sch/S)$. De plus $\pi$ est un conoyau dans la catégorie \EspAnn\ de tous les espaces annelés.
\item[(i')] $\pi$ est surjectif. De plus, si $p_1$ est ouvert (resp. universellement fermé, resp. plat) alors $\pi$ l'est aussi.
\item[(ii)] Supposons $S$ localement noethérien et $X$ localement de type fini (resp. de type fini) sur $S$. Alors $\pi : X \to Q$ et $Q\to S$ sont localement de présentation finie (resp. de présentation finie).
\item[(iii)] Le morphisme $R \to X\times_Q X$ de composantes $p_1$ et $p_2$ est un isomorphisme.
\end{itemize}
\end{souslem}
\begin{demo}
 Soit $U$ une quasi-section. Notons $V=p_1^{-1}(U)\subset R$ et $v : V \to X$ la restriction de $p_2$ à $V$. Notons aussi $\xymatrix@C=1pc{R_U \double & U}$ et $\xymatrix@C=1pc{R_V \double & V}$ les relations d'équivalence obtenues à partir de $R$ par image réciproque (\ref{image_reciproque}) \emph{via} les morphismes $\xymatrix@C=1pc{V\ar[r]^{p_1} &U \ar@{^(->}[r]^i &X}$ où $i$ est l'inclusion de $U$ dans $X$. Il est clair sur la définition que $R_U$ est l'intersection $p_1^{-1}(U)\cap p_2^{-1}(U)$. En particulier on a un carré cartésien :
$$\xymatrix{R_U \ar[r] \ar[d]_{p'_2={p_2}_{|_{R_U}}}& V \ar[d]^v \\ U \ar[r] &X.}$$
On en déduit que $p'_2$ est fini, localement libre et surjectif. La relation d'équivalence $R_U$ sur $U$ satisfait donc les hypothèses du théorème~\ref{cas_fini}, et elle admet un quotient catégorique. De plus, le morphisme $p_1 : V\to U$ a une section (car $p_1 : R \to X$ en a une) donc il est couvrant et l'exercice~\ref{descente_fpqc_existence_conoyau} montre que la relation d'équivalence $R_V$ sur $V$ a elle aussi un conoyau, qui de plus coïncide avec celui de $R_U$. Mais $R_V$ est obtenue par image réciproque à partir de $R$ \emph{via} $ip_1$, qui est égal au composé
$V\subset \xymatrix{R \ar[r]^{p_1} & X}$.
 L'exercice~\ref{images_reciproques_identiques} montre donc que $R_V$ s'obtient aussi par image réciproque \emph{via} le morphisme $v : V\to X$. Or ce dernier est fini et localement libre, donc couvrant, et l'exercice~\ref{descente_fpqc_existence_conoyau} donne l'existence d'un quotient catégorique $\pi: X \to Q$ pour la relation $R$ sur $X$. On a de plus un diagramme commutatif 
$$\xymatrix{X\ar[rd]_{\pi} & V
\ar[l]_v\ar[r]^{p_1}\ar[d]^{\pi_V} & U\ar[ld]^{\pi_U}\\
&Q}$$
où $\pi_U$ et $\pi_V$ sont les conoyaux des relations d'équivalence $R_U$ et $R_V$. Comme $\pi_U$ est fini et localement libre par~\ref{cas_fini}, on déduit sans difficultés de ce diagramme les assertions (i') et~(ii). Pour montrer (iii), par descente fidèlement plate il suffit de montrer que le morphisme naturel $R_V \to V\times_Q V$ est un isomorphisme, puisqu'il se déduit de $R\to X\times_Q X$ par le changement de base fidèlement plat et quasi-compact $v\times v$. Or il se déduit aussi par changement de base de $R_U\to U\times_Q U$, qui est un isomorphisme.

Il reste à démontrer que $\pi$ est un conoyau dans la catégorie de tous les espaces annelés. On déduit aisément de ce qui précède que $Q$ est bien l'\emph{ensemble} quotient de $X$ par $R$, puis l'\emph{espace topologique} quotient. En appliquant la propriété universelle du conoyau pour $\pi$ à la droite affine $\Ga$, on voit que $\Gamma(Q,\Oc_Q)$ s'identifie à l'ensemble des fonctions $a\in \Gamma(X, \Oc_X)$ telles que $p_1^*a=p_2^*a$. Ceci montre que le faisceau structural de $Q$ a les bonnes sections globales. Pour les sections au-dessus d'un ouvert quelconque $Q'$ de $Q$, on applique ce qui précède à la relation d'équivalence $R'$ sur $X\times_Q Q'$ induite par $R$ \emph{via} le changement de base $Q'\to Q$ (il faut juste vérifier que $U$ induit une quasi-section pour cette relation d'équivalence). 
\end{demo}

\begin{sousthm}[cas propre et plat]
\label{cas_propre_et_plat}
 Soient $S$ un schéma \emph{localement noethérien} et $(p_1,p_2) : \xymatrix@C=1pc{R \double &X}$ une relation d'équivalence dans $(Sch/S)$. On suppose que : 
\begin{itemize}
 \item[a)] $p_1$ est propre et plat (donc $p_2$ aussi) ;
\item[b)] $X$ est quasi-projectif sur $S$.
\end{itemize}
Alors :
\begin{itemize}
 \item[(i)] Il existe un morphisme $\pi : X \to Q$ qui est un conoyau de $(p_1,p_2)$ dans $(Sch/S)$. De plus $\pi$ est un conoyau dans la catégorie \EspAnn\ de tous les espaces annelés.
\item[(ii)] $\pi$ est propre, fidèlement plat et de présentation finie. De plus $Q\to S$ est de présentation finie.
\item[(iii)] Le morphisme $R \to X\times_Q X$ de composantes $p_1$ et $p_2$ est un isomorphisme.
\item[(iv)] On a les mêmes conséquences (iv) et (v) qu'en~\ref{cas_fini}.
\end{itemize}
\end{sousthm}
\begin{demo}
Pour prouver~(i), il suffit de montrer que tout point $z$ de $X$ possède un voisinage ouvert et saturé $U_z$ tel que le conoyau de la relation d'équivalence $R_z$ induite au-dessus de $U_z$ dans la catégorie des espaces annelés soit un conoyau dans la catégorie des schémas. Vu le lemme ci-dessus, il suffit de montrer que l'on peut choisir $U_z$ de telle sorte que $R_z$ ait une quasi-section. De plus, on peut supposer que $z$ est fermé dans sa fibre au-dessus de $S$, puisque dans chaque fibre les points fermés sont très denses. Fixons donc un tel $z$, et construisons $U_z$. 

 La question étant locale sur $S$, on peut le supposer affine. Le lemme \ref{un_bon_ferme} donne alors l'existence d'un fermé $F$ de $X$ tel que $F\cap p_1p_2^{-1}(z)$ soit fini et non vide, et tel que le morphisme composé $$\xymatrix{\tilde{p_2} : p_1^{-1}(F)\ar@{^(->}[r] & R \ar[r]^{p_2} & X}$$ soit plat aux points de $p_2^{-1}(z)$. Notons $F_R=p_1^{-1}(F)$ et $\tilde{p_1} : F_R \to F$ le morphisme induit par $p_1$.  Comme $\tilde{p_2}$ est de type fini, d'après SGA~1, IV~6.10 et EGA~3, 4.4.10, l'ensemble des points où il est plat et quasi-fini est un ouvert de $F_R$. Notons $\Phi$ le fermé complémentaire de cet ensemble dans $F_R$. C'est aussi un fermé de $R$. Comme $p_2$ est propre, $\tilde{p_2}(\Phi)$ est fermé dans $X$. On pose $U_z=\tilde{p_2}(F_R) \setminus \tilde{p_2}(\Phi)$. $U_z$ est un ouvert de $X$ car il est égal à $\tilde{p_2}(F_R\setminus \Phi) \setminus \tilde{p_2}(\Phi)$ et $\tilde{p_2}$ est plat et de type fini, donc ouvert, sur $F_R\setminus \Phi$.

Montrons que $U_z$ contient $z$. Pour cela, comme on sait déjà que $\tilde{p_2}$ est plat aux points de $p_2^{-1}(z)$, il suffit de montrer qu'il est aussi quasi-fini en ces points, c'est-à-dire que la fibre $\tilde{p_2}^{-1}(z)$ est finie. Pour cette question, on peut raisonner sur les fibres au-dessus de l'image $s$ de $z$ dans $S$, et supposer que $S$ est le spectre d'un corps. Alors $z$ est un point fermé de $X$. Comme l'ensemble $\tilde{p_1}(\tilde{p_2}^{-1}(z))$ est fini par construction de $F$, l'exercice~\ref{points_fermes_du_produit} montre que $X\times_S X$ ne contient qu'un nombre fini de points dont la première (resp. seconde) projection est dans cet ensemble (resp. est égale à $z$). Vu que $F_R$ s'injecte dans $X\times_S X$, on en déduit aussitôt que $\tilde{p_2}^{-1}(z)$ est fini.

En utilisant l'associativité de la relation d'équivalence, on montre que $U_z$ est saturé, et qu'il existe un ouvert $U$ de $F$ tel que $\tilde{p_2}^{-1}(U_z)=
\tilde{p_1}^{-1}(U)$ (ceci est admis ici, voir SGA~3 pour les détails). Notons que $U_z$ contient nécessairement $U$ puisqu'il est saturé. Il résulte alors de la construction de $U_z$ que le sous-schéma $U$ est une quasi-section pour la relation d'équivalence induite par $R$ au-dessus de $U_z$. Ceci achève la preuve de~(i).

Les assertions~(ii) et~(iii) sont locales sur $Q$. On peut donc supposer, vu ce qui précède, que la relation d'équivalence a une quasi-section. Le lemme~\ref{cas_quasi_section} donne alors toutes les conclusions sauf la séparation de $\pi$. Mais celle-ci résulte du fait que $X$ est quasi-projectif, donc séparé.
\end{demo}

\begin{souslem}
 \label{un_bon_ferme}
Soient $S$ un schéma affine et noethérien et $(p_1,p_2) : \xymatrix@C=1pc{R \double &X}$ une relation d'équivalence dans $(Sch/S)$. On suppose que $p_1$ est plat et de type fini (donc $p_2$ aussi) et que $X$ est quasi-projectif sur $S$. Soit $z$ un point de $X$ qui est fermé dans sa fibre au-dessus de $S$. Alors il existe un fermé $F$ de $X$ tel que :
\begin{itemize}
 \item[a)] $F\cap p_1p_2^{-1}(z)$ soit fini et non vide ;
\item[b)] le composé $\xymatrix@C=1pc{p_1^{-1}(F) \ar@{^(->}[r] & R \ar[r]^{p_2} & X}$ soit plat aux points de $p_2^{-1}(z)$.
\end{itemize}
\end{souslem}
\begin{demo}
 On construit une suite décroissante de fermés qui vérifient la propriété~b) et tels que $F_i\cap  p_1p_2^{-1}(z)$ soit non vide. On peut prendre $F_0=X$. Supposons $F_n$ construit. Si $F_n\cap  p_1p_2^{-1}(z)$ est fini, alors $F_n$ convient et on s'arrête. Sinon nous allons construire $F_{n+1}$. Comme $X$ est noethérien le processus doit s'arrêter donc ceci achèvera la preuve du lemme. Notons $s\in S$ l'image de $z$, et $X_s$, $R_s$ les fibres au-dessus de~$s$. On peut supposer que $X$ est un sous-schéma de $\Proj \Ac$ où $\Ac$ est l'algèbre symétrique d'un $\Gamma(S,\Oc_S)$-module de type fini. L'ensemble $p_2^{-1}(z)\cap p_1^{-1}(F_n)$ est fermé dans $R_s$. Notons $y_1, \dots, y_l$ les points génériques de ses composantes irréductibles. L'image $F_n\cap  p_1p_2^{-1}(z)$ de cet ensemble est une partie constructible de $X_s$, et infinie par hypothèse, donc elle contient une infinité de points fermés. On peut donc y choisir un point $x$ fermé dans la fibre~$X_s$ et distinct des points $p_1(y_1), \dots, p_1(y_l)$. Comme $x$ est fermé dans $X_s$, son adhérence $\ov{\{x\}}$ dans $\Proj \Ac$ ne contient pas les $p_1(y_i)$, donc il existe $f\in \Ac$ homogène de degré $d$ qui s'annule en $x$ mais pas aux points $p_1(y_i)$ (exercice analogue à~\ref{lemme_affine}). On pose alors $F_{n+1}=F_n\cap V_+(f)$. C'est un fermé de $X$, strictement inclus dans $F_n$ car il ne contient pas les $p_1(y_i)$, et non vide car il contient $x$.

Il reste à montrer que la restriction de $p_2$ à $p_1^{-1}(F_{n+1})$ est plate aux points de $p_2^{-1}(z)$. Soit $y$ un tel point. Notons $\Oc_z$ (resp. $\Oc_y$, $\ov{\Oc_y}$) l'anneau local de $z$ dans $X$ (resp. de $y$ dans $p_1^{-1}(F_n)$, de $y$ dans $p_1^{-1}(F_n)\cap p_2^{-1}(z)$ vu comme sous-schéma fermé de $R_s$). On sait par hypothèse que $\Oc_y$ est plat sur $\Oc_z$. On peut décrire l'anneau local $\Oc_y'$ de $y$ dans $p_1^{-1}(F_{n+1})$ comme suit. Soit $g\in \Ac$ homogène de degré 1 tel que $D_+(g)$ soit un voisinage de $p_1(y)$ dans $\Proj \Ac$. Au voisinage de $p_1(y)$, $F_{n+1}$ a pour équation $f/g^d$ (dans $F_n$). Donc au voisinage de $y$, $p_1^{-1}(F_{n+1})$ a pour équation (dans $p_1^{-1}(F_n)$) l'image $\varphi$ de $f/g^d$ dans $\Oc_y$. L'anneau local $\Oc_y'$ est donc $\Oc_y/\varphi \Oc_y$. Il résulte par ailleurs de la construction de $f$ que $\varphi$ ne divise pas 0 dans $\Oc_y$. On en déduit que $\Oc_y'$ est plat sur $\Oc_z$ d'après l'exercice~\ref{exercice_platitude_locale}.
\end{demo}

\begin{sousthm}[cas plat]
\label{cas_plat}
 Soient $S$ un schéma \emph{localement noethérien} et $(p_1,p_2) : \xymatrix@C=1pc{R \double &X}$ une relation d'équivalence dans $(Sch/S)$. On suppose que : 
\begin{itemize}
 \item[a)] $p_2$ est plat et de type fini ;
\item[b)] $X$ est de type fini sur $S$.
\end{itemize}
Il existe alors un ouvert $W$ de $X$ dense, saturé et satisfaisant aux propriétés suivantes :
\begin{itemize}
 \item[(i)] La relation d'équivalence $(q_1,q_2) : \xymatrix@C=1pc{R_W \double &W}$ induite par $(p_1, p_2)$ sur $W$ possède un conoyau $\pi : W \to V$ dans $(Sch/S)$. De plus $\pi$ est un conoyau dans la catégorie \EspAnn\ de tous les espaces annelés.
\item[(ii)] $\pi$ est fidèlement plat et de présentation finie, et $V\to S$ est de présentation finie.
\item[(iii)] Le morphisme $R_W \to W\times_V W$ de composantes $q_1$ et $q_2$ est un isomorphisme.
\end{itemize}
\end{sousthm}
\begin{demo}
 D'après le lemme~\ref{cas_quasi_section}, il suffit de montrer que l'on peut choisir $W$ de telle sorte que la relation d'équivalence induite sur $W$ ait une quasi-section. Il suffit en fait de montrer que, pour $z\in X$ fermé dans sa fibre au-dessus de $S$, il existe un ouvert saturé $W_z$ qui possède une quasi-section et rencontre toutes les composantes irréductibles de $X$ passant par $z$ (mais $W_z$ ne contient pas nécessairement $z$). En effet on peut alors construire un ouvert $W$ dense qui soit union disjointe de tels $W_z$.

Pour construire $W_z$, on peut supposer $S$ affine. Pour simplifier l'exposé, supposons $X$ quasi-projectif sur $S$ (voir SGA~3 pour la variante permettant de s'affranchir de cette hypothèse). On peut alors appliquer le lemme~\ref{un_bon_ferme}, qui nous donne un fermé $F$ de $X$ tel que $F\cap p_1p_2^{-1}(z)$ soit fini et non vide et que le composé $\xymatrix@C=1pc{\tilde{p_2} : p_1^{-1}(F) \ar@{^(->}[r] & R \ar[r]^{p_2} & X}$ soit plat aux points de $p_2^{-1}(z)$. De même qu'en~\ref{cas_propre_et_plat}, il résulte de~\ref{points_fermes_du_produit} que la fibre $\tilde{p_2}^{-1}(z)$ est finie. Notons $U_R$ l'ouvert de $p_1^{-1}(F)$ formé des points où $\tilde{p_2}$ est à la fois plat et quasi-fini. On note alors $W_z$ le plus grand ouvert de $\tilde{p_2}(U_R)$ au-dessus duquel $\tilde{p_2}$ est \emph{fini} et plat. Cet ouvert $W_z$ ne contient pas nécessairement $z$, mais il contient les points génériques des composantes irréductibles passant par $z$. En utilisant l'associativité de la relation d'équivalence et des arguments de descente \emph{fpqc}, on montre alors que $W_z$ est saturé, puis que $\tilde{p_2}^{-1}(W_z)$ est de la forme $\tilde{p_1}^{-1}(U)$ où $U$ est un ouvert de $F$ qui est une quasi-section pour la relation d'équivalence induite par $R$ au-dessus de $W_z$ (voir SGA~3 pour les détails). 
\end{demo}

\begin{sousexo}
 \label{points_fermes_du_produit}
Soient $k$ un corps et $X, Y$ deux $k$-schémas avec $X$ de type fini. Soient $x$ un point fermé de $X$ et $y$ un point de $Y$. Alors le produit fibré $X\times_k Y$ ne contient qu'un nombre fini de points dont les projections sont $x$ et $y$. 
\end{sousexo}

\begin{sousexo}[SGA~1, IV, 5.7]
 \label{exercice_platitude_locale}
Soit $A\to B$ un morphisme local d'anneaux locaux noethériens, $u : M'\to M$ un morphisme de $B$-modules de type fini, avec $M$ plat sur $A$, et $u\otimes_A k$ injectif (où $k$ est le corps résiduel de $A$). Alors $u$ est injectif, et $\Coker u$ est plat sur $A$. 
\end{sousexo}

\begin{sousexemple}\rm
 Soient~$S$ un schéma noethérien et~$X$ un $S$-schéma projectif, plat et à fibres géométriquement intègres. Notons $\Pic_{X/S}$ le foncteur de Picard de $X$. Grothendieck a montré dans~FGA que $\Pic_{X/S}$ est représentable par un schéma qui est union disjointe de sous-schémas ouverts, dont chacun est réunion d'ouverts quasi-projectifs. L'existence de quotients par des relations d'équivalence propres et plates est un ingrédient essentiel de la démonstration. Voici l'idée. Étant donné un polynôme~$\phi\in \Q[X]$, on note~$\Pic_{X/S}^{\phi}$ le sous-faisceau étale de $\Pic_{X/S}$ associé au préfaisceau dont les $T$-points sont les fibrés en droites $\Lc$ sur $X_T$ tels que pour tout $t\in T$ et pour tout $n$,
$$\chi(X_t, \Lc_t^{-1}(n))=\phi(n)\, .$$
On montre que $\Pic_{X/S}$ est une union disjointe des $\Pic_{X/S}^{\phi}$ ainsi définis. Notons $P_m^{\phi}$ le sous-faisceau étale de $\Pic_{X/S}^{\phi}$ associé au préfaisceau dont les $T$-points sont les fibrés en droites~$\Lc$ sur $X_T$ qui vérifient la condition précédente, et tels que de plus le faisceau~$R^if_{T*}\Lc(n)$ soit nul pour tout $i\geq 1$ et tout $n\geq m$. On montre alors que les $P_m^{\phi}$ sont des ouverts de~$\Pic_{X/S}^{\phi}$, dont $\Pic_{X/S}^{\phi}$ est la réunion croissante.  On montre enfin que $P_m^{\phi}$ est le quotient par une relation d'équivalence propre et lisse d'un ouvert convenable de $\Div_{X/S}$, et que cet ouvert est quasi-projectif (en plongeant $\Div_{X/S}$ dans le schéma de Hilbert de~$X$ sur~$S$). Le théorème~\ref{cas_propre_et_plat} montre alors que $P_m^{\phi}$ est représentable. Nous renvoyons à l'exposé de Kleiman dans~\cite{FGA_explained}, théorème~9.4.8, pour les détails de la preuve et pour la quasi-projectivité de~$P_m^{\phi}$.
\end{sousexemple}

\subsection{Espaces algébriques}
\label{espaces_algebriques}

Nous nous autorisons ici une petite digression que le lecteur désireux de rester dans le monde des schémas peut passer. La théorie des espaces algébriques est due à Michael Artin et fut développée avec Donald Knutson. Un espace algébrique est un faisceau étale qui, localement pour la topologie étale, est isomorphe à un schéma (voir la définition précise ci-dessous). De plus, la plupart des concepts et des techniques disponibles pour les schémas se généralisent aux espaces algébriques (voir~\cite{Knutson}), si bien que l'on peut travailler avec eux \og presque comme si\fg\ l'on travaillait avec des schémas. La catégorie des espaces algébriques contient celle des schémas quasi-séparés, et le difficile théorème~\ref{thm_espaces_alg_champ} montre en particulier que cette catégorie des espaces algébriques est un champ pour la topologie \emph{fppf}. Ce résultat remarquable a notamment la conséquence suivante. Soit $\Uc=\{U_i\to U\}$ une famille couvrante \emph{fppf} et soit $(X_i)$ une famille de $U_i$-schémas quasi-séparés. Alors toute donnée de descente sur cette famille est effective dans la catégorie des espaces algébriques. Autrement dit, à défaut d'obtenir un schéma en recollant les $X_i$ le long de cette donnée de descente, on obtient toujours au moins un espace algébrique. Pour la question de l'existence du quotient d'un schéma $X$ sous l'action d'un groupe $G$, on retiendra le corollaire~\ref{repres_par_esp_alg} : si l'action est libre et le groupe plat (avec des hypothèses de finitude raisonnables), alors $X/G$ est au moins un espace algébrique. 

\begin{sousdefi}
 Un espace algébrique (sous-entendu, quasi-séparé) est un faisceau étale
$$X : (Sch/S)^o \flechelongue \Ens$$ vérifiant les propriétés :
\begin{itemize}
 \item[(i)] Le morphisme diagonal $\Delta : X \fleche X\times_S X$ est schématique et quasi-compact.
\item[(ii)] Il existe un schéma $X'$ et un morphisme de faisceaux $\pi : X'\to X$ (automatiquement schématique par~(i)) qui est étale et surjectif.
\end{itemize}
\end{sousdefi}

\begin{sousexo}\rm
 Soient $X$ un $S$-schéma et $Y, Z$ deux $X$-schémas (ou plus généralement $X$, $Y$, $Z$ des foncteurs de $(Sch/S)^o$ vers $\Ens$ avec des morphismes $Y\fleche X$ et $Z\fleche X$). Vérifier que le carré
$$\xymatrix{Y\times_X Z \ar[r] \ar[d] & Y\times_S Z\ar[d] \\
X \ar[r]^-{\Delta} & X\times_S X}$$
est cartesien (où $\Delta$ est le morphisme diagonal). En déduire qu'un morphisme $\pi$ comme en~(ii) est automatiquement schématique sous la condition~(i).
\end{sousexo}

De manière équivalente, on peut définir un espace algébrique comme étant le quotient d'une relation d'équivalence étale.

\begin{sousprop}\ 

\begin{itemize}
\item[a)] Si $X$ est un espace algébrique et si $\pi : X' \to X$ en est une présentation (\emph{i.e.}  $X'$ est un schéma et $\pi$ est étale et surjectif) alors $X$ est le faisceau quotient de la relation d'équivalence
$$\xymatrix{R=X'\times_X X' \ar@<2pt>[r]^-{p_1} \ar@<-2pt>[r]_-{p_2}& X'}.$$ On vérifie facilement qu'ici les projections $p_1, p_2$ sont étales et le monomorphisme $R\to X'\times_SX'$ est quasi-compact.
\item[b)] Soit $$\xymatrix{R \ar@<2pt>[r]^{p_1} \ar@<-2pt>[r]_{p_2}& X'}$$ une relation d'équivalence dans $(Sch/S)$ avec $p_1, p_2$ étales et $\delta : R\to X'\times_S X'$ quasi-compact. Alors le faisceau quotient $X'/R$ est un espace algébrique et le morphisme canonique $\pi : X' \to X'/R$ en est une présentation.
\end{itemize}
\end{sousprop}
\begin{demo}
Voir \cite{Knutson}~II~1.3. On remarquera que a) est facile. En revanche, la réciproque b) nécessite des arguments de descente et un peu de travail. 
\end{demo}

\begin{sousexo}\rm Soit $G \to S$ un schéma en groupes étale et de présentation finie qui agit librement sur un $S$-schéma $X$ quasi-séparé. Montrer que la relation d'équivalence associée $R=G\times_S X \to X\times_S X$ vérifie les conditions de la proposition précédente et en déduire que le faisceau quotient $G\backslash X$ est un espace algébrique. 
\end{sousexo}

\begin{sousthm}[Artin]
\label{thm_espaces_alg_champ}
 Soit $\xymatrix@C=1pc{R \ar@<2pt>[r]^{p_1} \ar@<-2pt>[r]_{p_2}& X}$ une relation d'équivalence dans $(Sch/S)$ (ou plus généralement dans la catégorie des espaces algébriques sur $S$) avec $p_1$ et $p_2$ plats et localement de présentation finie et $\delta : R \fleche X\times_S X$ quasi-compact. Alors le faisceau quotient $X/R$ est un espace algébrique.
\end{sousthm}
\begin{demo}
 \cite{LMB}~10.4
\end{demo}

\begin{souscor}
\label{repres_par_esp_alg}
 Soit $G\fleche S$ un espace algébrique en groupes plat et de présentation finie qui agit librement sur un espace algébrique $X$ quasi-séparé. Alors le faisceau quotient $X/G$ est un espace algébrique.
\end{souscor}
\begin{demo} C'est juste un cas particulier du théorème précédent.
\end{demo}

\begin{sousexo}[\cite{Stacks_Project}, Example 03FN]\rm
\label{exemple_espace_alg}
 Soient $S=\Spec \R[x]$ et $U=\Spec \C[x]$.\\
1) Montrer que $U\times_S U$ est isomorphe à $U_1\coprod U_2$ avec $U_i=U$, les projections sur $U$ s'identifiant aux applications $\id\coprod \id$ et $\id\coprod \sigma$ où $\sigma$ est la conjugaison complexe.\\
2) On note $R=U_1\coprod (U_2\setminus\{0_{U_2}\})$. Montrer que l'inclusion de $R$ dans $U\times_S U$ est une relation d'équivalence sur $U$ et que le quotient $X=U/R$ est un espace algébrique.\\
3) En utilisant~\ref{prop_faisceau_quotient}, montrer que le morphisme $f : X \to S$ induit un isomorphisme $f^{-1}(S\setminus\{0_S\}) \to S\setminus\{0_S\}$, que la fibre $f^{-1}(0_S)$ s'identifie à $0_U=\Spec \C$, et que l'espace algébrique obtenu à partir de $X$ par le changement de base $U \to S$ (ou, ce qui revient au même, $\Spec \C \to \Spec \R$) est la droite affine complexe avec origine dédoublée.\\
4) Montrer que $X$ n'est pas un schéma (supposer le contraire, et regarder l'anneau local de l'unique point $0_X$ au-dessus de $0_S$).\\
5) Montrer que $X\to S$ n'est pas séparé, mais qu'il est tout de même localement séparé et quasi-séparé (\emph{i.e.} sa diagonale est une immersion quasi-compacte). Il est de plus étale.
\end{sousexo}

Le célèbre exemple de Hironaka (voir~\cite{Hironaka_nonprojective_example}), dont nous avons déjà parlé en~\ref{exemple_Hironaka}, d'une variété propre et lisse de dimension~3 sur~$\C$ qui n'est pas projective, a été utilisé à maintes reprises pour construire des exemples de données de descente non effectives. On obtient ainsi un espace algébrique lisse de dimension~3 sur $\C$ qui n'est pas un schéma (\emph{cf.} par exemple~\cite{Knutson}, \cite[4.4.2]{FGA_explained} ou~\cite{Mumford_GIT}). Historiquement, c'est Nagata qui a donné dans~\cite{Nagata_Existence_theorems} le premier exemple de surface algébrique $X$ propre et normale mais non projective (sur un corps assez gros). La surface est munie d'une involution $\tau$ qui échange deux points $x$ et~$x'$ qui ne sont simultanément contenus dans aucun ouvert affine. Ceci empêche le faisceau quotient $X/\tau$ d'être un schéma. C'est cependant un espace algébrique.

\subsection{Quotients de groupes sur un corps (voire un anneau artinien)}
\label{cas_ou_la_base_est_un_corps}

Dans toute cette section, $k$ désignera un corps. Il convient de mentionner que la plupart des résultats ci-dessous sont encore valables si l'on remplace $k$ par un anneau artinien, modulo quelques hypothèses de platitude sur les objets en jeu. Pour simplifier un peu, nous avons préféré nous contenter de résultats sur un corps. Compte tenu de~\ref{representabilite_et_nilpotents}, cette restriction n'est pas très gênante, et nous renvoyons à SGA 3 \cite{SGA3} pour les variantes sur un anneau artinien. 

\begin{sousthm}[\cite{SGA3} SGA 3 VI$_A$ 3.2]
 Soient $F$ et $G$ des $k$-schémas en groupes localement de type fini et soit $u : F \to G$ un monomorphisme quasi-compact de $k$-groupes. Alors le faisceau quotient $G/F$ (au sens \emph{fppf}) est représentable par un $k$-schéma.
\end{sousthm}

\begin{sousremarque}\rm
 On rappelle que par~\ref{comparaison_conoyau_faisceau} et les exercices qui suivent, le quotient $G/F$ dont ce théorème donne l'existence jouit d'un certain nombre de propriétés agréables. Ainsi le morphisme $\pi : G\to G/F$ est fidèlement plat et localement de présentation finie, le morphisme canonique $F\times_k G \to G\times_{(G/F)}G$ est un isomorphisme, $G/F$ est séparé et $u : F\to G$ est une immersion fermée.\footnote{Les deux dernières assertions demandent encore un peu de travail...} Si $G$ est lisse, alors $G/F$ l'est aussi. Enfin, si $F$ est un sous-groupe invariant de $G$, il existe sur $G/F$ une unique structure de $k$-groupe telle que $\pi$ soit un morphisme de groupes. 
\end{sousremarque}
\begin{demo}
 Nous donnons seulement la preuve dans le cas où $F$ et $G$ sont \emph{de type fini} sur $k$. L'idée est de montrer que le théorème est valable après une extension finie de~$k$, puis de conclure par des arguments de descente.

\begin{etape}{Assertion 1 : Si $G/F$ est représentable, alors toute partie finie de $G/F$ est contenue dans un ouvert affine.}
Notons $X=G/F$ et $\pi$ la projection de $G$ sur $X$. Soit $V$ un ouvert affine dense dans $X$ (il en existe car $X$ est de type fini sur $k$). Soient $x_1, \dots, x_n$ des points de $X$. Supposons dans un premier temps qu'il existe pour tout $i$ un point $k$-rationnel $g_i\in G$ au-dessus de $x_i$, et que l'ouvert $\cap_{i=1}^n g_i.(\pi^{-1}(V))^{-1}$ de $G$ (automatiquement dense car $\pi$ est ouvert) contienne un point $k$-rationnel $g$. Alors pour tout $i$, $g\in g_i.(\pi^{-1}(V))^{-1}$ donc $g_i\in g.\pi^{-1}(V)$ et $x_i\in g.V$. Donc l'ouvert affine $g.V$, image de $V$ par la translation à gauche par $g$ dans~$X$, convient.

Dans le cas général, on peut supposer les $x_i$ fermés dans $X$. Pour tout $i$ soit $g_i$ un point fermé au-dessus de $x_i$. Soit $K$ une extension finie de $k$ telle que tous les $g'_j\in G_K$ au-dessus des $g_i$ soient $K$-rationnels (prendre par exemple pour $K$ une extension normale de $k$ qui contient les corps résiduels des $g_i$). Alors $\cap_{j=1}^p g'_j.(\pi_K^{-1}(V_K))^{-1}$ est un ouvert dense de $G_K$, donc contient un point fermé $g$. Quitte à agrandir $K$ on peut supposer que $g$ est $K$-rationnel. Le cas traité précédemment montre alors qu'il existe un ouvert affine $V'$ de $X_K$ contenant les images $x'_j$ des $g'_j$. Comme les $x'_1, \dots, x'_p$ sont \emph{tous} les points de $X_K$ au-dessus des $x_i$, il forment une réunion d'orbites pour la relation d'équivalence finie et localement libre sur $X_K$ définie par la projection $X_K \to X$. En raisonnant comme dans~\ref{cas_fini} (construction d'un voisinage ouvert affine et saturé), on trouve un ouvert affine et saturé $W'\subset V'$ qui contient tous les $x'_j$. Son image~$W$ dans $X$ contient tous les $x_i$, et c'est un ouvert affine car quotient de l'affine $W'$ par une relation d'équivalence finie et localement libre (voir le cas affine de~\ref{cas_fini}).
\end{etape}

\begin{etape}{Assertion 2 : Il existe une extension finie $L$ de $k$ telle que $(G/F)_L$ soit représentable.}
 Ici $(G/F)_L$ désigne le foncteur sur $\Spec L$ obtenu à partir de $G/F$ par le changement de base $\Spec L\to \Spec k$. Si $k'$ est une extension finie de $k$, on note $U[k']$ la réunion des ouverts $W \subset G_{k'}$, stables sous l'action à droite de $F_{k'}$, et tels que le faisceau quotient~$W/F_{k'}$ soit un schéma. Si $W$ est un tel ouvert, ses translatés à gauche par les points de $G(k')$ vérifient encore les mêmes propriétés, donc $U[k']$ est stable par multiplication à gauche par les $k'$-points de $G_{k'}$. D'après~\ref{cas_plat}, on voit que $U[k]$ est dense dans $G$ et contient en particulier un point fermé. Quitte à remplacer $k$ par une extension finie, on peut supposer que $U[k]$ contient un $k$-point. Alors pour toute extension $k'/k$, $U[k']$ contient tous les $k'$-points, et l'exercice~\ref{fin_preuve_assertion2} donne le résultat.
\end{etape}

\begin{etapefinale}{Conclusion (descente)}
Le morphisme $\varphi : \Spec L \to \Spec k$ est en particulier fidèlement plat et de présentation finie. Alors d'après~\ref{descente_representabilite}, le schéma $(G/F)_L$ est muni d'une donnée de descente relativement à $\varphi$. D'après l'assertion 1, toute partie finie de $(G/F)_L$ est contenue dans un ouvert affine. Par SGA 1 \cite{SGA1}~VIII~7.6, la donnée de descente sur $(G/F)_L$ est alors effective (ici on pourrait aussi utiliser~\ref{cas_fini}), ce qui par~\ref{descente_representabilite} entraîne la représentabilité de $G/F$.
\end{etapefinale}
\end{demo}

\begin{sousexo}\rm
\label{fin_preuve_assertion2}
 Soit $X$ un schéma de type fini sur un corps $k$. Pour toute extension $k'$ de $k$, on suppose donné un ouvert $U[k']$ de $X_{k'}$, avec les propriétés suivantes :
\begin{itemize}
 \item[(i)] Si $k\subset k'\subset k''$, alors $U[k'']$ contient $U[k']\otimes_{k'}k''$.
\item[(ii)] $U[k']$ contient tous les points $k'$-rationnels de $X_{k'}$.
\end{itemize}
Montrer par récurrence sur la dimension du fermé $X\setminus U[k]$ qu'il existe une extension finie $k'$ de $k$ telle que $U[k']=X_{k'}$. [Indication : Notons $Z[k']=X_{k'}\setminus U[k']$. On choisit un point fermé $x_i$ dans chaque composante irréductible de $Z[k]$. Soit $K/k$ une extension finie et normale de $k$ qui contient les corps résiduels de tous les $x_i$. Alors tout point $x'_j$ de $X_K$ au-dessus d'un $x_i$ est $K$-rationnel. En déduire que $\dim Z[K] < \dim Z[k]$ et conclure grâce à l'hypothèse de récurrence.]
\end{sousexo}

Nous signalons ci-dessous quelques conséquences de ce théorème. Nous renvoyons à SGA 3~\cite{SGA3} pour les preuves, pour d'autres énoncés, ou le cas échéant pour des énoncés analogues valables sur un anneau artinien.

\begin{sousexo}[\cite{SGA3} IV 5.2.7 et VI$_A$ 5.3.1]\rm
 Soient $k$ un corps, $G$ un $k$-groupe localement de type fini et $F$ un sous-groupe invariant fermé\footnote{Cette dernière hypothèse est en réalité superflue. En effet, si $G$ est un schéma en groupes sur un corps, ou plus généralement sur un anneau artinien, alors tout sous-groupe de $G$ est fermé, \emph{cf.} \cite{SGA3}~VI$_A$~0.5.2.} de $G$. Alors les applications $H\mapsto H/F$ et $H'\mapsto H'\times_{(G/F)} G$ définissent une correspondance bijective entre les $k$-sous-groupes de $G/F$ et les $k$-sous-groupes de $G$ contenant $F$.
\end{sousexo}

\begin{sousexo}[\cite{SGA3} IV 5.2.9 et VI$_A$ 5.3.2]\rm
 Soient $k$ un corps, $G$ un $k$-groupe localement de type fini et $F\subset H\subset G$ des sous-groupes (fermés) de $G$ avec $F$ invariant dans $H$. Alors :
\begin{itemize}
 \item[(i)] $H/F$ opère librement à droite sur $G/F$ ;
\item[(ii)] le quotient $(G/F)/(H/F)$ existe ;
\item[(iii)] on a un isomorphisme canonique de $k$-schémas (munis d'actions de $G$) :
$$(G/F)/(H/F)\simeq G/H.$$
\end{itemize}
\end{sousexo}

\begin{sousprop}[\cite{SGA3} VI$_A$ 5.4.1]
 Soient $G,H$ des $k$-groupes localement de type fini et $u : G\to H$ un morphisme quasi-compact. On note $N=\Ker u$. Alors on a la factorisation :
$$\xymatrix{G \ar[r]^u \ar[d]_{\pi} & H \\ G/N \ar[ru]_i}$$
où $\pi$ est fidèlement plat, localement de présentation finie, et $i$ est une immersion fermée.
\end{sousprop}

\begin{sousremarque}\rm
 Cet énoncé se généralise aux $k$-groupes quelconques (non nécessairement localement de type fini) de la manière suivante. Soit $u : G\to H$ un morphisme quasi-compact de $k$-groupes. On note $N=\Ker u$. Alors le faisceau quotient \emph{fpqc} $G/N$ est représentable par un $k$-groupe, et $u$ se factorise comme ci-dessus avec $\pi$ fidèlement plat et $i$ une immersion fermée. (\emph{cf.} \cite{Perrin_groupes_qc}~V~3.3 à~3.4) En particulier, si $u$ est schématiquement dominant, il est fidèlement plat.
\end{sousremarque}

\begin{sousthm}[\cite{SGA3} VI$_A$~5.4.2, VI$_A$ 5.4.3 et \cite{Perrin_approximation_groupes_qc}~4.2.6]\ 
\label{categories_abeliennes}

\noindent Les catégories suivantes sont abéliennes :
\begin{itemize}
\item la catégorie des $k$-groupes commutatifs et quasi-compacts ;
 \item la catégorie des $k$-groupes commutatifs et de type fini ;
\item la catégorie des $k$-groupes commutatifs, de type fini et affines.
\end{itemize}
\end{sousthm}

\begin{sousremarque}\rm
 Pour conclure, mentionnons rapidement deux cas particuliers de quotients. Soit $G$ un $k$-groupe localement de type fini. Il est montré en~\cite{SGA3}~VI$_A$~5.5.1 que~$G/G^0$ est étale sur $k$, et même constant lorsque $k$ est algébriquement clos. Par ailleurs, si~$k$ est parfait, alors $G/G_{\text{réd}}$ est le spectre d'une $k$-algèbre finie et locale, de corps résiduel~$k$ (\cite{SGA3}~VI$_A$~5.6.1).
\end{sousremarque}

\subsection{Quotient par le normalisateur d'un sous-groupe lisse}
\label{section_normalisateur}

Le théorème ci-dessous donne un exemple non trivial d'existence de schéma quotient. La preuve donnée repose en partie sur l'existence d'un \emph{espace algébrique} quotient (donc sur les théorèmes d'Artin), qui remplace les techniques de Murre utilisées dans~SGA 3.

\begin{sousdefi}
 Soient $S$ un schéma et $H\to G$ un monomorphisme de $S$-groupes. On définit le centralisateur de $H$ dans $G$, noté $C=\Centr_G(H)$, et le normalisateur de $H$ dans $G$, noté $N=\Norm_G(H)$, comme étant les sous-foncteurs de $G$ sur $(Sch/S)^o$ définis fonctoriellement par
\beqn
N(U)&=&\{g\in G(U)\ |\ \forall V\fleche U\, , \ gH(V)g^{-1}=H(V) \} \\
C(U)&=&\{g\in G(U)\ |\ \forall V\fleche U\, ,\ \forall h\in H(V)\, , \ ghg^{-1}=h \}
\eeqn
\end{sousdefi}

\begin{sousthm}[SGA~3, tome II, exp. XVI, cor. 2.4]\ 

\label{quotient_par_norm}
\noindent Soient $S$ un schéma et $G$ un $S$-schéma en groupes de présentation finie sur $S$. Soit $H$ un sous-groupe de $G$, c'est-à-dire un $S$-schéma en groupes muni d'un morphisme de $S$-groupes $H\fleche G$ qui est un monomorphisme. On suppose $H$ lisse sur $S$ et à fibres connexes (il est alors automatiquement de présentation finie sur $S$ par~SGA~3~VIB~5.3.3). On note $N=\Norm_G(H)$ le normalisateur de $H$ dans $G$.
\begin{itemize}                                                                                                                                                                                                                                                                                                                   \item[(i)] Alors $N$ est représentable par un sous-schéma en groupes fermé de $G$, de présentation finie sur $S$. (Voir aussi SGA~3~XI~6.11.) 
\item[(ii)] On suppose de plus $N$ plat sur $S$. Alors $G/N$ est représentable par un $S$-schéma, de présentation finie sur $S$, et quasi-projectif sur $S$.                                                                                                                                                                                                                                                                                                            \end{itemize}
\end{sousthm}

 Nous allons seulement donner la preuve de ce théorème dans le cas particulier où $S$ et $G$ sont affines et où $H$ est un sous-groupe fermé de~$G$. Par des arguments standard de passage à la limite, on peut supposer $S$ noethérien. On note $S=\Spec A$, $G=\Spec B$ et $H=\Spec B/I$. On peut aussi supposer $S$ connexe. Alors $H$ est de dimension relative constante sur $S$. Notons $G^{(n)}$ (resp. $H^{(n)}$) le~$n$\up{ième} voisinage infinitésimal de la section unité dans $G$ (resp. dans $H$). Autrement dit, si~$J$ désigne l'idéal de $B$ correspondant à la section unité, on a $G^{(n)}=\Spec B/J^{n+1}$ et~$H^{(n)}=\Spec B/(J^{n+1}+I)$. Comme $H$ est lisse, le module $B/(J^{n+1}+I)$ est libre sur $A$ de rang fini (déterminé par $n$ et par la dimension relative de $H$ sur $S$). Notons encore\footnote{On rappelle que si $\Ec$ est un module quasi-cohérent sur un schéma $S$ et si $n$ est un entier, la grassmannienne $\Grass_n(\Ec)$ est définie fonctoriellement de la manière suivante. Pour tout $S$-schéma $T$, $\Grass_n(\Ec)(T)$ est l'ensemble des $\Oc_T$-modules quotients de $\Ec_T$ qui sont localement libres de rang $n$. Le foncteur $\Grass_n(\Ec)$ est toujours représentable par un schéma séparé sur $S$. Si de plus $\Ec$ est cohérent, ce schéma est projectif sur $S$ (\emph{cf.}~EGA~1~\cite{EGA1}~9.7 et~9.8).} $X_n=\Grass_l(B/J^{n+1})$ où $l$ est le rang sur $A$ de $B/(J^{n+1}+I)$.
L'action de $G$ sur lui-même par automorphismes intérieurs induit une action naturelle de $G$ sur $X_n$.
Par ailleurs le quotient $B/(J^{n+1}+I)$ de $B/J^{n+1}$ définit un $S$-point (qui est une immersion fermée) de $X_n$
$$\xi_n : S \flechelongue X_n.$$
On note $N^{(n)}=\Norm_G(H^{(n)})$ le normalisateur dans $G$ de $H^{(n)}$. On vérifie alors que l'on a un carré cartésien
$$\xymatrix{N^{(n)} \cartesien\ar[r] \ar[d] & G\ar[d]^{\omega(\xi_n)} \\
S \ar[r]_{\xi_n} & X_n
}$$
où le morphisme $\omega(\xi_n) : G \fleche X_n$ est l'orbite de $\xi_n$. 
Ceci prouve en particulier que~$N^{(n)}$ est représentable par un sous-schéma fermé de $G$. Les $N^{(n)}$ ainsi construits forment une suite décroissante de sous-foncteurs de $G$ :
$$N^{(0)} \supset N^{(1)} \supset N^{(2)} \supset \dots \supset N\, .$$
Soit $R$ l'intersection des sous-foncteurs $N^{(n)}$. On a évidemment l'inclusion $N\subset R$. Il se trouve que l'inclusion réciproque est vraie aussi\footnote{En effet, soit $g\in R(S)$ et montrons que $g\in N(S)$. Comme le raisonnement qui suit vaut après tout changement de base ceci prouvera que $R\subset N$. Notons $Y$ le sous-schéma en groupes fermé de~$H$ intersection de $H$ et de $gHg^{-1}$. Il faut montrer que $Y=H$. Pour tout $n$, on sait par hypothèse que~$g\in N^{(n)}$ donc $Y$ contient les voisinages infinitésimaux $H^{(n)}$ de la section unité. Ceci implique que le morphisme~$Y\to H$ est étale en tout point de la section unité de $Y$ sur $S$. Notons $Y_{\textrm{ét}}$ l'ouvert de~$Y$ formé des points où~$Y\to H$ est étale.
C'est un sous-groupe ouvert de $Y$ par~\cite{SGA3}~SGA~3~VI$_B$~2.2. Le morphisme~$Y_{\textrm{ét}}\to H$ est un monomorphisme étale, donc une immersion ouverte. Maintenant, pour tout~$s\in S$, la fibre $(Y_{\textrm{ét}})_s$ est un sous-groupe ouvert de $H_s$, mais il est aussi fermé (comme tout sous-groupe d'un schéma en groupes sur un corps, \cite{SGA3}~SGA~3~VI$_A$~0.5.2).
Comme $H_s$ est par hypothèse connexe on en déduit que l'immersion ouverte $Y_{\textrm{ét}}\to H$ est surjective. A fortiori $Y=H$ et le résultat est démontré.
}.
De plus, la suite $N^{(n)}$ stationne puisque $G$ est noethérien. Il en résulte que $N$ est égal à $N^{(n)}$ pour $n$ assez grand. En particulier $N$ est un sous-schéma fermé de $G$, ce qui prouve le point~(i).

Supposons maintenant $N$ plat et montrons que le quotient $G/N$ est représentable. Commençons par remarquer que d'après les résultats d'Artin (\ref{repres_par_esp_alg}) $G/N$ est un espace algébrique de présentation finie sur $S$. Il reste donc seulement à montrer qu'il est quasi-projectif. Fixons un entier $n$ tel que $N=N^{(n)}$. Le morphisme d'orbite $\omega(\xi_n)$ induit alors un monomorphisme
$$\varphi : G/N=G/N^{(n)} \flechelongue X_n\, .$$
Comme $G/N$ et $X_n$ sont de présentation finie sur $S$, $\varphi$ est lui-même de présentation finie. D'après le théorème principal de Zariski (généralisé aux espaces algébriques, \emph{cf.} par exemple~\cite{LMB}~A.2.2), le morphisme $\varphi$ est quasi-affine, ce qui achève la preuve. $\square$
\vskip 3mm

\begin{sousremarque}[représentabilité des centralisateurs]\rm
 Par le même genre d'arguments, on peut montrer que, sous les mêmes hypothèses, le centralisateur $\Centr_G(H)$ est lui aussi représentable par un sous-schéma en groupes fermé de $G$, de présentation finie sur $S$ (voir \cite{SGA3}~XI~6.11).
\end{sousremarque}

 Le théorème~\ref{quotient_par_norm} s'appliquera en particulier au cas du quotient d'un schéma en groupes réductif par un sous-groupe parabolique ou un sous-groupe de Lévi\footnote{Rappelons quelques définitions. Soit $G$ un schéma en groupes affine, lisse et de présentation finie sur un schéma $S$. Si $S$ est le spectre d'un corps algébriquement clos, un sous-groupe de Borel de $G$ est un sous-groupe algébrique lisse résoluble connexe, et maximal pour ces propriétés. Un sous-groupe parabolique est un sous-groupe qui contient un sous-groupe de Borel. Maintenant si $S$ est une base quelconque, un sous-groupe de Borel est un sous-schéma en groupes lisse et de présentation finie dont les fibres géométriques sont des sous-groupes de Borel des $G_{\ov{s}}$. De même un sous-groupe parabolique est un sous-schéma en groupes lisse et de présentation finie dont les fibres géométriques sont des sous-groupes paraboliques.}  d'icelui. En effet, nous verrons dans les cours de Conrad, Gross et Yu le résultat suivant.

\begin{sousprop}
 Soit $G$ un $S$-schéma en groupes réductif. Soit $P$ un sous-groupe parabolique de $G$. Alors $P$ est identique à son propre normalisateur.
\end{sousprop}

On en déduit\footnote{La preuve de~(ii) ci-dessous donnée dans SGA~3~XXVI~3 repose sur l'étude du schéma des sous-groupes critiques de~$G$. La preuve que nous esquissons est une alternative et s'appuie sur le théorème ci-dessus mais aussi sur des travaux de Seshadri pour obtenir le complément affine.} :

\begin{souscor}[SGA 3~\cite{SGA3}~XXVI~1.2 et~3.13] Soient $G$ un $S$-schéma en groupes réductif et $P$ un sous-groupe parabolique de $G$.
 \begin{itemize}
  \item[(i)] Alors $P$ est un sous-schéma fermé de $G$, et le quotient $G/P$ est représentable par un schéma projectif et lisse sur~$S$.
\item[(ii)] Si $L$ est un sous-groupe de Lévi de $P$, les quotients $G/L$ et $G/\Norm_G(L)$ sont représentables par des  schémas affines et lisses sur~$S$.
 \end{itemize}
\end{souscor}
\begin{demo}
 (i) Vu que $P$ est lisse sur $S$ et à fibres connexes, on peut appliquer le théorème~\ref{quotient_par_norm} et on voit que $G/\Norm_G(P)=G/P$ est un $S$-schéma quasi-projectif et de présentation finie sur $S$. Il est lisse par~\ref{exercice_pptes_1}. Enfin le morphisme $G/P\to S$ a une section, et ses fibres sont propres et géométriquement connexes, donc il est propre par~\cite{EGA}~EGA~IV$_3$, 15.7.11.

(ii) Par~\cite{SGA3}~XXII~5.10.2, le normalisateur $\Norm_G(L)$ est lisse sur $S$, et le groupe quotient $\Norm_G(L)/L$ est fini étale sur $S$. Il résulte de~\ref{quotient_par_norm} et~\ref{exercice_pptes_1} que $G/\Norm_G(L)$ est un $S$-schéma quasi-projectif lisse et de présentation finie sur $S$. De plus $\Norm_G(L)/L$ agit à droite sur $G/L$, et comme on vient de voir que le quotient est représentable, on en déduit facilement que la projection $G/L \to G/\Norm_G(L)$ est étale et finie (car $\Norm_G(L)/L$ l'est). Il reste à montrer que ces quotients sont affines. Par~\cite{EGA}~EGA~II~6.7.1 il suffit de le faire pour $G/L$. On peut ici appliquer directement un raisonnement de Colliot-Thélène et Sansuc (voir~\cite[6.12]{Colliot-Thelene_Sansuc_Fibres_quadratiques}) reposant sur un résultat de Seshadri. Voici l'idée. La question est locale sur $S$ pour la topologie étale donc on peut supposer $S$ affine et $G$ déployé (\cite[XXII, 2.3]{SGA3}). Alors $G$ provient d'un $\Z$-groupe réductif déployé (\cite[XXV, 1.1]{SGA3}). Utilisant~\cite[I.1, proposition 3 p. 236 et lemme 1 p. 230]{Seshadri_Geometric_reductivity}, on en déduit l'existence d'une immersion fermée $i : G \to \GL_{n,S}$ qui est un morphisme de groupes. L'action de $G$ sur lui-même par translations à droite est alors linéarisable (car il en est ainsi pour $\GL_{n,S}$). En particulier l'action de $L$ sur $G$ est linéarisable. On pose $W=\Spec((p_*\Oc_G)^L)$ où $p : G\to S$ est le morphisme structural de $G$. D'après~\cite[théorème~3 p.268 et remarque~8 p.269]{Seshadri_Geometric_reductivity}, $W$ est un quotient catégorique de $G$ sous l'action de $L$. Comme le faisceau quotient $G/L$ en est aussi un puisque c'est un schéma, ils sont isomorphes et $G/L$ est affine.
\end{demo}

\subsection{Quotient d'un schéma affine par un groupe diagonalisable opérant librement}
\label{cas_groupe_diago}

Soient $S$ un schéma et $P$ un schéma affine sur $S$. On note $P=\Spec \Ac$ avec $\Ac$ une $\Oc_S$-algèbre quasi-cohérente. Soit $G=D_S(M)$ un groupe diagonalisable. Alors se donner une action de $G$ sur $P$ revient à se donner une $M$-graduation de $\Ac$, \emph{i.e.} une décomposition de $\Ac$ en
$$\Ac = \bigoplus_{i\in M} \Ac_i$$
avec $\Ac_i.\Ac_{j}\subset \Ac_{i+j}$ pour tous $i, j\in M$.
En notant $\Oc_S[M]$ l'algèbre de $G$ et $X^i$, $i\in M$, ses générateurs canoniques comme $\O_S$-module, l'action associée à une telle décomposition de $\Ac$ est donnée par le morphisme :
$$\fonction{\varphi}{\Ac}{\Ac[M]}{a}{\sum_i a_iX^i}.$$

L'énoncé ci-dessous donne des conditions nécessaires et suffisantes pour qu'un tel $P$ muni d'une action de $G$ soit un $G$-torseur sur $S$.

\begin{sousprop}[\cite{SGA3}~SGA~3~VIII~4.6]
\label{torseurs_sous_gpe_diagonalisable}
 Avec les notations ci-dessus, $P$ est un $G$-torseur (au sens \emph{fpqc}) si et seulement si les conditions suivantes sont vérifiées :
\begin{itemize}
 \item[a)] Le morphisme $\Oc_S \to \Ac_0$ est un isomorphisme.
\item[b)] Pour tout $i\in M$, on a $\Ac_i.\Ac_{-i}=\Ac_0$.  
\end{itemize}
\end{sousprop}

\begin{sousthm}[SGA 3~\cite{SGA3}~VIII~5.1]
\label{quotient_par_gpe_diagonalisable}
 Soient $S$ un schéma, $M$ un groupe commutatif ordinaire, $G=D_S(M)$ le groupe diagonalisable associé, et $P$ un schéma affine sur $S$ sur lequel $G$ opère librement à droite.

Alors le quotient $X=P/G$ existe et $P$ est un $G_X$-torseur sur $X$ (où $G_X=G\times_S X$). De plus, $P/G$ est affine sur $S$. Plus précisément, si $P=\Spec \Ac$ où $\Ac$ est une $\Oc_S$-algèbre quasi-cohérente, alors $P/G$ est isomorphe à $\Spec \Ac^G$ où $\Ac^G$ est l'algèbre des invariants sous $G$. (Remarque : $\Ac^G$ est aussi le composant $\Ac_0$ de degré 0 de $\Ac$ avec les notations ci-dessus.)

Enfin, si $P$ est de type fini (resp. de présentation finie) sur $S$, il en est de même de~$X$.
\end{sousthm}
\begin{demo}
On note $X=\Spec \Ac_0$. Le morphisme $P\to X$ est alors $G$-invariant, donc l'action de $G$ sur $P$ (par $S$-morphismes) induit une action de $G_X$ sur $P$ (par $X$-morphismes). Comme l'action de départ est libre, on voit facilement que l'action de $G_X$ sur $P$ (au-dessus de $X$) est libre. D'après~\ref{comparaison_conoyau_faisceau}, il suffit pour conclure de montrer que cette action fait de $P$ un $G_X$-torseur sur $X$. On peut ainsi supposer $X=S$. De plus, être un torseur étant une propriété locale sur $S$, on peut supposer $S$ affine. On utilise la proposition~\ref{torseurs_sous_gpe_diagonalisable} ci-dessus. La condition a) est évidente. Donc pour conclure la première partie il suffit de montrer que pour tout $i\in M$ on a $A_i.A_{-i}=A_0$. C'est le point-clé de la preuve, et l'objet de l'exercice ci-dessous.

Il reste à montrer la dernière assertion. Supposons donc $P$ de type fini sur $S$. Alors $P$ est de type fini sur $X$. En particulier, après un changement de base qui trivialise le torseur $P$ sur $X$, on voit par \ref{prop_gpes_diagonalisables} que $M$ est de type fini\footnote{On peut suppose $X$ non vide, ce qui est évidemment loisible.}. Le groupe $G_X$ est donc de présentation finie sur $X$. Comme $P$ est un $G_X$-torseur sur $X$, on en déduit par descente que $P$ lui-même est de présentation finie sur $X$. On conclut par EGA~IV~\cite{EGA}~(11.3.16).
\end{demo}

\begin{sousexo} \rm
\label{exercice_cas_diag}
On garde les notations de la preuve. On montre ici que pour tout $i\in M$ on a $A_i.A_{-i}=A_0$. Le cas $M=0$ étant trivial, on peut supposer $M\neq 0$.

\noindent {\bf 1.} Montrer qu'il existe $i\in M\setminus\{0\}$ tel que $A_i.A_{-i}\neq 0$. [\emph{Indication :} Raisonner par l'absurde, et montrer que dans le cas contraire, la projection canonique $\pi : A \to A_0$ est un morphisme d'anneaux, ce qui définit un point $\pi\in P(S)$ fixe sous $G$ et contredit la liberté de l'action.]

\noindent {\bf 2.}  Montrer qu'il existe un sous-ensemble $I\subset M\setminus\{0\}$ fini tel que
$$\sum_{i\in I} A_iA_{-i}=A_0.$$
[\emph{Indication :} Dans le cas contraire, considérer un idéal maximal $\mgo$ de $A_0$ contenant l'idéal $J=\sum_{i\in M\setminus\{0\}} A_iA_{-i}$. On note $s\in S$ le point fermé correspondant à $\mgo$. Montrer que la fibre $P_s$ est stable sous $G$, et que l'action induite de $G$ sur $P_s$ est encore libre. Appliquer alors le point 1. et obtenir une contradiction.]

\noindent {\bf 3.} On note $N=\{i\in M\ |\ A_iA_{-i}=A_0\}$. Montrer que $N$ est un sous-groupe de $M$.

\noindent {\bf 4.} Montrer le résultat escompté lorsque $A_0$ est un corps. [\emph{Indication :} On raisonne par l'absurde. Si $M\neq N$, alors $M/N\neq 0$. On regarde $G'=D_S(M/N)$. C'est un sous-groupe fermé de $G$. Il agit donc librement sur $P$ (via $G$). Appliquer le résultat de la question 2. à cette action et obtenir une contradiction.]

\noindent {\bf 5.} Pour $i\in M$, montrer que $i\in N$ si et seulement si pour tout idéal maximal $\mgo$ de $A_0$, on a $A_0=\mgo+A_iA_{-i}$.

\noindent {\bf 6.} Pour un point fermé $s\in S$, correspondant à un idéal maximal $\mgo$ de $A_0$, montrer que l'on a $A_0=\mgo+A_iA_{-i}$ pour tout $i\in M$. [\emph{Indication :} Considérer comme dans la question~2. l'action induite de $G$ sur la fibre $P_s$ et utiliser la question 4.]

\noindent {\bf 7.} Conclure.
\end{sousexo}

\begin{sousexemple}\rm
 Le théorème~\ref{quotient_par_gpe_diagonalisable} montre que le groupe $PGL_{n,S}$ est affine sur $S$. En effet, par définition, ce groupe est le quotient de $GL_{n,S}$ (qui est affine) par l'action libre de ${\gm}_{,S}$ par homothéties. Plus précisément, si $S=\Spec A$, alors $PGL_{n,S}$ est le spectre du sous-anneau de $A[x_{11},\dots, x_{nn},\Delta^{-1}]$ formé des éléments homogènes de degré 0.
\end{sousexemple}

\begin{sousexemple}\rm
 Soit $k$ un corps. On fait agir ${\gm}_{,k}$ sur $\A^2_k\setminus\{0\}$ par $(x,y) \mapsto (\lambda x, \lambda y)$. Le quotient est représentable et s'identifie à $\P^1_k$, qui n'est pas affine. On retrouve donc grâce au théorème ci-dessus le fait que $\A^2_k\setminus\{0\}$ n'est pas affine. On peut aussi utiliser le théorème~\ref{quotient_par_gpe_diagonalisable} pour vérifier que le quotient est bien $\P^1_k$. En effet, $\A^2_k\setminus\{0\}$ est la réunion de deux ouverts affines stables $D(x)$ et $D(y)$. D'après le théorème, le quotient de $D(x)$ par l'action donnée de $\gm$ est représentable par le sous-anneau de $k[x,y,x^{-1}]$ formé des éléments invariants sous $\gm$, c'est-à-dire $k[\frac{y}x]$. De même le quotient de $D(y)$ s'identifie au spectre de $k[\frac{x}y]$, si bien que le quotient de $\A^2_k\setminus\{0\}$ s'obtient en recollant deux copies de la droite affine $\Spec k[t]$ par $t\mapsto t^{-1}$.
\end{sousexemple}

\begin{sousremarque}\rm
 Le cas particulier où $P$ est lui-même un groupe \emph{diagonalisable} et \emph{de type fini} sur $S$ dont $G$ est un sous-groupe mérite une digression. Dans ce cas, on peut écrire $P=D_S(N)$ où $N$ est un groupe abélien de type fini (\emph{cf.} \ref{prop_gpes_diagonalisables}). Les sorites de SGA~3 sur les groupes diagonalisables (voir~\cite{SGA3}~VIII~1.5, 3.1 et~3.2) montrent que le morphisme de $S$-groupes $G\to P$ provient d'un morphisme de groupes ordinaires $u : N\to M$, qui est même surjectif, et que $G$ est en fait automatiquement un sous-groupe \emph{fermé} de $P$. Alors, en notant $M'$ le noyau de $u$, le faisceau quotient $P/G$ est représentable par $D_S(M')$. 
\end{sousremarque}

\begin{sousexemple}\rm
 Le quotient de $\gm$ par son sous-groupe~$\mu_n$ est représentable par un schéma affine. On le savait déjà : le quotient n'est autre que $\gm$ lui-même d'après la suite exacte de Kummer. La suite exacte de Kummer correspond à la suite exacte de groupes abéliens
$$\xymatrix{0 \ar[r] & \Z \ar[r]^{.n} & \Z \ar[r]& \Z/n\Z \ar[r]& 0.}$$
\end{sousexemple}

\begin{sousremarque}\rm
Dire que pour tout $i\in M$, on a $A_iA_{-i}=A_0$, est en fait équivalent à dire que le (mono)morphisme $P\times_S G \to P\times_S P$ est une immersion fermée\footnote{En effet, c'est équivalent à l'égalité $A_iA=A$ pour tout $i$ dans $M$, ce qui équivaut encore à la surjectivité du morphisme $A\otimes_{\Oc_S}A \to A[M]$, $a\otimes b \mapsto \sum a_i bX^i$.}. En particulier, si l'on sait \emph{a priori} que ce morphisme est une immersion fermée, le théorème~\ref{quotient_par_gpe_diagonalisable} est plus ou moins trivial. Mais surtout, on a le
\end{sousremarque}

\begin{souscor}\ 

 \begin{itemize}
\item[a)] Sous les hypothèses de~\ref{quotient_par_gpe_diagonalisable}, le morphisme
$$\Psi : P\times_S G \lto P\times_S P$$ est une immersion fermée. 
\item[b)] Pour toute section $\sigma : S\to P$ de $P$, le morphisme d'orbite $G \to P$ donné fonctoriellement par $g \mapsto \sigma.g$, est une immersion fermée. (En particulier, dans un groupe affine, tout sous-groupe diagonalisable est automatiquement fermé.)
\item[c)] Si $M$ est de type fini, le quotient $X=P/G$ est un \og quotient géométrique universel\fg\ au sens de Mumford (\emph{cf.}~\cite{Mumford_GIT}). 
 \end{itemize}
\end{souscor}
\begin{demo}
 Le point a) a déjà été démontré et b) en est une conséquence immédiate puisque le morphisme d'orbite se déduit du morphisme de a) par le changement de base $P\to P\times_S P,\ p\mapsto (\sigma, p)$. Le cas particulier donné entre parenthèses s'obtient en appliquant b) à la section neutre du groupe affine considéré. Montrons c). Avec la terminologie de \emph{loc. cit.} on sait déjà que $X$ est un \og quotient catégorique universel\fg\ puisqu'il représente le faisceau quotient.
 Donc d'après la remarque (3) de \cite{Mumford_GIT}~0.~\S 2, pour conclure il suffit de montrer que $P\to X$ est universellement submersif\footnote{Un morphisme $\pi : P\to X$ est dit submersif si pour tout sous-ensemble $U\subset X$, on a $U$ ouvert dans $X$ si et seulement si $\pi^{-1}(U)$ est ouvert dans $P$.} et que l'image du morphisme $\Psi$ de a) est exactement $P\times_X P$. Or, il est évident que $P\to X$ est universellement submersif puisqu'il est fidèlement plat et de présentation finie (donc surjectif et universellement ouvert). D'autre part, le morphisme $\Psi$ se factorise en
$$P\times_S G \simeq P\times_X G_X \lto P\times_XP \lto P\times_S P$$
où la première flèche est un isomorphisme puisque $P$ est un $G_X$-torseur sur $X$, ce qui prouve c). D'ailleurs, on retrouve aussi a) puisque la seconde flèche de cette factorisation est une immersion fermée (car $X$ est affine donc séparé sur $S$). 
\end{demo}

\subsection{Passage au quotient par des actions non libres}
\label{groupoides}

Nous allons conclure avec quelques mots sur les actions non libres. Dans ce cas la situation est plus compliquée. Il y a tout de même quelques résultats positifs. Tout d'abord, les théorèmes de SGA 3 donnés plus haut pour l'existence d'un conoyau admettent une version \og groupoïde\fg, certes aux conclusions moins satisfaisantes. Ils donnent seulement l'existence d'un quotient catégorique. Par exemple on a 

\begin{sousthm}[SGA 3, exposé~V]
\label{thm_conoyau_gpd}
 Soient $S$ un schéma localement noethérien et $X_*$ un $S$-groupoïde (voir~\ref{def_groupoides} pour les notations) tel que :
\begin{itemize}
 \item les flèches $s$ et $b$ soient propres et plates ;
\item $X_0$ soit quasi-projectif sur $S$ ;
\item le morphisme $(s,b) : X_1\to X_0\times_S X_0$ soit quasi-fini.
\end{itemize}
Alors le couple de flèches $(s,b)$ admet un conoyau $\pi : X_0 \to Y$ dans la catégorie des schémas. De plus, ce conoyau est un conoyau dans la catégorie des espaces annelés. Le morphisme $\pi$ est surjectif, ouvert et propre, et les morphismes $\pi$ et $Y\to S$ sont de présentation finie. Enfin, le morphisme $(s,b) : X_1\to X_0\times_Y X_0$  est surjectif.
\end{sousthm}

Ce théorème ne donne pas la représentabilité du faisceau quotient. De fait, dès qu'il y a de l'inertie, le faisceau quotient a tendance à ne pas être représentable.

\begin{sousexemple}\rm
\label{exemple_droite_affine}
 Considérons par exemple la droite affine $X=\A_S^1$ sur $S=\Spec \Z$, sur laquelle agit le groupe $G=\Z/2\Z$ par $x\mapsto -x$. Dans ce cas particulier, on peut aisément calculer le conoyau schématique dont le théorème précédent donne l'existence. Il s'agit du morphisme $\pi : \A_S^1 \to \A_S^1$ donné par $x\mapsto x^2$. Soit $X/G$ le faisceau quotient \emph{fppf}. C'est le faisceau associé au préfaisceau $P$ qui à un schéma $U$ associe $\Gamma(U,\Oc_U)/\{\pm 1\}$. Le morphisme $\pi$ se factorise en $g\circ f$ où $f : \A_S^1 \to P$ est le morphisme de passage au préfaisceau quotient et $g$ est défini fonctoriellement par $g(x)=x^2$. Supposons que $X/G$ soit représentable. C'est alors automatiquement un quotient dans la catégorie des schémas (par propriété universelle) donc il doit être égal à $\pi$. On en déduit que $g$ vérifie la propriété universelle du faisceau associé. Soit $U$ le spectre de $\Z[\eps]/(\eps^2)$. Alors $g(\eps)=0$, donc il existe une famille couvrante \emph{fppf}, que l'on peut supposer réduite à un élément $V\to U$ avec $V$ affine, telle que $\eps_{|_V}$ soit nul dans $\Gamma(V,\Oc_V)/\{\pm 1\}$, donc dans $\Gamma(V,\Oc_V)$. C'est impossible car $V\to U$ est fidèlement plat donc $\Gamma(U, \Oc_U)$ s'injecte dans $\Gamma(V,\Oc_V)$, ce qui prouve que le faisceau quotient n'est pas représentable.
\end{sousexemple}

\begin{sousexo}\rm
\label{exemple_droite_affine_esp_alg}
 Montrer que le quotient catégorique $\pi : X\to Y$ de l'exemple précédent est aussi un quotient dans la catégorie des espaces algébriques. En déduire que le faisceau quotient $X/G$ n'est pas représentable par un espace algébrique.
\end{sousexo}

\begin{sousexo}\rm
 On fait toujours agir $G=\Z/2\Z$ sur la droite affine $\A_S^1$ comme ci-dessus, mais au-dessus de la base $S=\Spec \F_2$. Montrer que l'action est triviale et que le conoyau donné par le théorème~\ref{thm_conoyau_gpd} n'est autre que l'identité de $\A_S^1$. En déduire que la formation du quotient catégorique ne commute pas au changement de base.
\end{sousexo}

Dans le cas particulier d'une action de groupe, l'énoncé~\ref{thm_conoyau_gpd} devient le suivant. Soient $S$ un schéma localement noethérien, $X$ un $S$-schéma quasi-projectif et $G$ un schéma en groupes agissant sur $X$. On suppose que la projection $p_2 : G\times_S X \to X$ est propre et plate, et que le morphisme $(\mu,p_2) : G\times_S X \to X\times_S X$ est quasi-fini. Alors il existe un quotient catégorique $\pi : X\to Q$ dans la catégorie des schémas. Sous l'hypothèse que $X$ est fidèlement plat et quasi-compact sur $S$ (par exemple si $S$ est le spectre d'un corps), on voit que le théorème~\ref{thm_conoyau_gpd} ne s'applique que lorsque le groupe $G$ est \emph{propre et plat} sur~$S$. Si l'hypothèse de platitude semble difficile à éviter\footnote{Si la dimension des fibres de $G$ varie, le théorème de semi-continuité pour la dimension des fibres montre qu'en général il ne peut exister de quotient géométrique.}, il est fâcheux de devoir supposer~$G$ \emph{propre}. Il y a bien sûr dans SGA 3 un analogue de~\ref{cas_plat}, où l'on suppose seulement $G$ plat et de type fini et où l'on obtient l'existence d'un quotient catégorique \emph{génériquement}. Mais en pratique, on a souvent besoin de construire le quotient d'une variété sous l'action d'un groupe affine, et ceci globalement. Par exemple, lorsque l'on cherche à construire un espace de modules pour un certain type d'objets algébriques, disons certaines variétés, on procède souvent de la manière suivante. On commence par construire un espace de modules pour nos variétés munies d'un plongement dans un espace projectif~$\P_n$ fixé, puis on essaye de se débarrasser du plongement. Il faut alors généralement quotienter par le groupe d'automorphismes de~$\P_n$.  

Si l'action est \emph{libre}, le théorème de représentabilité d'Artin vu plus haut donne une réponse tout à fait satisfaisante : le faisceau quotient $X/G$ est un espace algébrique. Dans le cas général, Keel et Mori ont démontré le théorème suivant.

\begin{sousthm}[\cite{Keel_Mori}]
 Soit $G$ un $S$-schéma en groupes plat, séparé et de présentation finie, qui agit sur un espace algébrique $X$ quasi-séparé et localement de présentation finie. On suppose que \og les stabilisateurs sont finis\fg, \emph{i.e.} que le morphisme $$j^{-1}(\Delta_X) \lto X$$ est fini, où $j$ est le morphisme de $G\times_S X$ vers $X\times_S X$ défini fonctoriellement par $j(g, x)=(g.x, x)$. Alors il existe un $S$-espace algébrique $Q$ et un morphisme $\pi : X\to Q$ qui fait de $Q$ un quotient catégorique de $X$ par $G$ (dans la catégorie des espaces algébriques). De plus :
\begin{itemize}
 \item[a)] Pour tout point géométrique $\xi$, $\pi$ induit un isomorphisme $X(\xi)/G(\xi) \simeq Q(\xi)$. 
\item[b)] Pour tout morphisme $Q'\to Q$ \emph{plat}, $Q'$ est un conoyau du couple de flèches $\xymatrix@C=1pc{R'\double &X'}$ où $R'=(G\times_S X)\times_Q Q'$ et $X'=X\times_Q Q'$. (Ceci implique en particulier que le faisceau étale $\Oc_Q$ est formé des fonctions $G$-invariantes sur $X$.)
\item[c)] $\pi$ est surjectif et universellement ouvert.
\end{itemize}
Enfin, si l'on suppose de plus que l'action est propre, \emph{i.e.} que le morphisme $j$ ci-dessus est propre, alors $Q$ est un espace algébrique \emph{séparé}.
\end{sousthm}

\begin{sousexemple}\rm
 Dans le cas de l'exemple~\ref{exemple_droite_affine}, le morphisme $\pi$ est aussi le quotient au sens de Keel et Mori (d'après l'exercice~\ref{exemple_droite_affine_esp_alg}). 
\end{sousexemple}

\begin{sousexemple}\rm
Mumford construit dans~\cite[chap.~V]{Mumford_GIT} un espace de modules grossier~$M_g$ pour les courbes de genre $g$ ($g\geq 2$) comme quotient géométrique d'un certain schéma quasi-projectif par une action naturelle de $PGL(n+1)$. Voici une esquisse très rapide de la construction, nous renvoyons au texte original pour plus de détails. Un entier naturel $\nu\geq 3$ étant fixé, Mumford définit un sous-schéma $H_{\nu}$ du schéma de Hilbert~$Hilb_{\P_n}^{P(x)}$ où
$$P(x)=(2x\nu-1).(g-1)\quad \textrm{ et } \quad n=P(1)-1.$$
Heuristiquement $H_{\nu}$ est le sous-schéma des courbes avec plongement $\nu$-canonique dans~$\P_n$. Le groupe $PGL(n+1)$ agit naturellement sur $H_{\nu}$. Mumford montre alors dans le \S 2 du chapitre~V que, s'il existe un quotient géométrique de $H_{\nu}$ par $PGL(n+1)$, c'est un espace de modules grossier pour les courbes de genre~$g$, \emph{cf.} prop. 5.4 (\emph{i.e.} il est universel pour les morphismes vers un schéma, et il a \og les bons\fg\ points géométriques). Puis il prouve (\S 3 et \S 4) l'existence d'un quotient géométrique à l'aide des techniques ``GIT'' développées dans les chapitres précédents. Le théorème de Keel et Mori ci-dessus permet d'obtenir directement l'existence d'un quotient géométrique en tant qu'espace algébrique. 
\end{sousexemple}

Enfin, il convient de signaler qu'il existe un autre objet quotient qui peut rendre de nombreux services. Il s'agit du champ quotient. Il a le mérite d'exister bien plus souvent (en tant que champ algébrique) que les quotients mentionnés ci-dessus. Le quotient de Keel et Mori, lorsqu'il existe, est alors un espace de modules grossier pour ce champ algébrique.

\begin{sousdefi}
 Soient $S$ un schéma, $X$ un $S$-schéma, et $G$ un $S$-schéma en groupes qui agit sur $X$ à droite. On note $[X/G]$ (resp. $[X/G]_{\text{pl}}$) la catégorie fibrée en groupoïdes suivante sur $(Sch/S)$. Pour tout $U$, $[X/G](U)$ est la catégorie des couples $(P,\alpha)$ où $P$ est un $G_U$ torseur étale (resp. \emph{fppf}) sur $U$ et $\alpha : P\to X\times_S U$ est un $U$-morphisme qui est $G_U$-équivariant.
\end{sousdefi}

\begin{sousthm}
 Soient $S$ un schéma, $X$ un schéma quasi-séparé et $G$ un $S$-schéma en groupes séparé, plat et de présentation finie qui agit à droite sur $X$. Alors le $S$-champ $[X/G]_{\text{pl}}$ est un $S$-champ algébrique.
\end{sousthm}

\begin{sousremarque}\rm
 Si $G$ est lisse alors $[X/G]$ et $[X/G]_{\text{pl}}$ coïncident.
\end{sousremarque}

\bibliographystyle{../../../texmf/tex/latex/monlatex/plain-fr}
\addcontentsline{toc}{section}{Bibliographie}
\bibliography{../../../texmf/tex/latex/monlatex/mabiblio}
\end{document}